\documentclass[11pt]{article}
\usepackage{amsmath}
\usepackage{amssymb}
\usepackage{amscd}

\newtheorem{proposition}[subsubsection]{Proposition}
\newtheorem{definition}[subsubsection]{Definition}
\newtheorem{lemma}[subsubsection]{Lemma}
\newtheorem{corollary}[subsubsection]{Corollary}
\newtheorem{remark}[subsubsection]{Remark}

\newcommand{\Ext} {{\rm Ext}\,}
\newcommand{\EXT} {{{\cal E}\kern-.5pt{\it xt} }\,}
\newcommand{\Hom} {{\rm Hom}\,}
\newcommand{\HOM} {{{\cal H}\kern-.5pt{\it om} }\,}
\newcommand{\red} {{\rm red}\,}
\newcommand{\Aut} {{\rm Aut}\,}

\newcommand{\Bil} {{\rm Bil}\,}
\newcommand{\Biext} {{\rm Biext}\,}
\newcommand{\Div} {{\rm Div}\,}
\newcommand{\Ker} {{\rm Ker}\,}
\renewcommand{\Im} {{\rm Im}\,}

\newcommand{\Alb} {{\rm Alb}\,}
\newcommand{\Pic} {{\rm Pic}\,}
\newcommand{\Spec} {{\rm Spec}\,}
\newcommand{\Cospec} {{\rm Cospec}\,}
\newcommand{\Prim} {{\rm Prim}\,}
\newcommand{\cu}{{\rm cu}}
\newcommand{\rk}{{\rm rk}\,}
\newcommand{\F}{{\rm F}}
\newcommand{\V}{{\rm V}}
\newcommand{\Gr}{{\rm Gr}\,}
\newcommand{\Sym}{{\rm Sym}\,}
\renewcommand{\H}{{\rm H}}
\newcommand{\Lie}{{\bf Lie}\,}
\newcommand{\fppf}{{\rm fppf}}

\newcommand{\Inf}{{\rm Inf}}
\newcommand{\Aff}{{\bf A}}
\newcommand{\ZZ}{{\bf Z}}
\newcommand{\NN}{{\bf N}}
\newcommand{\QQ}{{\bf Q}}
\newcommand{\CC}{{\bf C}}
\newcommand{\DD}{{\bf D}}

\newcommand{\TT}{{\bf T}}
\newcommand{\WW}{{\bf W}}
\newcommand{\E}{{\mathbb E}}
\newcommand{\HH}{{\mathbb H}}
\newcommand{\GG}{{\bf G}}
\newcommand{\X}{{\mathbb X}}
\newcommand{\Y}{{\mathbb Y}}
\newcommand{\T}{{\mathbb T}}
\newcommand{\G}{{\mathbb G}}
\newcommand{\A}{{\mathbb A}}
\newcommand{\B}{{\mathbb B}}

\newcommand{\M}{{\mathbb M}}
\newcommand{\N}{{\mathbb N}}
\newcommand{\W}{{\mathbb W}}

\renewcommand{\P}{{\mathcal P}}
\newcommand{\Q}{{\mathcal Q}}
\renewcommand{\u}{{\mathrm u}}
\renewcommand{\v}{{\mathrm v}}
\newcommand{\Frac}{{\rm Frac}}
\newcommand{\crys}{{\rm crys}}
\newcommand{\logcrys}{{\rm logcrys}}
\newcommand{\Crys}{\mbox{$\rm Pic$}^{\rm crys, 0}}
\newcommand{\bCrys}{\mbox{$\bf Pic$}^{\rm crys, 0}}
\newcommand{\LogCrys}{\mbox{$\rm Pic$}^{\rm logcrys, 0}}
\newcommand{\bLogCrys}{\mbox{$\bf Pic$}^{\rm logcrys, 0}}
\newcommand{\bLogCrysno}{\mbox{$\bf Pic$}^{\rm logcrys}}
\newcommand{\barLogCrys}{\overline{\rm Pic}^{\rm logcrys, 0}}
\newcommand{\barbLogCrys}{\overline{\bPic}^{\rm logcrys, 0}}
\newcommand{\InfDef}{{\rm InfDef}\,}
\newcommand{\Tcrys}{{{\mathbb T}_\crys}\,}
\newcommand{\Gal}{{\rm Gal}\,}

\newcommand{\bPic}{\mbox{$\bf Pic$}} 
\newcommand{\bDiv}{\mbox{$\bf Div$}}

\renewcommand{\exp}{{\rm exp}\,}

\newcommand{\by}[1]{\stackrel{#1}{\rightarrow}}
\newcommand{\longby}[1]{\stackrel{#1}{\longrightarrow}}

\newcommand{\tensor}{\otimes}
\newcommand{\into}{\hookrightarrow}

\renewcommand{\hat}{\widehat}

\renewcommand{\leadsto}{\mbox{ $\longmapsto$ }}
\newcommand{\df}{\mbox{\,$:=$}\,}
\newcommand{\ie}{{\it i.e.\/},\ }
\newcommand{\cf}{{\it cf.\/}\ }
\newcommand{\eg}{{\it e.g.\/},\  }

\newcommand{\et}{\mbox{\scriptsize{\'{e}t}}}

\renewcommand{\bar}{\overline}
\newcommand{\sext}{\mbox{${\cal E}xt\,$}}  
\newcommand{\shom}{\mbox{${\cal H}om\,$}}  

\newcommand{\veq}{\mbox{\large $\parallel$}}  
\newcommand{\longto}{\longrightarrow}
\newcommand{\isomarrow}{\stackrel{\sim}{\longrightarrow}}
\newcommand{\llongrightarrow}{\mbox{{\Large{$\longto$}}}}

\newcommand{\liminv}[1]{{\displaystyle{\mathop{\rm
lim}_{\buildrel\longleftarrow\over{#1}}}}\,}

\newcommand{\boxtensor}{{\Box\kern-9.03pt\raise1.42pt\hbox{$\times$}}}

\renewcommand{\d}{\mbox{\LARGE $\cdot $}}

\newcommand{\arc}[1]{{#1}\kern-16pt\raise7pt\hbox{\Large{$\frown$}}}

\newcommand{\longinto}{\lhook\joinrel\kern-3pt\hbox to
100pt{\rightarrowfill}}           

\newcommand{\propsubset}{\mbox{$\textstyle{
\subseteq_{\kern-5pt\raise-1pt\hbox{\mbox{\tiny{$/$}}}}}$}}

\newcounter{elno}                

\newenvironment{proof}{{\bf Proof}:\quad
                         }{\hfill$\odot$\par\vspace*{0.5cm}}

\newcommand{\cA}{{\cal A}}
\newcommand{\cB}{{\cal B}}
\newcommand{\cC}{{\cal C}}
\newcommand{\cD}{{\cal D}}
\newcommand{\cE}{{\cal E}}
\newcommand{\cF}{{\cal F}}
\newcommand{\cG}{{\cal G}}

\newcommand{\cI}{{\cal I}}

\newcommand{\cK}{{\cal K}}
\newcommand{\ccL}{{\cal L}}
\newcommand{\cM}{{\cal M}}

\newcommand{\cO}{{\cal O}}
\newcommand{\cP}{{\cal P}}

\begin{document}

\title{\bf Crystalline realizations of 1-motives}

\author{Fabrizio {\sc Andreatta} and Luca {\sc Barbieri Viale}}

\date{}

\maketitle

\begin{abstract}
We consider the crystalline realization of Deligne's 1-motives in
positive characteristics and prove a comparison theorem with the
De Rham realization of (formal) liftings to zero characteristic.

We then show that one dimensional crystalline cohomology of an algebraic
variety, defined by forcing universal cohomological descent {\em via}\,
de Jong's alterations, coincides with the crystalline realization of the
Picard 1-motive, over perfect fields of cahracteristic $>2$.
\end{abstract}

\section*{}

Grothendieck's vision of crystals is well explained in \cite{MG} (see
also \cite{MazurMessing} and the subsequent papers by Berthelot and
Messing
\cite{BerMess}). Grothendieck pointed out that one can recover the first
De Rham cohomology of an abelian scheme in characteristic zero {\em
via}\, the Lie algebra of the universal $\GG_a$-extension of the dual.
Moreover, in positive characteristics, this universal $\GG_a$-extension
and the Poincar\'e biextension are crystalline in nature and depend only
on the $p$-divisible group associated to the abelian scheme. Recall the
following classical results.

If $G$ is any group scheme over a base scheme  $S$ such that
$\shom (G, \GG_a)=0$
and $\sext (G,
\GG_a)$ is a locally free $\cO_S$-module of finite rank, the universal
$\GG_a$-extension $\E
(G)$ exists: it is an extension of $G$ by the vector group $\sext (G,
\GG_a)^{\vee}$.\\

\noindent {\bf Theorem A} (\cite[Prop. I.4.1.4 \& \S\ II.1.5]
{MG}):\ {Let $A_0$ be an abelian scheme over a perfect field $k$
of characteristic $p>0$. Let $A$ be a (formal) lifting to the Witt
vectors $\WW (k)$ and let $A^{\vee}$ be the dual. Let
$\TT^{\crys}(A_0)$ be the contravariant Dieudonn\'e module
associated to the  Barsotti-Tate group $A_0[p^{\infty}]$. Then
$$\TT^{\crys}(A_0)\cong \Lie \E(A^{\vee}) \cong \H^1_{\rm DR}(A/\WW
(k))$$}

Providing a good definition of crystalline topology (\cf
\cite{BerthelotOgus})
one can recover one dimensional crystalline cohomology from the above.
We then also have:\\

\noindent {\bf Theorem B} (\cite[II.3.11.2]{Illusie}):\ {Let $X$ be
smooth and proper over a perfect field $k$ of characteristic $p>0$. Let
$\Pic^{0, \red}(X)$ be the abelian
Picard scheme. Let $\TT_{\crys}(-)$ denote the covariant Dieudonn\'e
module. Then
$$\TT_{\crys}(\Pic^{0, \red}(X))\cong \H^1_{\crys}(X/\WW (k))$$}

A geometric way to prove it is by applying the former to the
universal extension of a lifting, considering the Albanese variety
$\Alb (X) = \Pic^{0, \red}(X)^{\vee}$. Let $\Crys (X)$ be the
sheaf on the (small) fppf site on $\WW (k)$ given by the functor
associating to $T$ the group of isomorphism classes of crystals of
invertible  $\cO^\crys_{X\times_{\WW (k)} T/T}$-modules (which are
algebraically equivalent to~$0$ on the Zariski site on
$X\times_{\WW (k)} T$). By construction we have that $\Lie \Crys =
\H^1_{\crys}$, inducing a canonical isomorphism
$$\Lie\Crys (\Alb (X))\longby{\simeq} \Lie\Crys(X)$$
by the Albanese mapping. Furthermore, $\Lie\Crys (\Alb (X))$ can
be identified to the Lie algebra of the universal extension of a
(formal) lifting of $\Pic^{0, \red}(X)$ to the Witt vectors. Note
that $\Crys$ is the natural substitute of the usual functor
$\Pic^{\natural , 0}$ given by invertible sheaves with an
integrable connection, and recall that, in general, in positive
characteristics, $\H^1_{\rm DR}(X/k)$ cannot be recovered from the
Picard scheme, as illustrated by Oda \cite{Oda}.

Another application of Grothendieck's vision is provided by Deligne's
definition of the De Rham realization of a 1-motive (see \cite{Deligne}
and \cite{BS}). In fact, the universal $\GG_a$-extension $\E (\M)$ of a
1-motive $\M$ over any base scheme exists and the De Rham realization is
naturally defined {\it via}\, the Lie algebra (see
\cite[10.1.10-11]{Deligne}, \cf Section~2 below).

\subsection*{The results}
The purpose of this paper, after providing a 1-motive of the {\em
crystalline realization}, is to draw the picture above in the
motivic world, over a perfect field $k$.

The target category of our realizations is the category of filtered
$F$-crystals of level $1$ over $k$, for simplicity referred to as
filtered
$F$-crystals in the sequel. Actually, we consider the category of
filtered
$F$-$\WW(k)$-modules consisting of finitely generated
$\WW(k)$-modules endowed with an increasing filtration and a
$\sigma$-linear operator, the Frobenius $F$, respecting the
filtration (here $\sigma$ is the Frobenius on $\WW(k)$). Such category
is, in an obvious way, a tensor category. Filtered $F$-crystals are the
objects whose underlying $\WW(k)$-modules are free and there exists a
$\sigma^{-1}$-linear
operator, the Verschiebung $V$, such that $V \circ F=F \circ V=p$.
We let $\WW(k)(1)$ be the filtered $F$-crystal $\WW(k)$, with
filtration $W_n=\WW(k)$ if $n \geq -2$ and $W_n=0$ for $n<-2$ and
with the $\sigma$-linear operator $F$ given by $1 \mapsto 1$ and
the $\sigma^{-1}$-linear operator $V$ defined by $1\mapsto p$.

For a 1-motive $\M$ over a perfect field define $\M[p^n]$ as
$\H^{-1}\bigl(\M/p^n\M\bigr)$, where~$\M/p^n\M$ is the cone of
multiplication by~$p^n$ on~$\M$. Then get a Barsotti-Tate group
$\M[p^{\infty}]$ taking the direct limit. Passing to its
contravariant Dieudonn\'e module we get a filtered
$F$-crystal over $k$ which we denote $\TT^\crys (\M)$ (see
Section~1 for details). We note that such a realization is an
instance (the good reduction case) of a more general theory of
$p$-adic realizations of $1$-motives over a local field due to
J.-M.~Fontaine (to appear in \cite{FJ}).

On the other hand, following Deligne, define $\TT^{\rm DR} (\M)$
by $\Lie \E (\M^{\vee})$ where $\M^{\vee}$ is the Cartier dual
(see Section~2, \cf\ref{dual} for duals and \ref{vectorext} for
extensions). We show that for $p>2$ the universal extension and
the Poincar\'e biextension are crystalline (see Section~3).
Furthermore,  $\TT^{\rm DR} (-)$ yields a filtered crystal (no
restriction on $p$), providing the Dieudonn\'e crystal of the 1-motive
(in the
terminology of \cite{MazurMessing}, see~\ref{crysviaext}). We then
have the following comparison theorem.\\

\noindent {\bf Theorem A${}^{\prime}$}:\ {\em Let~$\M_0$ be a
$1$-motive over a perfect field~$k$ of positive
characteristic~$p>0$.  Choose a (formal) $1$-motive $\M$
over~$\WW(k)$ lifting~$\M_0$. Then, there is a canonical
isomorphism
$$\TT^\crys(\M_0)\cong \TT^{\rm DR} (\M)$$
of filtered $F$-crystals. Let $\TT_\crys (\M_0)$ be the crystal
given by $\TT^\crys(\M_0^\vee)$. There is a bilinear perfect
pairing of filtered $\WW(k)$-modules
$$\TT^\crys (\M_0) \otimes_{\WW (k)} \TT_\crys (\M_0) \to \WW (k)
(1).$$}

See \ref{vectorext} for the definition of the filtration on
$\TT^{\rm DR} (\M)$ and \ref{Mpinfty} for the definition of the
filtration on $\TT^\crys(\M_0)$.  The proof of
Theorem~A${}^{\prime}$ is given in Section~4 (see~\ref{propTcrys}
for a list of properties of $\TT_\crys$) and it `mostly' adorns
Grothendieck's
original proof for abelian schemes, \cf \cite{MazurMessing}. We
stress that the description of~$\TT^\crys(\M_0)$ in terms of
universal extensions is {\it essential} to prove Theorem~B${}^{\prime}$
below.\\

Now let $V$ be an algebraic variety over a perfect field $k$.
Applying de Jong's method \cite{DE} we obtain a pair $(X_{\d},
Y_{\d})$ where $X_{\d}$ is a smooth proper simplicial scheme,
$Y_{\d}$ is a normal crossing divisor in $X_{\d}$ and $V_{\d}\df
X_{\d}-Y_{\d}$ is a smooth hypercovering of $V$. Note that we may
and will assume that $V_0\to V$ is generically \'etale by
\cite[Thm. 4.1]{DE}. Let $\Pic^{+}(V)\df \Pic^{+}(X_{\d}, Y_{\d})$
be the Picard 1-motive of $V$ (see Section~5 and
Appendix~\ref{appendix}, \cf \cite{BS} and \cite{RA}). In the same
way, as de Jong suggested in \cite[p. 51-52]{DE}, we set
$\H^*_\crys (V_{\d}/\WW (k)) \df \HH^*_\logcrys (X_{\d}, Y_{\d})$
where $(X_{\d}, Y_{\d})$ here denotes the simplicial logarithmic
structure on $X_{\d}$ determined by $Y_{\d}$ (see Section~6). We
show in \ref{verygory} that $\HH^1_\logcrys (X_{\d}, Y_{\d})$ is
naturally a free filtered $F$-$\WW (k)$-module. We have
the following link.\\

\noindent {\bf Theorem B${}^{\prime}$}:\ {\em Let $V$ be an
algebraic variety over a perfect field $k$ of characteristic
$p\geq 3 $. There is a functorial  isomorphism of filtered
$F$-$\WW (k)$-modules
$$\TT_\crys (\Pic ^{+}(V))  \longby{\simeq}
\H^1_\crys (V_{\d}/\WW (k))(1).$$}

The proof of Theorem~B${}^{\prime}$ is the full Section~7 and, very 
roughly speaking, goes as follows.
By the comparison Theorem~A${}^{\prime}$, we are left to compute
the universal extension crystal of the $1$-motive $\Pic ^{+}(V)$.
Thus, one gets to deal with a simplicial version of the functor of
invertible crystals $\LogCrys$ showing that its Lie algebra
coincides with $\HH^1_\logcrys$ (see 7.3--7.5).

It is shown in the Appendix~\ref{appendix} that $\Pic ^{+}(V)$ is
independent of the choices made, \ie of the pair $(X_{\d},
Y_{\d})$, and in particular functorial in $V$. Thus {\em via}\, 
Theorem~B${}^{\prime}$, see Corollary~\ref{independenceH1crysV},  
$\H^1_\crys(V_{\d}/\WW(k))$
is independent of the chosen hypercovering $V_{\d}\to V$ so
that one can define $\H^1_\crys(V/\WW(k)):=
\H^1_\crys(V_{\d}/\WW(k))$ forcing descent for
crystalline cohomology. This answers a question raised in \cite[p.
52]{DE} concerning the independence of
$\H^i_\crys(V_{\d}/\WW(k))\otimes \QQ$ from the choice of
hypercoverings, at least for $i=1$ and $p\geq 3$.

For the sake of exposition we omit the case $p=2$ deserving an {\it ad
hoc}\, explanation. Actually, in this case, the independence of
$\Pic ^{+}(V)$ from the choices made and $\H^1_\crys(V_{\d}/\WW(k))$ are still in place. Some technical difficulties in proving Theorem~B${}^{\prime}$ arise from the fact that the standard divided power structures on the ideal $2\WW(k)$ are {\em not} topologically nilpotent; see \ref{casep=2}. We will treat these matters elsewhere.\\

Note that the abelian quotient of $\Pic ^{+}(V)$ is the (reduced)
identity component of the kernel of the canonical map $\Pic^{0,
\red}(X_0)\to \Pic^{0, \red}(X_1).$ Moreover, for $V$ a normal
proper scheme $\Pic ^{+}(V) = \Pic^{0, \red}(V)$ is abelian, \eg
for every prime $p$ we obtain from Theorem~B${}^{\prime}$
$$\H^1_\crys (V/\WW (k)) \cong \Ker \bigl(\H^1_\crys(X_0)\to
\H^1_\crys(X_1)\bigr)$$ if $V$ is a normal projective variety.

Recall that the Cartier dual of $\Pic ^{+}(V)$ is the
(homological) Albanese 1-motive $\Alb^{-}(V)$ see \cite{BS} (\cf
\cite{RA}), \eg for $V$ a normal proper scheme is the cokernel of
the map $\Alb (X_1)\to \Alb (X_0)$.\\

\noindent {\bf Corollary}:\ {\em There is a functorial $\WW
(k)$-linear isomorphism
$$\TT_\crys (\Alb^{-}(V)) \longby{\simeq} \Hom (\H^1_\crys (V/\WW
(k))(1),
\WW (k)(1)).$$}

This suggests that one can define the first crystalline homology
group (modulo torsion!) as $\H_1^\crys (V/\WW (k))\df \Hom
(\H^1_\crys (V/\WW (k))(1), \WW (k)(1))$. Finally note that
Theorem B${}^{\prime}$ can be reformulated in terms of the first
rigid cohomology group of~$V$, \ie
$$\TT_\crys (\Pic ^{+}(V)) \otimes \QQ \longby{\simeq}
\H^1_\crys (V/\WW (k))(1)\otimes \QQ\longby{\simeq}\H^1_{\rm rig}
(V/\Frac(\WW (k)))(1).$$The last isomorphism follows using the
recent proof of cohomological descent of rigid cohomology for
proper coverings due to N.~Tsuzuki~\cite{Tsu} and the comparison
between log-crystalline cohomology and rigid cohomology (for
proper and smooth schemes with strict normal crossing divisors)
due to A. Shiho~\cite{Shi}. This suggests, as remarked by
B.~Chiarellotto and K.~Joshi, that one could approach Theorem
B${}^{\prime}$ using rigid cohomology. On the other hand, our
approach allows to prove the version of Theorem B${}^{\prime}$
with integral coefficients.

In the same spirit, the Corollary can be rephrased for rigid
homology. However, if we compare to \cite{Pet}, we find that
$\H_*^\crys (V/\Frac (\WW (k)))\neq \H_*^{\rm rig}(V/\Frac(\WW
(k)))$ in general, since the rigid homology defined in \cite[\S
2]{Pet} is the Borel-Moore variation of crystalline homology, \eg
they coincide for proper schemes.

\subsection*{Open problems}
A natural question is to extend our crystalline realization
functor to mixed motives. Recall that the category of 1-motives
(up to isogeny) is embedded in Voevodsky's triangulated category
of motives. However, for 1-motives over a perfect field, the more
striking question is to prove the following claim.

Let $(X, Y)$ be a pair such that $Y$ is a closed subvariety of an
algebraic variety $X$ over a perfect field. In \cite{BM}, over an
algebraically closed field of characteristic zero, effective
1-motives with torsion $\Pic^{+}(X,Y; i)$ (= $M_{i+1} (X, Y)$ in
\cite{BM}) are constructed showing Deligne's conjecture
\cite[10.4.1]{Deligne}, \eg showing the algebraicity of the
maximal mixed Hodge structure of type $\{(0,0), (0,-1), (-1,0),
(-1,-1)\}$ contained in $\H^{i+1}(X,Y; \QQ (1))$ if $k =\CC$. By
using an appropriate resolution it is easy to modify the
construction yielding $\Pic^{+}(X,Y; i)$ over a perfect field,
such that $\Pic^{+}(X,\emptyset ; 0) = \Pic^{+}(X)$ above
(however, it is not clear,  for $i>0$, if $\Pic^{+}(X,Y; i)$  is
{\it integrally} well defined!). Furthermore, it is easy to obtain
a Barsotti-Tate crystal of an effective 1-motive, using Fontaine's
theory \cite{FO}. Note that now we deal with finite groups and
therefore we can lift to the Witt vectors only the free part of an
effective 1-motive. However, we may define $\TT_\crys (-)$ by
means of the covariant Dieudonn\'e module, \ie given by the
Dieudonn\'e module of the
dual of the associated formal $p$-group.\\

\noindent {\bf Conjecture~C}:\ {\em Provide a weight filtration
$W_*$ on the crystalline cohomology $\H^{*}_{\crys} ((X, Y)/\WW
(k))$ of the pair $(X, Y)$. Let $\H^{*}_{\crys ,(1)} ((X, Y)/\WW
(k))$ denote the submodule of $W_2\H^{*}_{\crys} ((X, Y)/\WW (k))$
whose image in $\Gr_W^2$ is generated by the image of the discrete
part of $\Pic ^{+}(X,Y; i)$ under a suitable cycle map. Then there
is a canonical isomorphism (eventually up to $p$-power isogenies)
$$\TT_\crys (\Pic ^{+}(X,Y; i)) \longby{\simeq} \H^{i+1}_{\crys
,(1)} ((X, Y)/\WW (k))(1)$$ of filtered $F$-$\WW (k)$-modules.\\}

The corresponding statement for De Rham cohomology over a field of
characteristic zero is Theorem~3.5 in \cite{BM}.
Theorem~B${}^{\prime}$ above corresponds to this statement for $i
=0$ and $Y =\emptyset $. Note that only the free part of $\Pic
^{+}(X,Y; i)$ is independent of the hypercoverings in char.~$0$
(\cf \cite[2.5]{BM} and \cite[4.4.4]{BS}). Here $\H^{*}_{\crys}
((X, Y)/\WW (k))$ is defined following \cite{DE}. However, by
dealing with rigid cohomology Conjecture C, can be rephrased
switching crystalline to rigid. Note that a cycle map for rigid
cohomology is fully described in \cite{Pet}.

\subsubsection*{Acknowledgements}
We would like to thank B.~Chiarellotto and V.~Cristante for
several useful discussions. We also thank J.-M.~Fontaine and
K.~Joshi for informing us about the contents of their forthcoming
book \cite{FJ}. We finally thank the referee for several useful
comments which helped us to improve the exposition.

\subsection*{Notation}
We let our base schemes be locally noetherian. For a perfect field
$k$ we let $\WW (k)$ (resp. $\WW_n(k)$ with $n\in\NN$) denote the
ring of (truncated) Witt vectors with the standard divided power
structure on its maximal ideal. Note that for $p=2$ the standard
divided power structure on the maximal ideal of $\WW_n(k)$ is {\it
not} nilpotent. For $S_0$ a base scheme such that $p$ is locally
nilpotent we let $S_0\into S$ be a thickening defined by an ideal
$I$ with nilpotent divided powers.

Let $\G$ denote a group scheme over $S$. We consider $\G$ as a
sheaf for the fppf topology on $S$. We sometimes denote $g\in \G$
an $S$-point. Denote $\GG_a$ and $\GG_m$ the usual additive and
multiplicative structure $S$-group schemes.

If $\underline{G}$ is a $p$-divisible group over a perfect field $k$ of
characteristic $p$, we denote by $\DD(\underline{G})$ the {\it
contravariant
Dieudonn\'e module} of $\underline{G}$ defined as the module over the
Dieudonn\'e ring $D_k:=\WW(k)[F,V]/(FV=VF=p)$ of homomorphisms
from $\underline{G}$ to the group of Witt covectors over $k$; see
\cite[\S
III.1.2]{FO}. By \cite[\S III.6.1]{FO} such functor defines an
antiequivalence from the categroy of $p$-divisible groups over $k$
to the category of $D_k$-modules which are finite and free as
$\WW(k)$-modules. By the comparison Theorem \cite[Thm.
15.3]{MazurMessing} $\DD(\underline{G})$ coincides with the Lie algebra
of
the universal extension of a lifting of the dual $p$-divisible
group $\underline{G}^\vee$ to $\WW(k)$. For the latter approach to
contravariant
Dieudonn\'e theory we refer the reader to \cite[IV.2.4.3]{Messing} or 
\cite[\S 9.2 ]{MazurMessing}.

For a simplicial $S$-scheme $X_{\d}$ we denote $d_i^j$ the faces $X_j\to
X_{j-1}$, \eg $d_0^2, d_1^2, d_2^2\colon X_2\rightarrow X_1$, over $S$,
and sometimes we omit the $j$-index when it is clear from the context.

\section{Barsotti-Tate crystal of a 1-motive}
First recall some definitions and constructions. We refer to \cite[\S
10]{Deligne} and \cite[\S 1]{BS} for more details on 1-motives.

\subsection{Deligne's 1-motives}\label{1motives}
Recall that algebraic 1-motives are originally defined in
\cite[D\'efinition 10.1.2]{Deligne} over an algebraically closed field.
However, in \cite[Variante 10.1.10]{Deligne} (\cf \cite{BS}), the
definition of {\it 1-motif lisse}\, is taken over arbitrary
base schemes, as follows.  Let\/~$S$ be a scheme. A $1$-motive
$\M$ over~$S$, is

\begin{itemize}
\item[\rm a)] a group scheme~$\X$ over~$S$ which \'etale locally
on~$S$ is  constant, free and of finite type as $\ZZ$-module;

\item[\rm b)] a semi-abelian scheme $\G$ over~$S$, extension of an
abelian scheme~$\A$ over~$S$ by a torus~$\T$ over~$S$;

\item[\rm c)] a homomorphism of $S$-group schemes~$\u\colon\X
\rightarrow \G$.

\end{itemize}

As customary, we view the category of commutative $S$-group
schemes as a full sub-category of the derived category~${\rm
D}^{\rm b}_\fppf(S)$, identifying an $S$-group scheme with the
complex having its underlying fppf sheaf concentrated in
degree~$0$. Analogously, we identify a $1$-motive $\M$ with the
complex in ${\rm C}^{\rm b}_\fppf(S)$ $$\cdots \to 0 \rightarrow
\X\stackrel{\u}{\longrightarrow} \G\rightarrow 0\to\cdots$$ of
fppf sheaves over~$S$ concentrated in degree~$-1$ and~$0$. One
sees (\cf \cite{BS}) that
$$\Hom_{{\rm C}^{\rm b}_\fppf(S)}(\M_1 , \M_2) =\Hom_{{\rm D}^{\rm
b}_\fppf(S)}(\M_1 , \M_2).$$By \cite[Prop. 2.3.1]{Ray} the
morphisms of $1$-motives, from
$\left[\X_1\stackrel{\u_1}{\longrightarrow} \G_1\right]$ to
$\left[ \X_2\stackrel{\u_2}{\longrightarrow} \G_2\right]$, also as
objects of~${\rm D}^{\rm b}_\fppf(S)$, are morphisms of complexes,
\ie given by pairs of homomorphisms $\X_1 \rightarrow \X_2$ and
$\G_1 \rightarrow \G_2$ of $S$-group schemes making the following
diagram
\begin{displaymath}
\begin{CD}
\X_1  @>{\u_1}>>  \G_1 \\ @VVV  @VVV \\
\X_2 @>{\u_2}>> \G_2
\end{CD}
\end{displaymath}
commute. Given a $1$-motive $\M$ we define a weight filtration
$W_0(\M):=\M$, $W_{-1}(\M):=[0\rightarrow \G]$,
$W_{-2}(\M):=[0\rightarrow \T]$ and $W_{n<-2}(\M)=0$.

\subsection{Cartier duality}\label{dual}
Cartier duality is naturally extended to 1-motives over a field
(see \cite[10.2.11]{Deligne} and \cf \cite[\S 1.5]{BS}). Over a
base scheme it is further extended as follows. Let $S$ be a
locally noetherian base scheme. Let~$\A$ be an abelian scheme
over~$S$. Let~$\A^\vee:=\Pic^0_{\A/S}$ be the abelian scheme dual
to~$\A$. Let~$\P \rightarrow \A\times_S \A^\vee$ be the Poincar\'e
$\GG_m$-bundle rigidified along~$\{0\}\times_S \A^\vee$
and~$\A\times_S\{0\}$: it defines a biextension of~$\A$
and~$\A^\vee$ by~$\GG_{m,S}$ \cite[VIII.3.2]{SGA7}. For every
$S$-scheme $T$, an extension of~$\A\times_S T$ by~$\GG_{m,T}$
defines a $\GG_{m,T}$-bundle over~$\A\times_S T$ and, hence, a
$T$-valued point of~$\A^\vee$. The map associating to~$x\in
\A^\vee(T)$ the $\GG_{m,T}$-bundle $\P\vert_{A \times_T \{y\}}$
rigidified over~$0\in \A(T)$ defines a homomorphism
form~$\A^\vee(T)$ to~$\Ext^1_T\bigl(\A\times_S T,\GG_{m,T}\bigr)$.
These two maps are inverse one of the other and define an
isomorphism $\A^\vee \isomarrow \Ext^1\bigl(\A,\GG_{m,S}\bigr)$.

\begin{itemize}
\item[\rm 1)] let $\T$ be a torus over~$S$ and let $\Y$ be its
character group. To give an extension~$\G$ of~$\A$ by~$\T$ is
equivalent to give a homomorphism~$\u_{\A^\vee}\colon \Y
\rightarrow \A^\vee$ over~$S$.

\item[] If~$y \in\Y$, then the image of~$-y$ (note the minus sign!)
in~$\A^\vee$, identified with~$\Ext^1(\A,\GG_{m,S})$, is the
unique rigidified extension of~$\A$ by~$\GG_{m,S}$ obtained as the
push-out of $$0 \llongrightarrow \T \llongrightarrow \G
\llongrightarrow \A \llongrightarrow 0$$ by~$y$;

\item[\rm 2)] to give a homomorphism $\u\colon\X \rightarrow
\G$ is equivalent to give  \item a homomorphism $\u_{\A}\colon \X
\rightarrow \A$
\item a trivialization, as  biextension, of
the pull back of~$\P$ by the homomorphism $$\u_{\A} \times
\u_{\A^\vee}\colon \X\times_S \Y \llongrightarrow \A \times_S
\A^\vee.$$ \end{itemize}

\noindent In particular, let
$\M:=\bigl[\X\stackrel{\u}{\longrightarrow} \G\bigr]$ be a
$1$-motive. Define the dual $1$-motive
$$\M^\vee:=\bigl[\X^\vee\stackrel{\u^\vee}{\longrightarrow}
\G^\vee\bigr]$$as follows:
\begin{itemize}
\item the group scheme~$\X^\vee$ is the character group~$\Y$
of~$\T$;
\item the semiabelian scheme~$\G^\vee$ is the
extension of~$\A^\vee$ by the torus~$\X^\vee\tensor_\ZZ \GG_{m,S}$
defined by the composite homomorphism
$\u_\A\colon\X\stackrel{\u}{\rightarrow}\G\rightarrow
\A=\left(\A^\vee \right)^\vee$;
\item the homomorphism~$\u^\vee$ is defined by
the trivialization of the biextension
$\bigl(\u_\A\times\u_{\A^\vee}\bigr)^*\left(\P\right)$
(\cf \cite[10.2.11]{Deligne}).
\end{itemize}

\subsection{The functors $\TT^{\crys}$ and $\TT_{\crys}$}\label{Mpinfty}
Let~$p$ be a prime number. Let~$n\in\NN$. Define the group
scheme~$\M[p^n]$ as $\H^{-1}\bigl(\M/p^n\M\bigr)$,
where~$\M/p^n\M$ is the cone of multiplication by~$p^n$ on~$\M$.
More explicitly,
$$\M[p^n]:=\frac{\Ker\Bigl( \u + p^n\colon \X \times_S \G
\llongrightarrow \G \Bigr)}{ {\rm Im}\Bigl( (p^n, -\u)
\colon\X\llongrightarrow \X\times_S \G\Bigr)}.$$
It is a finite and flat group scheme over~$S$ and it sits in the
exact sequence

\begin{equation}\label{exactseq}
0\llongrightarrow \G[p^n] \llongrightarrow \M[p^n]
\llongrightarrow \Bigl(\X/p^n\X\Bigr)\llongrightarrow
0.
\end{equation} Note also that we have the exact sequence $$
0\llongrightarrow \T[p^n] \llongrightarrow \G[p^n]
\llongrightarrow \A[p^n] \llongrightarrow 0.$$Define
$$\M[p^\infty]:=\lim_{n \rightarrow\infty} \M[p^n];$$the direct
limit being taken using the natural inclusions~$\M[p^n] \subset
\M[p^m]$ for~$m\geq n$. Then,~$\M[p^\infty] $ defines a
Barsotti-Tate group in the sense of~\cite[Def I.2.1]{Messing}.
It sits in the exact sequence

\begin{equation}\label{exactseqinfty}
0\llongrightarrow \G[p^\infty] \llongrightarrow \M[p^\infty]
\llongrightarrow \X[p^\infty]\llongrightarrow 0,
\end{equation}
\noindent where~$\X[p^\infty]:=\X\tensor_\ZZ
\Bigl(\QQ_p/\ZZ_p\Bigr)$. We also get the exact sequence
\begin{equation}\label{exactseqinfty2} 0\llongrightarrow
\T[p^\infty] \llongrightarrow \G[p^\infty] \llongrightarrow
\A[p^\infty] \llongrightarrow 0.
\end{equation}

Note that we clearly get a splitting of~(\ref{exactseqinfty}) over a
suitable faithfully flat extension of $S$. In fact, choose a faithfully
flat extension~$S'\rightarrow S$ and a compatible set of homomorphisms
$$\bigl\{\u_n\colon \frac{1}{p^n}\Bigl(\X\times_S S'\Bigr)
\rightarrow \G\times_S S'\bigr\}_{n\in\NN}$$such that
$\u_0=\u\times_S S'$. Such a choice determines a splitting
of~(\ref{exactseq}) over~$S'$, for any~$n$, defining
$\Bigl(\X/p^n\X\Bigr)\times_S S' \rightarrow \M[p^n]\times_S S'$
by~$\bar{x}\mapsto \bigl(x,-\u_n(p^{-n}x)\bigr)$. Hence, we get a
splitting of~(\ref{exactseqinfty}) over~$S'$.\\

Let\/~$S$ be a scheme where~$p$ is locally nilpotent. Denote by
\begin{equation}\label{dieu}
\DD\left(\M[p^\infty] \right)
\end{equation} the contravariant Dieudonn\'e crystal, on the crystalline
site of\/~$S$, associated to the Barsotti-Tate
group~$\M[p^\infty]$ (\cf \cite[IV.2.4.3]{Messing}).

Let\/~$k$ be a perfect field of characteristic~$p$ and
let\/~$\WW(k)$ be the Witt vectors of\/~$k$. Suppose that $\M$ is
defined over $k$.
\begin{definition}\label{defTcrysupanddwn} We
call\ {\em crystalline realizations}\, of $\M$ the following
$\WW(k)$--modules
$$\TT^{\crys}\left(\M\right)\df
\lim_{\infty\leftarrow
n}\DD\left(\M[p^\infty]\right)\Bigl(\Spec(k) \hookrightarrow
\Spec\left(\WW_n(k)\right)\Bigr)$$and
$$\TT_{\crys}\left(\M\right)\df
\lim_{\infty\leftarrow
n}\DD\left(\M[p^\infty]^\vee\right)\Bigl(\Spec(k) \hookrightarrow
\Spec\left(\WW_n(k)\right)\Bigr).$$
\end{definition}

Call $\TT^{\crys}\left(\M \right)$ the {\it Barsotti-Tate
crystal}\, of the 1-motive $\M$. The functor associating to a
$1$-motive its $p$-divisible group is exact and covariant. The
Dieudonn\'e functor is exact and contravariant. It follows from
\eqref{exactseqinfty} and \eqref{exactseqinfty2} that
$\TT_{\crys}(\M)$ admits Frobenius and Verschiebung operators and
a filtration (respected by Frobenius and Verschiebung): $W_{\geq
0}\left(\TT_{\crys}(\M)\right):=\TT_{\crys}(\M)$,
$W_{-1}\left(\TT_{\crys}(\M)\right):=\TT_{\crys}(\G)$,
$W_{-2}\left(\TT_{\crys}(\M)\right):=\TT_{\crys}(\T)$ and $W_{\leq
-3}\left(\TT_{\crys}(\M)\right):=0$. Hence, $\TT^{\crys}$
(resp.~$\TT_{\crys}$) defines a contravariant (resp.~covariant)
functor from the category of $1$-motives over~$k$ to the category
of filtered $F$-crystals.\\

We remark that the Poincar\'e biextension in \ref{dual} yields a
perfect pairing (\cf \cite[10.2.5]{Deligne})
$$\M[p^n] \otimes \M^\vee[p^n] \to \ZZ/p^n (1)$$
identifying~$\M^\vee[p^n]$ with the Cartier dual~$\M[p^n]^\vee$
of $\M[p^n]$. We thus get a perfect pairing of Barsotti-Tate groups
$$\M[p^\infty]\otimes \M^\vee[p^\infty]\to \QQ_p/\ZZ_p (1)$$
We therefore obtain a perfect pairing of  filtered $\WW (k)$-modules
\begin{equation}\label{crysdual}
\TT^{\crys}(\M)\otimes_{\WW(k)} \TT^{\crys}(\M^\vee)\to \WW(k)
(1).
\end{equation}
We then may regard the Barsotti-Tate crystal of the motivic
Cartier dual as the Dieudonn\'e module of the usual Cartier dual.

\section{Universal vector extension of a 1-motive}
We give some generalities on vector extensions of $1$-motives;
see~\cite[10.1.7]{Deligne}. As in~\cite{Messing}
and~\cite{MazurMessing} this concept will be essential to define
the Dieudonn\'e module of a $1$-motive directly in terms of the
$1$-motive, without using the associated Barsotti-Tate group.

\subsection{Vector extensions of 1-motives}\label{vectorext}
A vector group scheme over~$S$ is a group scheme endowed with an
action of the group ring $\Aff^1$ and isomorphic to~$\GG_{a,S}^r$
(compatibly with the $\Aff^1$-action), for some~$r\in\NN$, locally
on~$S$. If $\cE$ is a (fixed, Zariski) locally free
$\cO_S$-module, denote by~$\W$ the vector group scheme over~$S$
whose sections over an $S$-scheme~$T$ are
$$\W(T):=\Gamma\left(T,\cO_T\tensor_{\cO_S} \cE \right),$$
\eg if $\cE = \cO_S$ then $\W = \GG_{a,S}$.

\noindent A vector extension
$$ 0\llongrightarrow \W \llongrightarrow \E \llongrightarrow  \M
\llongrightarrow 0$$of a $1$-motive
$\M=\left[\X\stackrel{\u}{\longrightarrow} \G\right]$ is an
extension of the complex given by $\M$
(here $\X$ is in degree $-1$ and $\G$ in degree $0$) by a complex
consisting of a vector group scheme~$\W$ over~$S$ concentrated in
degree~$0$.

To give a vector
extension~$\E$ of\/~$\M$ is equivalent to give a complex $\left[\X
\stackrel{\u_\E}{\longrightarrow} \E_\G\right]$, where
\begin{itemize}
\item[\rm i)] $\E_\G$ is a vector extension $$
0\longrightarrow \W \longrightarrow \E_\G \longrightarrow  \G
\longrightarrow 0$$of~$\G$  by a vector group scheme~$\W$
over~$S$;
\item[\rm ii)] $\u_\E\colon\X\longrightarrow \E_\G$  is a group
homomorphism so that the map $\X \stackrel{\u_\E}{\longrightarrow}
\E_\G \longrightarrow \G$ coincides with~$\u$.
\end{itemize}

\noindent Let~$\M=\left[\X\stackrel{\u}{\longrightarrow} \G\right]$ be a
$1$-motive over~$S$. A vector extension
\begin{equation}\label{exactvect}
0\rightarrow
\W\left(\M\right) \rightarrow \E\left(\M\right) \rightarrow  \M
\rightarrow 0
\end{equation}
of\/~$\M$ is called {\it universal}\, if, for any vector extension
$0\rightarrow \W \rightarrow \E \rightarrow  \M \rightarrow 0$,
there exists a unique homomorphism of\/ $S$-vector group schemes
$\phi\colon \W\left(\M\right) \rightarrow \W$ such that~$\E$ is
the push-out of~$\E\left(\M\right)$ by\/~$\phi$.

In the following Sections \ref{vectextA} and \ref{psi} we will show (\cf
\cite[10.1.7]{Deligne}) that there exists a universal vector extension
$\E\left(\M\right)$ over
$S$. Moreover, we clearly have (\cf \cite[1.4]{BS}) an exact sequence:
\begin{equation}\label{exactseqvect}
0\llongrightarrow \E\left(\G\right) \llongrightarrow
\E\left(\M\right)_\G \llongrightarrow \E\left([\X\rightarrow 0]
\right)\llongrightarrow 0.\end{equation}
Following Deligne's notation (see \cite[10.1.11]{Deligne}) we let
$$\TT_{\rm DR} (\M) \df \Lie \left(\E\left(\M\right)_\G
\right).$$The weight filtration on $\M$, defined in
\ref{1motives}, induces a weight filtration on $\TT_{\rm DR}$:
$W_{\geq 0}(\TT_{\rm DR} (\M)):=\TT_{\rm DR} (\M)$,
$W_{-1}(\TT_{\rm DR} (\M)):=\TT_{\rm DR}(\G)$, $W_{-2}(\TT_{\rm
DR} (\M)):=\TT_{\rm DR}(\T)$ and $W_{\leq -3}(\TT_{\rm DR}
(\M)):=0$.

The fact that $W_n(\TT_{\rm DR} (\M))$, as defined above, are
submodules of $\TT_{\rm DR} (\M)$ follows  from \ref{vectextA} and
\ref{psi} below. The short exact sequence (\ref{exactseqvect}) yields a
corresponding short exact sequence of De Rham realizations (\cf
\cite[1.4]{BS}). Denote $\TT^{\rm DR} (\M) \df \TT_{\rm DR}
(\M^\vee)$.

\subsection{Construction of~$\E\left(\A\right)$,
$\E\left(\G\right)$ and~$\E\left(\bigl[\X\rightarrow
0\bigr]\right)$} \label{vectextA} By~\cite{Messing}
or~\cite{MazurMessing} the abelian scheme~$\A$ admits a universal
vector extension~$\E\left(\A\right)$ over~$S$. Let\/~$\A^\vee$ be
the abelian scheme dual to~$\A$. Let~$\Pic^{\natural, 0}$ be the
functor of invertible sheaves in $\Pic^{0}$ endowed with an
integrable $S$-connection (see \cite{MazurMessing}). The map which
forgets the connection gives to this functor $\Pic^{\natural,
0}_{\A^\vee/S}$ the structure of a functor
over~$\A=\Pic^0_{\A^\vee/S}$. One proves (see \cite[I 2.6 and
3.2.3]{MazurMessing}) that this functor is representable and
$$\E\bigl(\A\bigr)\cong
\Pic^{\natural,0}_{\A^\vee/S}.$$In particular, the kernel of the
forgetful map
$$\Pic^{\natural,0}_{\A^\vee/S}\llongrightarrow
\Pic^0_{\A^\vee/S}$$ classifies all possible connections
on the structure sheaf of\/~$\A^\vee$. We conclude that,
if~$\omega_{\A^\vee/S}$ is the $\cO_S$-module of invariant
differentials on~$\A^\vee$ dual to~$\A$, then
$$\W\left(\A\right)\isomarrow \W\left( \omega_{\A^\vee/S}\right)$$in the
notation of~\ref{vectorext}.

\begin{lemma}\label{vectextG} The universal
extension~$\E\left(\G\right)$ exists and is isomorphic
to~$\G\times_\A \E\left(\A\right)$.
\end{lemma}\begin{proof}
Applying the functor~$\Hom_\fppf\left(\_\,
,\GG_{a,S}\right)$ to the exact sequence
$$0 \llongrightarrow \T \llongrightarrow \G \llongrightarrow  \A
\llongrightarrow 0$$ we get the long exact sequence
$$0 \llongrightarrow \Hom\left(\A,\GG_{a,S}\right)
\llongrightarrow \Hom\left(\G,\GG_{a,S}\right) \llongrightarrow
\Hom\left(\T,\GG_{a,S}\right) \llongrightarrow$$
$$\Ext^1\left(\A,\GG_{a,S}\right) \llongrightarrow
\Ext^1\left(\G,\GG_{a,S}\right) \llongrightarrow
\Ext^1\left(\T,\GG_{a,S}\right).$$ Note that
\begin{itemize}
\item[i)] $\Hom\left(\A,\GG_{a,S}\right)=\{0\}$,
since~$f\colon \A\rightarrow S$ is proper, smooth and
geometrically irreducible so that $f_*(\cO_\A)=\cO_S$;
\item[ii)] $\Hom(\T,\GG_{a,S})=\{0\}$, since $\Hom(\GG_{m,S},
\GG_{a,S})=\{0\}$ by a direct computation comparing the
comultiplications on the associated Hopf algebras;
\item[iii)] $\Ext^1\left(\T,\GG_{a,S}\right)=\{0\}$ (\cf
\cite[10.1.7.b)]{Deligne}).
\end{itemize}

\noindent  We conclude that $\Ext^1\left(\A,\GG_{a,S}\right)
\isomarrow \Ext^1\left(\G,\GG_{a,S}\right)$. The map is defined by
sending an extension~$\E$ of~$\A$ by~$\GG_{a,S}$ to~$\G\times_\A
\E$. Hence, $\G$ admits a universal vector
extension~$\E\left(\G\right)$ over~$S$ defined by~$\G\times_\A
\E\left(\A\right)$.\end{proof}

\begin{lemma}\label{vectextX}
{\rm (i)} \enspace  The sheaf\/~$\X\tensor_\ZZ\GG_{a,S} $ is
represented by a vector group scheme over~$S$.\\
{\rm (ii)}\enspace The universal extension~$\E\left(\bigl[\X
\rightarrow 0\bigr]\right)$ is defined by the homomorphism $$\X
\llongrightarrow \X\tensor_\ZZ\GG_{a,S}
$$which sends~$x\mapsto x\tensor 1$.
\end{lemma}
\begin{proof}
Let\/~$S'$ be an \'etale cover of\/~$S$ such that\/~$\X\times_S
S'$ is split. Then, $\bigl(\X\tensor_\ZZ \GG_{a,S}\bigr)\times_S
S'$ is the spectrum of  the symmetric algebra $\Sym(\cI')$
with~$\cI'=\X\times_\ZZ \cO_{S'}$. By descent theory the latter
descends to a locally free $\cO_S$-module $\cI$. Hence,
$\X\tensor_\ZZ \GG_{a,S}$ is represented by the vector group
scheme associated to $\Hom(\cI,\cO_S)$. This proves part~(i).

\noindent The homomorphism in~(ii) is well defined over~$S$.  Due
to the uniqueness of the universal extensions and descent theory
for morphisms, it suffices to prove~(ii)  over~$S'$. Hence, we may
assume that~$\X$ is split. The proof is trivial.\end{proof}

\subsection{Construction of\/~$\E(\M)$}\label{psi}
See \cite[\S 1.4]{BS} for a construction of $\E(\M)$. For our
purposes the following less canonical construction will be useful.
Let\/~$S'$ be an \'etale cover of\/~$S$ such that\/~~$\X\times_S
S' \cong \ZZ^n $, for some $n\in\NN$. Assume that a universal
vector extension~$\E\left(\M\times_S S'\right)$ exists over~$S'$.
The universal property of~$\E\left(\M\times_S S'\right)$ gives the
necessary descent data to get~$\E\left(\M\right)$
descending~$\E\left(\M\times_S S'\right)$ from~$S'$ to~$S$.
Since~$\M$ is defined over~$S$, the descent is effective. Hence,
we may assume that~$\X= \ZZ e_1\oplus\cdots\oplus\ZZ e_n$. In the
notation of~\ref{vectorext}, let
$$\E\left(\M\right)_\G:=\E\left(\G\right)\times_S
\Bigl(\X\tensor_\ZZ \GG_{a,S}\Bigr).$$It is naturally endowed with
the structure of commutative group scheme over~$S$.

\noindent Note that
$$\E\left(\M\right)_\G/\E\left(\G\right) \isomarrow
\X\tensor_\ZZ\GG_{a,S}$$and $\E\left(\M\right)_\G\rightarrow \G$
is a vector extension of~$\G$. Let $$\psi\colon\X \llongrightarrow
\E\left(\G \right)$$be a homomorphism lifting~$\u$; by our
assumption on~$\X$ and possibly replacing $S$ with an \'etale
cover, it exists. Let
$$\u_{\E(\M)}\colon \X \llongrightarrow  \E\left(\M\right)_\G$$be
the homomorphism  $$x=\sum_i x_ie_i \mapsto (\sum_i x_i\cdot
\psi(e_i), \sum_i x_ie_i\tensor 1).$$It is a homomorphism of group
schemes and the composition with~$\E\left(\M\right)_\G \rightarrow
\G$ coincides with~$\u$. By the discussion in \ref{vectorext} the
homomorphism~$\u_{\E(\M)}$ defines a vector extension of~$\M$. We
leave to the reader the proof that it is (canonically isomorphic
to) the universal vector extension of~$\M$ (\cf\cite[1.4]{BS}). By
construction we then get the claimed exact sequence
(\ref{exactseqvect}).

\subsection{$\E(\M)$ via $\E([\X \to \A])$}\label{dualvectext}
Let\/~$\M$ be a $1$-motive and let $\u_\A \colon\X\rightarrow\A$
be the induced 1-motive such that
$$0 \llongrightarrow \T\llongrightarrow \M \llongrightarrow
[\X\rightarrow\A]\llongrightarrow 0$$ is exact. Taking
$\Ext^1\left(-,\GG_{a,S}\right)$ we see that
$$\Ext^1\left([\X\rightarrow\A],\GG_{a,S}\right)\isomarrow
\Ext^1\left(\M,\GG_{a,S}\right)$$ in such a way that $\E([\X \to
\A])$ pulls back to $\E(\M)$.
\bigskip

Let~$\M_1:=\left[\X_1\stackrel{\u_1}{\longrightarrow} \A_1\right]$
and~$\M_2:=\left[\X_2\stackrel{\u_2}{\longrightarrow} \A_2\right]$
be $1$-motives over~$S$ with~$\A_1$ and~$\A_2$ abelian schemes
over~$S$. Assume we are given
\begin{itemize}
\item[\rm 1)] vector
extensions $\E_1$ and~$\E_2$ of~$\M_1$ and~$\M_2$ respectively;
\item[\rm 2)] a $\GG_{m,S}$-biextension~$\P$ of~$\E_{1,\A_1}\times_S
\E_{2,\A_2}$;
\item[\rm 3)] a map $\alpha\colon \X_1 \times_S \X_2 \rightarrow \P$
lifting~$(\u_1,\u_2)$ and inducing a trivialization, as a
biextension, of~$(\u_1,\u_2)^*(\P)$.\end{itemize}

\noindent For every~$x_2\in \X_2$, the element~$\alpha(0,x_2)$
defines a rigidification of~$\P\vert_{\E_1\times_S \{x_2\}}$ over
the $0$ element of~$\E_1\cong \E_1\times_S \{x_2\}$. Because of
the biextension structure of~$\P$, we get
that~$\P\vert_{\E_1\times_S \{x_2\}}$ is endowed with the
structure of a commutative group scheme, having~$\alpha(0,x_2)$
as~$0$ element and sitting in an exact sequence
$$ 0\longrightarrow \GG_{m,S} \longrightarrow
\P\vert_{\E_1\times_S \{x_2\}} \longrightarrow \E_1
\longrightarrow 0.$$Moreover, the map $\X_1 \rightarrow
\P\vert_{\E_1\times_S \{x_2\}}$ defined by~$x_1 \mapsto
\alpha(x_1,x_2)$ is a group homomorphism. Hence, the data\ 1)\ --\ 3)
above define an extension of~$\E_1$ by a torus with character
group~$\X_2$.

\noindent Apply these considerations taking
$\M_1=\left[\X\rightarrow\A\right]$, $\M_2:=\left[\X^\vee
\rightarrow \A^\vee \right]$ (the notation is as in~\ref{dual})
and~$\P$ the pull back of the Poincar\'e biextension
on~$\A\times_S\A^\vee$ to $\E\left(\X \times_S \X^\vee
\rightarrow\A\times_S\A^\vee\right)_{\A\times_S\A^\vee}$. The
universal vector extension~$\E\left(\M\right)$ of\/~$\M$, regarded
as the symmetric avatar of \cite[\S 10.2.12]{Deligne}, is then
equivalent to:

\begin{itemize}
\item[\rm 1)] the universal extension
$$\u_{\E\bigl(\X \times_S \X^\vee
\rightarrow\A\times_S\A^\vee\bigr)}\colon\X \times_S \X^\vee
\rightarrow \E\left(\X \times_S \X^\vee
\rightarrow\A\times_S\A^\vee\right)_{\A\times_S\A^\vee}.$$
\item[\rm 2)] a map $\alpha\colon \X \times_S \X^\vee \rightarrow
\P$ lifting $\u_{\E\bigl(\X \times_S \X^\vee
\rightarrow\A\times_S\A^\vee\bigr)} $ and inducing a
trivialization, as biextension, of the pull back of~$\P$
to~$\X\times_S\X^\vee$. \end{itemize}

\subsection{Universal vector extension of a
Barsotti-Tate group}\label{univvectextpdiv} In~\cite{Messing}
and~\cite{MazurMessing} it is proven that, if~$p$ is locally
nilpotent on~$S$ and~$\underline{G}$ is a Barsotti-Tate group,
there exists a universal vector extension~$\E\left(\underline{G}
\right)$ of~$\underline{G}$.

\smallskip Suppose that $p^N S=0$. As in~\cite[Prop IV.1.10]{Messing},
we may classify the extensions of~$\underline{G}$  by~$\GG_{a,S}$
via the exact sequence $$0 \llongrightarrow \underline{G}[p^N]
\llongrightarrow \underline{G} \stackrel{p^N}{\llongrightarrow}
\underline{G} \llongrightarrow 0$$ applying the
functor~$\Hom_\fppf\left(\_\,,\GG_{a,S}\right)$. Then,
$$\Hom\bigl(\underline{G}[p^N],\GG_{a,S}\bigr)\llongrightarrow
\Ext^1\bigl(\underline{G},\GG_{a,S}\bigr)$$is an isomorphism.
Moreover, by~\cite[Prop IV.1.3]{Messing},
$$\Hom\bigl(\underline{G}[p^N],\GG_{a,S}\bigr)\cong
\Hom\bigl(\omega_{\underline{G}[p^N]^\vee},\GG_{a,S}
\bigr),$$where~$\underline{G}[p^N]^\vee$ is the Cartier dual
to~$\underline{G}[p^N]$. The universal map~$\underline{G}[p^N]
\rightarrow \omega_{\underline{G}[p^N]^\vee}$ is defined by
$$\underline{G}[p^N] \isomarrow
\Hom\left(\underline{G}[p^N]^\vee,\GG_{m,S} \right)
\longrightarrow
\Hom\left(\Inf^1\bigl(\underline{G}[p^N]^\vee\bigr),\GG_{m,S}
\right) \cong \omega_{\underline{G}[p^N]^\vee},$$where
$\Inf^1\bigl(\underline{G}[p^N]^\vee\bigr)\hookrightarrow
\underline{G}[p^N]^\vee $ is the first infinitesimal neighborhood
of the identity.\\

Let~$S_0$ be a scheme on which~$p$ is locally nilpotent and
let~$\underline{G}_0$ be a Barsotti-Tate group on~$S_0$. The
construction above can be extended to define a crystal in fppf
sheaves on the nilpotent crystalline site of~$S_0$ as follows.
Let~$S_0\hookrightarrow S$ be a locally nilpotent thickening with
nilpotent divided powers and let~$\underline{G}$ be a (any) lifting
of~$\underline{G}_0$ to~$S$. Define
\begin{equation}\label{cryspdiv}
\E\left(\underline{G}_0 \right)\bigl(S_0\hookrightarrow S\bigr)
:=\E\left(\underline{G} \right).
\end{equation}
It is proven in~\cite[Thm IV.2.2]{Messing} that indeed this
defines a crystal and is functorial in~$\underline{G}_0$.

Recall that the Dieudonn\'e crystal~$\DD\left(\underline{G}_0
[p^\infty] \right)$, on the crystalline site of\/~$S_0$,
associated to the Barsotti-Tate group~$\underline{G}_0 [p^\infty]$
in~\cite[IV.2.4.3]{Messing}, is by definition the Lie algebra of
the crystal $\E\left(\underline{G}_0\right)$.

\subsection{Another construction
of~$\E\left(\G\right)$ for $p$ nilpotent}\label{secondEG}
Let\/~$\G$ be a semiabelian scheme over~$S$. One can construct the
universal vector extension of~$\G$ exactly as
in~\ref{univvectextpdiv} substituting~$\G$ to~$\underline{G}$.

\smallskip Assume that~$p^N S=0$. The composite map
$$\Hom\bigl(\G[p^N],\GG_{a,S}\bigr)\llongrightarrow
\Ext^1\bigl(\G,\GG_{a,S}\bigr) \llongrightarrow
\Ext^1\bigl(\G\bigl[p^\infty\bigr],\GG_{a,S}\bigr)$$is the
isomorphism in~\ref{univvectextpdiv}. Since all these cohomology
groups are represented by locally free $\cO_S$-modules of the same
rank, all the maps are isomorphisms. Let~$\Delta$ be a splitting
of~$\E\left(\G\right)\rightarrow\G$ Zariski locally on~$\G$. Then,
$\rho:=p^N\cdot \Delta\colon \G \rightarrow \E(\G)$ is a well
defined homomorphism not depending on~$\Delta$. Hence, we have a
diagram

\begin{displaymath}
\begin{CD}
0 @>>> \G\bigl[p^N\bigr] @>>> \G @>{p^N}>> \G @>>> 0\\
&& @VVV   @V{\rho}VV  @V{\wr}VV \\
0 @>>> \W\bigl(\G\bigr) @>>> \E(\G) @>>> \G @>>> 0.\\
\end{CD}
\end{displaymath}Let~$\G[p^\infty]$ be the
Barsotti-Tate group associated to~$\G$. Restricting this diagram
to~$\G[p^\infty]$ and using~\ref{univvectextpdiv} we get:

\begin{lemma}\label{compeasy}{\rm (\cite[Thm V.2.1]{Messing})}
The map $\G[p^N] \llongrightarrow \W\bigl(\G\bigr)$ is universal
for vector extensions of\/~$\G$ over~$S$. Moreover, the canonical
map
$$\E\left(\G[p^\infty]\right)\llongrightarrow
\E\left(\G\right)\times_\G \G[p^\infty]$$is an isomorphism.
\end{lemma}

\subsection{Another construction of $\E\left(\M\right)$ for $p$ nilpotent}
\label{secondEM} Let\/~$\M$ be a $1$-motive over a scheme~$S$ on
which~$p^N=0$. Let~$\rho$ be as in~\ref{secondEG}. By construction
we have $\rho \circ \u= p^N \u_{\E(\M)}$. Let
$$q\colon \X\times_S \G \rightarrow \G$$be the map $(x,g)\mapsto
\u(x)+p^N g$.  Define
$$\Phi\colon \X \times_S \G \llongrightarrow
\E\left(\M\right)_\G$$by $\Phi(x,g):= \u_{\E(\M)}(x)+\rho(g)$.
Then,
$$\Phi\bigl(p^Nx,-\u(x)\bigr)=p^N\u_{\E(\M)}(x)-\rho(\u(x))=0.$$The
composition of~$\Phi$ with the projection $\pi\colon
\E\left(\M\right)_\G\rightarrow\G$ is $(x,g)\mapsto \u(x)+p^N x$.
Hence, $\pi\circ\Phi=q$ and
$$\Ker\left(\pi\circ\Phi\right)=\Ker(q)=\left\{(g,x)\in \G \times_S
\X \vert \u(x)=-p^Ng \right\}.$$By~\ref{Mpinfty} we get a
surjective map~$\Ker(q)\rightarrow \M[p^N]$. Summarizing we get a
diagram
\begin{displaymath}
\begin{CD}
0 @>>> \Ker(q) @>>> \X \times_S \G @>{q}>> \G @>>> 0\\
&& @VVV   @V{\Phi}VV  @V{\Vert}VV \\
0 @>>> \W\bigl(\M\bigr) @>>>
\E\left(\M\right)_\G @>>> \G @>>> 0.\\
\end{CD}
\end{displaymath}Moreover, the induced map $\Ker(q)\rightarrow
\W\bigl(\M\bigr)$ factors as
$$\Ker(q) \llongrightarrow \M[p^N] \llongrightarrow
\W\bigl(\M\bigr).$$This map is functorial in~$\M$. Taking $\G
\rightarrow \M$ and in virtue of~\ref{secondEG}, the composition
of $$\G[p^N] \llongrightarrow \M[p^N] \llongrightarrow
\W\bigl(\M\bigr)$$factors via~$\W\bigl(\G\bigr)$ and is the
universal map for vector extensions. Replacing~$S$ with a suitable
\'etale cover, we may assume that $\X=\ZZ e_1\oplus\cdots\oplus\ZZ
e_n$. By Lemma~\ref{vectextX} and functoriality of~$\Phi$, the
induced map   $$\M[p^N] \llongrightarrow \X/ p^N\X
\llongrightarrow \E\bigl([\X \rightarrow 0]
\bigr)=\X\tensor_\ZZ\GG_{a,S}$$is defined by $x\mapsto x\tensor 1$. Hence, it coincides with the universal map
$$\X/p^N\X \llongrightarrow
\omega_{(\X/p^N\X)^\vee}=\omega_{\X\tensor_\ZZ
\mu_{p^N}}=\X\tensor_\ZZ
\omega_{\mu_{p^N}}$$of~\ref{univvectextpdiv} (the computation is
left to the reader). By~\ref{univvectextpdiv} the map
$\M[p^N]\rightarrow \W\bigl(\M\bigr)$ induced by~$\Phi$ factors
via a unique map $\omega_{\M[p^N]^\vee} \rightarrow
\W\bigl(\M\bigr)$.  Putting together all the remarks given so far,
we conclude that this map is an isomorphism and, hence, it is the
universal one.

\subsection{Comparison between~$\E\left(\M\right)$
and~$\E\left(\M[p^\infty]\right)$ for $p$ nilpotent}
Fix~$n\in\NN$. By~\ref{Mpinfty} the group scheme~$\M[p^n]$ is the
quotient of the subgroup
$$\widetilde{\M[p^n]}:=\left\{(x,g)\in \X \times_S \G \vert
\u(x)=-p^ng \right\}\hookrightarrow\X\times_S \G$$by
$(p^n,-\u)(\X)$. The inclusion $\M[p^n]\hookrightarrow
\M[p^{n+m}]$ is induced by the inclusion~$\widetilde{\M[p^n]}
\hookrightarrow \widetilde{\M[p^{n+m}]}$ defined by~$(x,g)\mapsto
(p^mx,g)$. Multiplication by $p^m\colon \M[p^{n+m}]\rightarrow
\M[p^n]$ is induced by the map~$\widetilde{\M[p^{n+m}]}
\rightarrow \widetilde{\M[p^n]}$  defined by $(x,g)\mapsto (x,p^m
g)$. Let
$$\widetilde{\M[p^\infty]}:=\lim_{n\rightarrow
\infty}\widetilde{\M[p^n]}.$$Then, $\M[p^\infty]$ is the quotient
of~$\widetilde{\M[p^\infty]}$ by the image of~$\X$. Note
that~$\widetilde{\M[p^N]}=\Ker(q)$. Let~$\M[p^\infty]'$ be the
fiber product of~$p^N\colon \M[p^\infty]\rightarrow \M[p^\infty]$
and the inclusion~$\G[p^\infty]\subset\M[p^\infty]$
over~$\M[p^\infty]$. Note that~$\M[p^\infty]'$ surjects
onto~$\G[p^\infty]$ with kernel~$\M[p^N]$.
Let~$\widetilde{\M[p^\infty]}'$ be the fiber product
of~$\widetilde{\M[p^\infty]}\rightarrow \widetilde{\M[p^\infty]}
\rightarrow \M[p^\infty]$(the first map being multiplication
by~$p^N$) and~$\G[p^\infty]\subset\M[p^\infty]$
over~$\M[p^\infty]$. By construction we have a surjective
map~$q'\colon \widetilde{\M[p^\infty]}'\rightarrow \G[p^\infty]$
with kernel~$\widetilde{\M[p^N]}$. If~$(x,g) \in
\widetilde{\M[p^\infty]}' \cap \widetilde{\M[p^n]}$ with~$n\geq
N$, then~$\frac{x}{p^{n-N}}$ lies in~$\X$
and~$q'\bigl((x,g)\bigr)=\u\left(\frac{x}{p^{n-N}}\right)+p^N
g=q\left(\frac{x}{p^{n-N}},g\right)$. Thus, we obtain a
commutative diagram
\begin{displaymath}
\begin{CD}
0 @>>> \Ker(q) @>>> \widetilde{\M[p^\infty]}' @>{q'}>>
\G[p^\infty] @>>>
0\\
&& @V{\Vert}VV   @VVV  @VVV \\
0 @>>> \Ker(q) @>>> \X \times_S \G @>{q}>> \G @>>> 0.\\
\end{CD}
\end{displaymath} Note that~$\M[p^\infty]'$ is the quotient
of~$\widetilde{\M[p^\infty]}'$ by the image of~$\X$. By~\ref{univvectextpdiv} the
push-out of  $\M[p^\infty]'$ along~$\M[p^N] \rightarrow
\omega_{\M[p^N]^\vee}$ is the pull-back
of~$\E\bigl(\M[p^\infty]\bigr)$ via~$\G[p^\infty]\hookrightarrow
\M[p^\infty]$. Hence, the push-out of the extensions in the diagram above via
the composite of $\Ker(q) \rightarrow \M[p^N]\rightarrow
\omega_{\M[p^N]^\vee}$ yields, by \ref{secondEM}, the following
\begin{displaymath}
\begin{CD}
0 @>>> \omega_{\M[p^N]^\vee} @>>>
\E\bigl(\M[p^\infty]\bigr)\times_{\M[p^\infty]} \G[p^\infty]
@>>> \G[p^\infty] @>>> 0\\
&& @V{\wr}VV   @VVV  @VVV \\
0 @>>> \W(\M) @>>> \E(\M)_\G @>>> \G @>>> 0.\\
\end{CD}
\end{displaymath}

\begin{proposition}\label{comparison} Let~$S$ be a scheme on
which\/~$p$ is locally nilpotent. Let\/~$\M$ be a $1$-motive
over~$S$. There is a canonical and functorial isomorphism
$$\E\bigl(\M[p^\infty]\bigr)\times_{\M[p^\infty]} \G[p^\infty]
\isomarrow \E\bigl(\M\bigr)_\G \times_{\G} \G[p^\infty].$$In
particular, it induces an isomorphism on the level of formal
groups of $\E\bigl(\M[p^\infty]\bigr)$ and of $\E\bigl(\M\bigr)$
and of Lie algebras
$$\Lie\Bigl( \E\left(\M[p^\infty]\right)\Bigr) \llongrightarrow
\Lie\Bigl(\E\left(\M\right)_\G\Bigr) = \TT_{\rm DR} (\M).$$
\end{proposition}

\section{Universal extension crystal of a 1-motive}
Let\/~$S_0$ be a scheme such that\/~$p$ is locally nilpotent. Let
$\M_0:=\left[\u_0\colon\X_0\rightarrow \G_0\right]$ be a
$1$-motive over~$S_0$. Let~$S_0\hookrightarrow S$ be a locally
nilpotent pd thickening of~$S_0$: it is a closed immersion defined
by an ideal sheaf~$\cI$ endowed with locally nilpotent divided
powers structure~$\left\{\gamma_n\colon \cI\rightarrow
\cI\right\}_n$. Let $\M:=\left[\u\colon\X\rightarrow \G\right]$
and~$\M':=\left[\u'\colon\X'\rightarrow\G'\right]$ be two
$1$-motives over~$S$ lifting~$\M_0$. We will show that there is a
canonical isomorphism $\E (\M)\cong \E(\M')$.

\subsection{Exp and Log}\label{eXp}
    Let~$\Spec(\cD)$ be an affine flat scheme over~$S$. The
divided power structure on~$\cI$ extends uniquely to a divided
power structure on~$\cI \cD$ by the flatness of~$\cD$;
see~\cite[Cor 3.22]{BerthelotOgus}. We get a homomorphism
$$\exp\colon \Bigl\{a\in \cD\vert
a\equiv 0\, \hbox{{\rm mod }} \cI\,\cD\Bigr\}\to \Bigl\{m\in
\cD\vert m\equiv 1\, \hbox{{\rm mod }} \cI\,\cD\Bigr\}
$$defined by $a
\mapsto \sum_n \gamma_n(a)$. The inverse exists and is the
logarithm $m\mapsto \log\left(m\right):=\sum_n (n-1)!
\gamma_n(m)$. Thus, $\exp$ and $\log$ are isomorphisms.

\subsection{The crystalline nature
of\/~$\E\left(\bigl[\X\rightarrow\A\bigr]\right)$}\label{keyA} The
category of \'etale group schemes over~$S$ is equivalent to the
category of \'etale group schemes over~$S_0$. Hence, we get
canonical isomorphisms
\begin{equation}\label{keyX}
\rho\colon\X\cong \X',
\end{equation}
\noindent since their base change to $S_0$ are isomorphic to
$\X_0$. Moreover, there is a canonical isomorphism of $S$-group
schemes $\sigma\colon\E\left(\A\right)\isomarrow
\E\left(\A'\right)$. This follows from~\cite{Messing}
and~\cite{MazurMessing}: in loc.~cit. a crystal, over the
crystalline site of~$S_0$, is constructed whose value at
$S_0\hookrightarrow S$  is canonically $\E\left(\A\right)$
(equivalently $\E\left(\A'\right)$).

\begin{lemma}\label{key} There is a unique isomorphism of complexes of\/
$S$-group schemes
$$\xi\colon\E\left(\bigl[\X \rightarrow \A\bigr] \right)\isomarrow
\E\left(\bigl[\X'\rightarrow \A'\bigr]\right)$$making the
following diagram commute:

\begin{displaymath}
\begin{CD}
0 @>>>  \E\left(\A\right) @>>> \E\left(\bigl[\X \rightarrow
\A\bigr]\right) @>>> \E\left([\X \rightarrow 0] \right)
@>>>  0 \\
& & @V{\sigma}VV  @V{\xi}VV @V{\rho}VV \\
0 @>>>  \E\left(\A'\right) @>>> \E\left(\bigl[\X' \rightarrow
\A'\bigr] \right) @>>>
\E\left([\X' \rightarrow 0] \right)@>>> 0,\\
\end{CD}
\end{displaymath}

\noindent where the vertical map $\sigma$ on the left is the isomorphism
defined in \cite{MazurMessing} and the vertical map on the right is the
isomorphism
$$\E\left(\bigl[\X\rightarrow 0 \bigr]\right)
\cong \E\left(\bigl[\X'\rightarrow 0\bigr] \right)$$deduced
from~(\ref{keyX}).
\end{lemma}\begin{proof}
By the asserted uniqueness and by descent theory we may assume
that~$\X_0\cong\ZZ e_1\oplus\ldots\oplus\ZZ e_n$ and~$\cI$ is
nilpotent. Let $\pi'\colon\E\left(\A'\right)\rightarrow \A'$ be
the projection. We may assume that~$\Ker(\pi')$ is isomorphic
to~$\GG_{a,S}^{\dim(\A')}$. Choose homomorphisms $\psi\colon\X
\rightarrow \E\left(\A\right)$ and $\psi'\colon\X' \rightarrow
\E\left(\A' \right)$, lifting~$\u$ and~$\u'$ respectively, such
that $\psi\times_S S_0=\psi'\times_S S_0$. By~\ref{psi} we have
$$\E\left(\bigl[\X \rightarrow \A\bigr] \right)_{\A}
=\E\left(\A\right)\times_S \Bigl(\X\tensor_\ZZ
\GG_{a,S}\Bigr)$$and $$\E\left(\bigl[\X' \rightarrow \A'\bigr]
\right)_{\A'}=\E\left(\A'\right)\times_S \Bigl(\X'\tensor_\ZZ
\GG_{a,S}\Bigr).$$Let~$\xi:=\sigma\times \rho$. With the notation
of~\ref{psi}, let $\v:=\u_{\E([\X\rightarrow \A])}$
and~$\v':=\u_{\E([\X'\rightarrow \A'])}$. Consider
$d:=\xi\circ\v-\v'\circ\rho\colon\X \rightarrow \E\left([\X'
\rightarrow \A'] \right)_{\A'}$. Its projection to~$\X'\tensor_\ZZ
\GG_{a,S}$ is~$0$ by construction.

Let $i=1,\ldots,n$. Notice that~$(\pi\circ d)(e_i)\in
\Ker\left(\A(S) \rightarrow \A'(S_0)\right)$. Hence, $\pi(
d(e_i))$ is a point defined in the formal group of~$\A'$.
Since~$\cI$ is nilpotent, such point factors through a finite flat
subgroup scheme~$\A'[p^N]$ for suitable~$N$. Let~$S[\varepsilon]$
be the dual numbers over~$S$. Let~$\A'[p^N]^\vee=:\Spec(\cD)$,
let~$\Delta \colon \cD \rightarrow \cD\tensor_{\cO_S}\cD$ be the
comultiplication and let $\cu\colon \cD\rightarrow \cO_S$ be the
counit. We have
\begin{displaymath}
\begin{split}
\A'[p^N](S) &=\Hom_S\left(\A'[p^N]^\vee,\GG_{m,S} \right)\\
&=\left\{d\in \cD\vert \Delta(x)=x\tensor x,\, \cu(x)=1\right\}\\
&=:\Cospec(\cD).\\
\end{split}
\end{displaymath}

\noindent Note that
\begin{displaymath}
\begin{split}
\Lie\left(\A'[p^N]\right)&=\Bigl\{ \tau\colon\A'[p^N]^\vee
\times_S S[\varepsilon]\rightarrow\GG_{m,S[\varepsilon]}\big\vert
\tau\times_{S[\varepsilon]} S=1\Bigr\}\\ &=\left\{d\in \cD\vert
\Delta(x)=x\tensor 1+1\tensor x,\, \cu(x)=0 \right\}\\
&=:\Prim(\cD).\\
\end{split}
\end{displaymath}

\noindent Since~$\cD$ is a flat $\cO_S$-module, as explained in
\ref{eXp} the exponential and the logarithm define an isomorphism
$$\exp\colon \Prim(\cD)\cap \cI\, \cD \llongrightarrow
\Ker\Bigl(\Cospec(\cD) \rightarrow \Cospec(\cD/\cI\cD)\Bigr);$$
see~\cite[Rmk III.2.2.6]{Messing}. Hence, the exponential defines
an isomorphism
$$\exp\colon\Lie\left(\A'[p^N]\right)\cap \cI\, \cD \isomarrow
\Ker\Bigl(\A'[p^N](S) \rightarrow
\A'[p^N](S_0)\Bigr).$$Let~$y_i:=\exp^{-1}\left((\pi\circ
d)(e_i)\right)$. Let $\GG_{a,S}\rightarrow
\Lie\left(\A'[p^N]\right)$ be the homomorphism of group schemes
sending~$1$ to~$y_i$. The composition with~$\exp$ defines the
unique map $\GG_{a,S}\rightarrow \A'[p^N]$ sending~$1$
to~$\pi(d(e_i))$. Since~$\E\left(\A'\right)$ is the extension
of\/~$\A'$ by~$\GG_{a,S}^{\dim(\A')}$, there exists a unique map
from $\GG_{a,S}$ to $ \E\left(\A'\right)$ sending~$1$ to~$d(e_i)$.
Hence, there exists a unique map $t\colon\X\tensor_\ZZ
\GG_{a,S}\rightarrow \E\left(\A'\right)$ such that, if
\begin{displaymath} \xi:=\begin{pmatrix}\sigma & t\cr 0 & \rho
\end{pmatrix},
\end{displaymath} then $\xi\circ\v=\v'\circ\rho$.
\end{proof}

\subsection{Deformations of biextensions} Let\/~$S_0$ be a scheme.
Let\/~$S_0\hookrightarrow S$ be a thickening defined by an
ideal\/~$\cI$ with nilpotent divided powers.

\begin{lemma}\label{redondant-1} Let\/~$X$ be a flat and separated
group scheme over~$S$. The group of isomorphism classes of
extensions of $X$ by $\GG_{m,S}$, endowed with a trivialization
over $X\times_S S_0$, is naturally isomorphic to the group of
isomorphism classes of extensions of $X$ by $\GG_{a,S}$, endowed
with a trivialization over $X\times_S S_0$.
\smallskip

Given representatives of two corresponding isomorphism classes~$P_m$
and\/~$P_a$, there is a natural isomorphism between the group of
automorphisms of\/~$P_m$, preserving the trivialization
over~$X\times_S S_0$, and the group of automorphisms of\/~$P_a$,
preserving the trivialization over~$X\times_S S_0$.
\end{lemma}
\begin{proof} Let~$\GG=\GG_{m,S}$ or~$\GG_{a,S}$. Let~$P$
be an extension of~$X$ by~$\GG$  and let
$$\alpha\colon P\times_S S_0\isomarrow X\times_S
\GG\times_S S_0$$be a trivialization of the extension
over~$X\times_S S_0$. Choose an open affine covering~$\{U_i\}_i$
of~$X$ and, for each~$i$, a trivialization $$\beta_i\colon
P\times_S U_i\isomarrow  \GG\times_S U_i$$compatible
with~$\alpha$. For $i\neq j$, let~$\Spec(\cD_{i,j}):=U_{i,j}$. The
trivializations~$\beta_i$ and~$\beta_j$ restricted
to~$U_{i,j}:=U_i\times_X U_j$ differ by

\begin{itemize}
\item[\rm i)] a multiplicative cocycle~$m_{i,j}\in \cD_{i,j}^*$ such
that~$m_{i,j}$ is~$1$ mod~$\cI$ if~$\GG=\GG_{m,S}$;

\item[\rm ii)] an additive cocycle~$a_{i,j}\in \cD_{i,j}$ such
that~$a_{i,j}$ is~$0$ mod~$\cI$ if~$\GG=\GG_{a,S}$.
\end{itemize}

\noindent In case~(i) define the extension $$0 \llongrightarrow
\GG_{a,S} \llongrightarrow P_a \llongrightarrow X \llongrightarrow
0$$by the additive cocycle $a_{i,j}:=\log\bigl(m_{i,j} \bigr)\in
\cD_{i,j}$; see \ref{eXp} for $\log$. In case~(ii) define the
extension
$$0 \llongrightarrow \GG_{m,S} \llongrightarrow P_m
\llongrightarrow X \llongrightarrow 0$$by the multiplicative
cocycle $m_{i,j}:=\exp\bigl(a_{i,j} \bigr)\in \cD_{i,j}$. It is
easily checked that these maps define an equivalence between
$\GG_{m,S}$-torsors and~$\GG_{a,S}$-torsors with a trivialization
over~$S_0$.

\noindent Let~$m$, $p_1$ and~$p_2$ be the maps form~$X\times_S X$
to~$X$ defined by the multiplication, the first and the second
projection respectively. Let~$P$ be a $\GG$-torsor over~$X$.  To
give a multiplication law on~$P$, compatible with the one on~$X$
and with the action of~$\GG$ and inducing the standard group law
on~$P\times_S S_0$, is equivalent to give a map
$$X \times_S X \longrightarrow m^*(P) p_1^*(P)^{-1}
p_2^*(P)^{-1}$$reducing to the identity after base change
to~$S_0$. Proceeding as before, using  the logarithm and the
exponential, one passes from the case~$\GG=\GG_{m,S}$ to the
case~$\GG=\GG_{a,S}$ and viceversa. Clearly the commutativity and
the associativity of the multiplication is preserved. Hence, the
conclusion.

The proof of the second part of the
Lemma is left as an exercise for the reader.\end{proof}

\begin{corollary}\label{redondant} Let\/~$A$ and\/~$B$ be
flat and separated group schemes over~$S$. There is a natural
isomorphism between the group of isomorphism classes of
biextensions~{\rm \cite[Ex.~VII]{SGA7}} of\/~$A$ and\/~$B$
by~$\GG_{m,S}$, endowed with a trivialization over~$(A\times_S
S_0)\times_{S_0} (B\times_S S_0)$, and the group of isomorphism
classes of biextensions~{\rm \cite[Ex.~VII]{SGA7}} of\/~$A$
and\/~$B$ by~$\GG_{a,S}$, endowed with a trivialization
over~$(A\times_S S_0)\times_{S_0} (B\times_S S_0)$.\smallskip

Given representatives of two corresponding isomorphism classes~$P_m$
and\/~$P_a$ in this two groups, there is a natural isomorphism
between the group of automorphisms of\/~$P_m$, preserving the
trivialization over~$(A\times_S S_0)\times_{S_0} (B\times_S S_0)$,
and the group of automorphisms of\/~$P_a$, preserving the
trivialization over~$(A\times_S S_0)\times_{S_0} (B\times_S S_0)$.
Such group is isomorphic to
$$\Ker\Bigl(\Bil\left(A\times_S B,\GG_{a,S}\right) \rightarrow
\Bil\left(A_0\times_{S_0} B_0,\GG_{a,S_0}\right)\Bigr),$$where
$\Bil$ stands for the bilinear homomorphisms.
\end{corollary}

\begin{proposition}\label{reduction0}
Let\/~$\A$  be an abelian scheme over~$S$.  Let\/~$X$ be a flat
and separated group scheme over~$S$. Let\/~$\GG$ be~$\GG_{a,S}$
or~$\GG_{m,S}$. Then, the group of  biextensions of $ X$ and $
\E(\A)$ by $\GG$, endowed with a trivialization over $(X\times
\E(\A))\times_S S_0$, is identified with a subgroup of the
homomorphisms from
$X$ to the sheaf of isomorphism classes of extensions of $\W(\A)$
by $\GG$ endowed with a trivialization over $S_0$. Such
identification is compatible with the isomorphisms
in Lemma~\ref{redondant-1} and Corollary~\ref{redondant}.
\end{proposition}
\begin{proof} Let $\cB$ be a biextension of $ X$ and $
\E(\A)$ by $\GG$. For every $S$-scheme $T$ and any $T$-valued
point $x$ of $X$ the fiber $\cB\vert_{\{x\}\times_S \E(\A)}$
defines an extension $\cE_x$ of $\E(\A)\times_S T$ by $\GG\times_S
T$. This defines a map $\vartheta\colon
\Biext^1\Bigl(X,\E(\A);\GG\Bigr) \llongrightarrow
\Hom\Bigl(X,\EXT^1\bigl(\E(\A),\GG\bigr) \Bigr)$.
By~\cite[Ex.~VIII, 1.1.4]{SGA7} the kernel of this map is
$\Ext^1\Bigl(X, \HOM\bigl(\E,\GG\bigr) \Bigr)$.

A trivialization of $\cB$ over $(X\times \E(\A))\times_S S_0$ induces a
trivialization of $\cE_x$ over $\E(\A_0)$. This construction is
compatible with the isomorphisms of\/~\ref{redondant-1} and
\ref{redondant}. Hence, it suffices to consider the
case~$\GG=\GG_{a,S}$. The exact sequence~(\ref{exactvect}) induces
a long exact sequence
\begin{displaymath}
\begin{split}
0  \llongrightarrow &\HOM\bigl(\A,\GG_{a,S}\bigr) \llongrightarrow
\HOM\bigl(\E(\A),\GG_{a,S}\bigr) \llongrightarrow
\HOM\bigl(\W(\A),\GG_{a,S}\bigr)\\ &
\stackrel{\delta}{\llongrightarrow} \EXT^1\bigl(\A,\GG_{a,S}\bigr)
\llongrightarrow \EXT^1\bigl(\E(\A),\GG_{a,S}\bigr)
\llongrightarrow \EXT^1\bigl(\W(\A),\GG_{a,S}\bigr).
\end{split}
\end{displaymath}

\noindent Note that
\begin{itemize}
\item $\HOM\bigl(\A,\GG_{a,S}\bigr)=\{0\}$ since~$\A\rightarrow S$ is
proper;
\item the map~$\delta$ is defined by push forward
of~(\ref{exactvect}) and is an isomorphism by definition of
universal extension.
\end{itemize}

\noindent We conclude
that~$\HOM\bigl(\E(\A),\GG_{a,S}\bigr)=\{0\}$ so that $\vartheta$
is injective for $\GG=\GG_a$. Furthermore, we get an injective map
$\EXT^1\bigl(\E(\A),\GG_{a,S}\bigr) \hookrightarrow
\EXT^1\bigl(\W(\A),\GG_{a,S}\bigr)$. Since the homomorphisms from
$\E(\A_0)$ to $\GG_{a,S_0}$ are trivial, such map is injective
also considering extensions of $\E(\A)$ (resp.~$\W(\A)$) by
$\GG_{a,S}$ endowed with trivialization over $S_0$. \end{proof}

\begin{proposition}\label{reduction1} Let\/~$\M:=[\X\rightarrow\A]$
and\/~$\N:=[\Y\rightarrow \B]$ be $1$-motives over~$S$ with
trivial toric part. Then, kernel of the base change map
$$\Biext^0\left(\E\left(\N\right)_\B,\E\left(\M\right)_\A;
\GG_{m,S}\right) \llongrightarrow
\Biext^0\left(\E\left(\N_0\right)_{\B_0},
\E\left(\M_0\right)_{\A_0};\GG_{m,S_0}\right)$$is isomorphic to
the kernel of
$$\Hom\Bigl(\Y\tensor\X,\GG_{a,S}\Bigr)\llongrightarrow
\Hom\Bigl(\Y_0\tensor\X_0,\GG_{a,S_0}\Bigr).
$$
\end{proposition} \begin{proof} By Corollary~\ref{redondant}, it
suffices to prove the proposition for the kernel of the
homomorphism
$\Biext^0\left(\E\left(\N\right)_\B,\E\left(\M\right)_\A;
\GG_{a,S}\right) \llongrightarrow
\Biext^0\left(\E\left(\N_0\right)_{\B_0},
\E\left(\M_0\right)_{\A_0};\GG_{a,S_0}\right)$. Note that
$\Biext^0\left(\E(\N)_\B\times_S \E(\M)_\A, \GG_{a,S}\right)\cong
\Hom\Bigl(\E(\N)_\B,\Hom\left(\E(\M)_\A,\GG_{a,S}\right)\Bigr)$.
Proceeding as in the proof of Proposition~\ref{reduction0}, get
that $ \Hom\bigl(\E(\M)_\A,\GG_{a,S}\bigr)\cong
\Hom\bigl(\E(\X),\GG_{a,S}\bigr)\cong
\Hom\bigl(\X,\GG_{a,S}\bigr)\cong \X^*\otimes \GG_{a,S}$ with
$\X^*:=\Hom_\ZZ(\X,\ZZ)$. Analogously, since~$\E(\N)_\B$ is the
extension of~$\E(\B)$ by~$\Y\tensor\GG_{a,S}$, we get that
$\Hom\bigl(\E(\N)_\B,\X^*\tensor\GG_{a,S}\bigr)$ is isomorphic to
$\Hom\bigl(\Y, \X^*\tensor\GG_{a,S} \bigr)$. The latter is
$\Hom\bigl(\Y\tensor\X,\GG_{a,S}\bigr)$.
\end{proof}

\subsection{The crystalline nature of the Poincar\'e biextension}
Let\/~$\M_0:=[\X_0\rightarrow\G_0]$ be a
$1$-motive over~$S_0$. Let\/~$\P_0$ be the pull back to the universal
extension
$\E\Bigl(\bigl[\X_0\times_{S_0} \X_0^\vee \rightarrow \A_0
\times_{S_0} \A_0^\vee\bigr]\Bigr)_{\A_0 \times_{S_0} \A_0^\vee}$
of the Poincar\'e biextension on~$\A_0 \times_{S_0} \A_0^\vee$.
Let
$$\alpha_0\colon \X_0\times_{S_0}\X_0^\vee\llongrightarrow \P_0$$be
as in~\ref{dualvectext}. Let~$\M:=[\X\rightarrow \G]$
and\/~$\M':=[\X'\rightarrow \G']$ be $1$-motives over~$S$
lifting~$\M_0$. Let\/
$$\sigma\colon\E\Bigl(\bigl[\X\times_S \X^\vee \rightarrow \A \times_S
\A^\vee\bigr]\Bigr) \isomarrow \E\Bigl(\bigl[\X'\times_S
{\X^\vee}' \rightarrow \A' \times_S {\A'}^\vee\bigr]\Bigr) $$be
the isomorphism defined in Lemma~\ref{key}. Let
$$\P \llongrightarrow \E\Bigl(\bigl[\X\times_S \X^\vee
\rightarrow \A \times_S \A^\vee\bigr]\Bigr)_{\A \times_S
\A^\vee}$$and $$\P'\llongrightarrow \E\Bigl(\bigl[\X'\times_S
{\X^\vee}' \rightarrow \A' \times_S {\A'}^\vee\bigr]\Bigr)_{\A'
\times_S {\A'}^\vee}$$be the pull back of the Poincar\'e
biextension on~$\A \times_S \A^\vee$ (resp.~on $\A' \times_S
{\A'}^\vee$).  Let
$$\alpha\colon \X\times_S \X^\vee \llongrightarrow
\P\qquad\hbox{{\rm and}}\qquad \alpha_2\colon \X'\times_S
{\X^\vee}' \llongrightarrow \P'$$be the maps
defining~$\E\bigl(\M\bigr)$ (resp.~$\E\bigl(\M'\bigr)$) as in~2)
of\/~\ref{dualvectext}. Then,

\begin{proposition}\label{finalreduction} There is a unique
isomorphism
$$\P\isomarrow\sigma^*(\P'),$$as biextensions
of\/~$\E\Bigl(\bigl[\X\times_S \X^\vee\rightarrow \A\times_S
\A^\vee\bigr]\Bigr)_{\A\times_S \A^\vee}$, compatible with the
trivializations~$\alpha$ and\/~$\alpha'$ via the identification
$\X\times_S \X^\vee \cong \X'\times_S {\X^\vee}'$
of\/~(\ref{keyX}).
\end{proposition}
\begin{proof} By the claimed uniqueness and using descent,
we may assume that~$\X_0\cong\ZZ^r$ and that the toric part
of~$\G_0$ is split. This is equivalent to require that~$\X$,
$\X^\vee$, $\X'$ and~${\X^\vee}'$ are constant group schemes. In
particular, the map~$\E\Bigl(\bigl[\X\times_S \X^\vee \rightarrow
\A \times_S \A^\vee\bigr]\Bigr)_{\A\times_S \A^\vee}\rightarrow
\A\times_S \A^\vee$ can be factored via $$
\E\Bigl(\bigl[\X\times_S \X^\vee \rightarrow \A \times_S
\A^\vee\bigr]\Bigr)_{\A\times_S \A^\vee}
\stackrel{p}{\llongrightarrow}\E\bigl(\A\times
\A^\vee\bigr)\stackrel{\pi}{\llongrightarrow} \A\times_S
\A^\vee.$$We conclude that~$\P$ is also the pull back via~$p$ of
the pull back~$\Q$ of the Poincar\'e biextension on~$\A\times
\A^\vee $ via~$\pi$ and, hence, does not depend on the choice
of~$p$. Denote by~$\Q'$ the pull back to~$\E\bigl(\A'\times_S
{\A'}^\vee)$ of the Poincar\'e biextension
on~$\A'\times_S{\A'}^\vee$. Consider the  isomorphism
$$\varrho\colon \E\Bigl( \A\times_S \A^\vee\Bigr)\isomarrow
\E\bigl( \A\bigr)\times_S \E\bigl(\A^\vee\bigr)\isomarrow \E\bigl(
\A'\bigr)\times_S \E\bigl({\A'}^\vee\bigr) \isomarrow \E\Bigl(
\A'\times_S {\A'}^\vee\Bigr)$$defined in~\ref{keyA}.
Let~$\GG=\GG_{m,S}$. Consider the element $
\varrho^*\bigl(\Q'\bigr)\Q^{-1}$: it is a biextension of $\E(\A)$
and $\E(\A^\vee)$ by $\GG_{m,S}$ whose base change to
$\E(\A_0)\times \E(\A_0^\vee)$ is canonically trivialized. By
Proposition~\ref{reduction0} such biextensions form a subgroup of
the group of homomorphisms from $\E(\A)$ to the sheaf of
isomorphism classes of extensions of $\W(\A^\vee)$ by $\GG_{a,S}$
endowed with a trivialization over $S_0$. Since~$\A$ is a
divisible sheaf and $\W(\A^\vee)$ and $\GG_{a,S}$ are torsion
sheaves, such homomorphisms are a subgroup of the group $\cK$ of
homomorphisms from $\W(\A)$ to the sheaf of isomorphism classes of
extensions of $\W(\A^\vee)$ by $\GG_{a,S}$ (equivalently
$\GG_{m,S}$) endowed with a trivialization over $S_0$. Note
that~$\Q$ defines the trivial $\GG_{m,S}$-biextension of~$\W(\A)$
and~$\W(\A^\vee)$ since it is pulled-back from~$\A\times_S
\A^\vee$. Thus, the image of~$\varrho^*\bigl(\Q'\bigr)\Q^{-1}$
in~$\cK$ coincides with the pull-back of~$\P'$  via the map
$\W(\A)\times_S \W(\A^\vee) \rightarrow \A'\times_S {\A'}^\vee$
induced by~$\varrho$. Hence, the image
of~$\varrho^*\bigl(\Q'\bigr)$ in~$\cK$ lies in the image of
homomorphisms from $\A'$ to the sheaf of isomorphism classes of
extensions of $\W(\A^\vee)$ by $\GG_{a,S}$ endowed with a
trivialization over $S_0$, which is trivial. Therefore,
$\varrho^*\bigl(\Q'\bigr)$ is isomorphic to~$\Q$ as biextension.
We conclude that~$\P$ is isomorphic to~$\sigma^*(\P')$ as
biextensions
lifting $\P_0$.\\

By Proposition~\ref{reduction1} the automorphism group of a
biextension over~$\E\Bigl(\bigl[\X\times_S \X^\vee \rightarrow \A
\times_S \A^\vee\bigr]\Bigr)$, reducing to the identity after base
change to~$S_0$, is isomorphic to the kernel of
$$\Hom_\fppf\Bigl(\X\tensor\X^\vee,\GG_{a,S}\Bigr)\llongrightarrow
\Hom_\fppf\Bigl(\X_0\tensor\X_0^\vee,\GG_{a,S_0}\Bigr).
$$The set of maps~$\alpha\colon \X\times_S \X^\vee \rightarrow
\P$, as in 3) of~\ref{dualvectext}, deforming the map~$\alpha_0$
defining~$\M_0$, is a principal homogenous space under the kernel
of
$$\Hom_\fppf\Bigl(\X\tensor\X^\vee,\GG_{a,S}\Bigr)\llongrightarrow
\Hom_\fppf\Bigl(\X_0\tensor\X_0^\vee,\GG_{a,S_0}\Bigr).
$$Hence, the conclusion.
\end{proof}

\begin{corollary}\label{gory} There is a unique  isomorphism of
$S$-group
schemes
$$\zeta\colon\E\left(\M\right)\isomarrow
\E\left(\M'\right)$$making the following diagram commute
\begin{displaymath}
\begin{CD}
0 @>>> \T @>>> \E\left(\M \right) @>>> \E\left(\bigl[\X
\rightarrow \A\bigr] \right)@>>> 0 \\ && @V{\tau}VV   @V{\zeta}VV
@V{\xi}VV \\ 0 @>>> \T' @>>>
\E\left(\M'\right)@>>>\E\left(\bigl[\X' \rightarrow \A'\bigr]
\right)@>>> 0,
\end{CD}
\end{displaymath}where $\xi$ is defined in Lemma~\ref{key} and $\tau$
is the canonical isomorphism between the tori~$\T$ and~$\T'$ as
deformations of\/~$\T_0$.
\end{corollary}\begin{proof}
By~\ref{dualvectext} the universal extension of $\M$ is defined by
the universal extension $\E\bigl([\X\times_S\X^\vee\rightarrow \A
\times_S \A^\vee]\bigr)$ and a trivialization $\alpha$ of the
pull-back $\cP$ to $\X\times \X^\vee$. Analogous description
exists for $\M'$.  By Proposition~\ref{finalreduction} we get an
isomorphism $\E(\M)\to\E(\M')$ with the claimed properties. Since
the homomorphims from a vector group scheme to a torus and from an
abelian scheme to a torus are trivial, we have
$\Hom\bigl(\E([\X\to\A]),\T'\bigr)=0$. This implies the claimed
uniqueness.
\end{proof}

\subsection{Crystals via universal extensions}\label{crysviaext}
Let\/~$S_0$ be a scheme such that\/~$p$ is locally nilpotent.
Let\/~$\M_0$ be a $1$-motive over~$S_0$. Define the crystal of
group schemes $ \E\left(\M_0\right)$ on the nilpotent crystalline
site of\/~$S_0$. Let\/~$S_0\hookrightarrow S$ be a locally
nilpotent pd thickening of\/~$S_0$. Let\/~$\M$ be a $1$-motive
lifting~$\M_0$ to~$S$. Then,
\begin{equation}\label{E(M)cris}\E\left(\M_0\right)\left(S_0
\hookrightarrow S\right):= \E\left(\M\right).\end{equation} Here
$\M_0\rightarrow S_0$ can be lifted locally and by virtue
of\/~Corollary~\ref{gory} we get indeed a crystal which we call
the {\it universal extension crystal}\, of the 1-motive.

Define the contravariant functor associating to a $1$-motive
$\M_0$ over~$S_0$ the crystal\/~$\DD\left(\M_0\right)$ on the
nilpotent crystalline site of\/~$S_0$ as follows:
\begin{equation}\label{Lie(E(M))cris}
\DD\left(\M_0\right)\left(S_0\hookrightarrow S\right)
:= \TT^{\rm DR} (\M).
\end{equation}

See \ref{vectorext} for the notation $\TT^{\rm DR} (\M)$.
Furthermore, let\/~$\P_0$ be the pull back
to~$\E\Bigl(\bigl[\X_0\times_{S_0} \X_0^\vee\rightarrow \A_0
\times_{S_0} \A_0^\vee\bigr]\Bigr)$ of the Poincar\'e biextension
on~$\A_0 \times_{S_0} \A_0^\vee$. Define the {\it Poincar\'e
crystal}\, of biextensions $\P_0$ over the crystal~$\E(\M_0)$ on
the crystalline site of\/~$S_0$ as follows.
Let\/~$S_0\hookrightarrow S$ be as above and let\/~$\M$ be a
$1$-motive lifting~$\M_0$ to~$S$. Let\/~$\P$ be  the pull back to
$\E\Bigl(\X\times_S \X^\vee \rightarrow \A \times_S \A^\vee\Bigr)$
of the Poincar\'e biextension on~$\A \times_S \A^\vee$. Then,
\begin{equation}\label{Pcris}
\P_0\left(S_0\hookrightarrow S\right):=\P(S).\end{equation}
Due to Proposition~\ref{finalreduction} this defines a crystal.
Define the
$\cO_{S_0}^\crys$-bilinear pairing of crystals
$$\Phi\colon\DD\left(\M_0\right)\times \DD\left(\M_0^\vee\right)
\llongrightarrow \DD\left(\GG_{m,S_0}^\vee\right)$$as follows.
By~{\rm \cite[Prop 10.2.7.4]{Deligne}} there is a unique
$\natural$-structure on the biextension~$\P$. See~{\rm
\cite[10.2.7.2]{Deligne}} for this notion. Its curvature defines
an $\cO_S$-bilinear map {\rm \cite[10.2.7.3]{Deligne}}
\begin{equation}\label{dualityonLie(E(M))cris}
\Phi\colon \Lie\Bigl(\E\left(\M^\vee\right) \Bigr)\tensor_{\cO_S}
\Lie\Bigl(\E\left(\M\right) \Bigr) \llongrightarrow
\Lie\Bigl(\GG_{m,S}\Bigr)
\end{equation}
where the Lie algebras are applied to the degree zero component of
the complexes, \ie to connected group schemes, yielding a pairing
of filtered vector bundles over $S$.

\begin{remark} {\rm If $S_0 \subset S$ is a pd thickening defined by a
locally nilpotent ideal, but the pd structure is {\it not} locally
nilpotent, then the universal extension of a lifting of $\M_0$ to
$S$ need not be crystalline. This holds, for example, if $p=2$ and
$S$ is a scheme over $\WW(k)$ with pd structure compatible with
the one on the maximal ideal of $\WW(k)$. On the other hand,
thanks to Proposition~\ref{comparison} and \cite[\S 11]{MazurMessing},
its
formal completion at the origin is crystalline. Thus, the crystal
$\DD\left(\M_0\right)$ is defined on the full crystalline site of
$S_0$.}
\end{remark}

\section{Proof of Theorem A${}^\prime$}
Let $\M_0$ be a 1-motive defined over a perfect field $S_0 = \Spec (k)$
of
positive  characteristic $p>0$. Let $\M$ denote a lifting to $S = \Spec
(\WW_n (k))$. We first show that the crystal \eqref{Lie(E(M))cris}
is a filtered $F$-crystal, \ie we show how to get Frobenius and
Verschiebung compatibly with the weight filtration.

\subsection{Frobenius and Verschiebung}\label{FROBVER}
    Let $\sigma \colon \Spec(k) \to
\Spec(k)$ be the Frobenius map defined on $k$ by $x \mapsto x^p$.
By abuse of notation we denote by $\sigma\colon
\Spec\bigl(\WW_n (k)\bigr)\to \Spec\bigl(\WW_n (k)\bigr)$ also the map
associated to the Frobenius map on Witt vectors. Let $\M^{(p)}_0$ be
the pull-back of $\M_0$ via $\sigma$ (as complex of group schemes,
\ie $\M^{(p)}_0:=[\X^{(p)}_0\to \G^{(p)}_0]$). The Frobenius map on
$\M_0$
defines a morphism of $1$-motives $F\colon \M_0\to \M^{(p)}_0$
over~$k$. Given the lifting $\M$ of $\M_0$ to
$\WW_n (k)$, then $\sigma^*(\M)$ is a lift of $\M^{(p)}_0$ to
the truncated Witt vectors $\WW_n (k)$. Hence, we obtain a map, called
{\it
Verschiebung},
$$V\colon \Lie\left(\E(\M)\right)\llongrightarrow
\Lie\left(\E\bigl(\sigma^*(\M)\bigr)\right)\cong
\Lie\bigl(\E(\M)\bigr) \otimes_\sigma \WW_n (k).$$The Frobenius map
on $\M^\vee_0$ defines by duality a morphism of $1$-motives $V\colon
\M^{(p)}_0\to \M_0$ over~$k$. Note that $F \circ V= p$ and $V\circ
F=p$: this is clear on the lattice $\X_0$, on the torus $\T_0$ and it
is classical on the the abelian part $\A_0$. The Verschiebung on
$\M_0$ defines a homomorphism, called {\it Frobenius},
$$F\colon \Lie\left(\E(\M)\right)\otimes_\sigma \WW_n (k)
\llongrightarrow \Lie\bigl(\E(\M)\bigr).
$$Then, $F \circ V= p$ and $V\circ F=p$. Analogously, Frobenius
on $\M^\vee_0$ defines $V$ on $\Lie\bigl(\E(\M^\vee)\bigr)$ and
Frobenius on $\M_0$ defines, by duality, Verschiebung on $\M^\vee_0$
and Frobenius on $\Lie\bigl(\E(\M^\vee)\bigr)$. Clearly, $F$ and
$V$ preserve the weight filtrations. Passing to the limit over $n$ we
then obtain
$F$ and $V$ on $\liminv{n}\DD\left(\M_0\right)\bigl(\Spec(k) \subset
\Spec(\WW_n(k))\bigr)$ and $
\liminv{n}\DD\left(\M^\vee_0\right)\bigl(\Spec(k) \subset
\Spec(\WW_n(k))\bigr)$ over $\WW (k)$.\\

Note that, for~$\M:=\GG_m$ we get~$\T_\crys(\M)=\WW(k)$ as
$\WW (k)$-module, with~$F$ sending~$1\mapsto 1$ and~$V$
sending~$1\mapsto p$ ,\ie
$\T_\crys([0\to \GG_m])=\WW(k)(1)$ as $F$-crystals according with the
notation
adopted above.\\

Finally, we remark that Frobenius defines a map from the
Poincar\'e biextension $P_0$ on $\E\Bigl(\bigl[\X_0\times_{S_0}
\X_0^\vee\rightarrow \A_0 \times_{S_0} \A_0^\vee\bigr]\Bigr)$ to
the Poincar\'e biextension $P_0^{(p)}$ on
$\E\Bigl(\bigl[\X_0^{(p)}\times_{S_0} (\X_0^\vee)^{(p)}\rightarrow
\A_0^{(p)} \times_{S_0} (\A_0^\vee)^{(p)}\bigr]\Bigr)$. It follows
that, for $p\geq 3$, the pairing~\eqref{dualityonLie(E(M))cris}
satisfies $V_{\GG_m}\circ \Phi_\M= p \sigma^{-1} \circ
\Phi_{\M}=\Phi_{\M^{(p)}} \circ (V_{\M^\vee}\times V_{\M})$ and,
thus, $p F_{\GG_m} \circ \Phi_{\M^{(p)}}= \Phi_\M \circ
\left(F_{\M^\vee}\times F_{\M}\right)$.

\subsection{The comparison}
We then show that $\TT^{\crys}(\M_0) = \TT^{\rm DR}(\M)$.
By \cite[\S 15]{MazurMessing}
the contravariant Dieudonn\'e module of a $p$-divisible group over
$k$ is canonically isomorphic to the Lie algebra of the universal
extension of a lifting to $\WW (k)$ of the dual $p$-divisible
group. It follows from~\ref{comparison}, applied to
$S=\Spec(\WW_n(k))$ that
$\DD\left(\M_0\right)(S_0\subset S)= \TT^{\rm DR}(\M)$, defined
in~(\ref{Lie(E(M))cris}), and\/~$\DD\left(\M_0[p^\infty]
\right)(S_0\subset S)$, defined in~(\ref{dieu}), are canonically
isomorphic as filtered $\WW_n(k)$-modules. By definition of
$\TT^{\crys}$ and $\TT_{\crys}$, see Definition~\ref{defTcrysupanddwn},
we get
canonical isomorphisms $\TT^{\crys}\left(\M_0\right)\cong
\liminv{n}\DD\left(\M_0\right)\bigl(\Spec(k) \subset
\Spec(\WW_n(k))\bigr)$ and, since
$\M_0[p^\infty]^\vee=\M^\vee_0[p^\infty]$,
$\TT_{\crys}\left(\M_0\right)\cong
\liminv{n}\DD\left(\M_0^\vee\right)\bigl(\Spec(k) \subset
\Spec(\WW_n(k))\bigr)$ as filtered $\WW(k)$--modules. By the
functoriality claimed in Proposition~\ref{comparison}, they are
compatible
with Frobenius and Verschiebung. This proves the first claim of
Theorem A${}^\prime$.

The second claim of Theorem A${}^\prime$ follows from the above
and \ref{FROBVER}, at least for $p\geq 3$, or from
\eqref{crysdual} for any $p$. Note that, for $p\geq 3$, both the
pairings \eqref{dualityonLie(E(M))cris} and  \eqref{crysdual} are
induced by the Poincar\'e extension on a formal lifting of
$\M[p^\infty]$ and $\M$, respectively, to $\WW(k)$. One can prove
that, indeed, the two pairings are the same. Since we will not use
this compatibility in the sequel, we omit the proof and leave it
to the reader.\\[10pt]
We also have obtained the following result.
\begin{corollary} Let $R$ be a complete discrete valuation
ring  with perfect residue field\/~$k$ and with ramification index
smaller than $p-1$. Let\/~$\M$ be a $1$-motive over~$\Spec(R)$.
Then, there is a canonical isomorphism
$$\TT_{\crys}\left(\M_k\right) \tensor_{\WW(k)} R \cong
\TT_{\rm DR}\left(\M\right).$$\end{corollary}
\begin{proof} The hypothesis on the ramification implies that
the maximal ideal of $R$ is endowed with canonical divided powers
structure compatible with that on the maximal ideal of $\WW(k)$.
The $1$-motive $\M_n$ obtained base changing $\M$ to $R/(p^n)$ is
a lifting of $\M_k$. Let $\M_n'$ be a lifting of $\M_k$ to
$\WW_n(k)$. It follows from~\ref{comparison} that we have
canonical isomorphisms  $\Lie\left(\E(\M_n[p^\infty])\right)\cong
\Lie\left(\E(\M_n'[p^\infty])\otimes_{\WW(k)} R \right) \cong
\DD(\M_k[p^\infty])\bigl(\Spec(k)\subset
\Spec(\WW_n(k))\bigr)\otimes_{\WW(k)} R$.
\end{proof}

\subsection{Properties of $\TT_\crys(\M)$}\label{propTcrys} We
explicitly state some of the properties of
$\TT_{\crys}\left(\M\right)$ (similarly for
$\TT^{\crys}\left(\M\right)$), which can be deduced using either
its definition as the Dieudonn\'e module of $\M[p^\infty]^\vee$ or
via $\DD(\M)$:
\begin{description}
\item[\rm 1)] it is a free $\WW(k)$-module of rank equal
to~$\dim(\G)+\rk(\X)$;
\item[\rm 2)] it admits a weight filtration~$W$
\begin{description}
\item[\rm 2.a)] $W_{\geq
0}\Bigl(\TT_{\crys}\left(\M\right)\Bigr)=\TT_{\crys}\left(\M\right)$;
\item[\rm 2.b)] $W_{-1}\Bigl(\TT_{\crys}\left(\M\right)\Bigr)
=\TT_{\crys}\left(\G\right)$. It is a free $\WW(k)$-module of rank
$\dim(\G)$;
\item[\rm 2.c)] $W_{-2}\Bigl(\TT_{\crys}\left(\M\right)\Bigr)
=\TT_{\crys}\left(\T\right)$. It is a free $\WW(k)$-module of rank
$\dim(\T)$;
\item[\rm 2.d)] $W_{\leq
-3}\Bigl(\TT_{\crys}\left(\M\right)\Bigr)=0$;
\end{description}
\item[\rm 3)] the graded pieces of the weight filtration satisfy
\begin{description}
\item[\rm 3.a)] $ \Gr_{-2}^W\Bigl( \TT_{\crys}\left(\M \right)\Bigr):=
\TT_{\crys}\left([0\rightarrow \T] \right)$. It is a free
$\WW(k)$-module of rank $\dim(\T)$;
\item[\rm 3.b)]$ \Gr^W_{-1}\Bigl( \TT_{\crys}\left(\M \right)\Bigr):=
\TT_{\crys}\left([0\rightarrow\A]\right)$. It is a free
$\WW(k)$-module of rank $\dim(\A)$;
\item[\rm 3.c)] $\Gr^W_{0}\Bigl( \TT_{\crys}\left(\M \right)\Bigr):=
\TT_{\crys}\left([\X\rightarrow 0] \right)$. It is a free
$\WW(k)$-module of rank $\rk(\X)$.
\end{description}
\item[\rm 4)] Let $\sigma $ be the Frobenius homomorphism
on $\WW(k)$. Then $\TT_\crys(\M)$ is endowed with a
$\sigma$-linear morphism~$\F$ and a $\sigma^{-1}$-linear
morphism~$\V$ such that
\begin{description}
\item[\rm 4.a)] they respect the weight filtration;
\item[\rm 4.b)] $\F\circ\V=\V\circ\F=p$;
\item[\rm 4.c)] $\V$ is an isomorphism on $ \Gr_{0}^W\Bigl(
\TT_{\crys}\left(\M \right)\Bigr)$;
\item[\rm 4.d)] $\F$ is an isomorphism on $ \Gr_{-2}^W\Bigl(
\TT_{\crys}\left(\M \right)\Bigr)$;
\end{description}
\item[\rm 5)] there is a perfect bilinear pairing of filtered
$\WW(k)$--modules
\begin{equation}
\langle -,-\rangle\colon \TT_{\crys}(\M)\otimes_{\WW(k)}
\TT^{\crys}(\M)\to \WW(k)(1).
\end{equation} satisfying $\langle -,-\rangle^{\sigma^{-1}}=
\langle{ V(-),V(-) }\rangle$ and $p \langle -,-\rangle^{\sigma}=
\langle{ F(-),F(-) }\rangle$.
\end{description}

\section{Picard 1-motives and $\natural$-structures}
Let $S_0=\Spec (k)$ where $k$ is any field. Let $\pi\colon X_{\d}\to
S_0$ be a
simplicial $S_0$-scheme. Note that the simplicial structure
provide a complex of $S_0$-schemes in the sense of~\cite[\S 2]{BM}.
Let $\pi_i\colon X_i\to S_0$ denote the structure morphisms.  We
say that $\pi$ is proper, smooth, etc. if each $\pi_i$ is proper,
smooth, etc.  Assume $X_{\d}$ proper over $S_0$ so that, for each
component
$X_i$, the usual Picard fppf-sheaf $\Pic_{X_i/S_0}\df
R^1(\pi_i)_*\G_{m, X_i}$ is representable by a (commutative) group
scheme, locally of finite type over $S_0$. We refer to \cite{BL} for the
general framework of Picard functors.

If, in addition $S_0 =\Spec(k)$ is a perfect field, and  $X_i$ is
smooth over $S_0$ then the reduced connected component of the identity
$\Pic_{X_i/S_0}^{0, \red}$ is an abelian scheme over
$S_0$, \ie this is the case if $\pi$ is proper and smooth over
$S_0$ (see \cite{BL}).

\subsection{$\bPic^+$}\label{simpic}
For any simplicial $S_0$-scheme $X_{\d}$ let $\bPic (X_{\d})$ be
the group of isomorphism classes of simplicial line bundles on $X_{\d}$
(we refer to  \cite[\S 4.1 \& A.3]{BS} for the basic properties of the
simplicial
Picard functor). By descent we see that
$$\bPic (X_{\d})\cong \HH^1_{\et}(X_{\d},
\G_{m, X_{\d}})\cong \HH^1_{\fppf}(X_{\d}, \G_{m, X_{\d}}).$$
Denote $T \leadsto \bPic_{X_{\d}/S_0}(T)$ the Picard fppf-sheaf
obtained by sheafifying the functor
$$T\leadsto \bPic (X_{\d}\times_{S_0} T)$$ with respect to the
fppf-topology on $S_0$, \ie if $\pi: X_{\d}\times_{S_0}T\to T$,
then
$$\bPic_{X_{\d}/S_0}(T)\cong
\H^0_{\fppf}(T,R^1\pi_*\G_{m,X_{\d}\times_{S_0}T}).$$ Considering
the canonical spectral sequence (see \cite[5.2.7.1 ]{Deligne})
\begin{equation}\label{sscomp}
E^{p,q}_1 = R^q(\pi_p)_*\G_{m, X_p}\implies R^{p+q}\pi_*\G_{m,
X_{\d}}
\end{equation}
we obtain an exact sequence of fppf-sheaves:
\begin{equation}\label{semisimp}
0\to \frac{\ker ((\pi_1)_*\G_{m,X_1}\to (\pi_2)_*\G_{m,X_2})}{\Im
((\pi_0)_*\G_{m,X_0}\to (\pi_1)_*\G_{m,X_1})}\to
\bPic_{X_{\d}/S_0} \to \Ker (\Pic_{X_0/S_0}\to\Pic_{X_1/S_0})
\end{equation}
We then have the following result.
\begin{lemma}\label{reppic} Let $X_{\d}$ be a smooth proper
simplicial $S_0$-scheme, $S_0 =\Spec(k)$. Then the Picard $\fppf$-sheaf
$\bPic_{X_{\d}/S_0}$ is representable by a group scheme, locally of
finite type over $S_0$. Denote by $\bPic_{X_{\d}/S_0}^{0,\red}$ the
connected component of the identity of the simplicial Picard scheme of
$X_{\d}$ over
a perfect field $S_0 =\Spec(k)$, endowed with its reduced structure: it
is a semi-abelian scheme, extension of the abelian scheme
$$\Ker^{0, \red}(\Pic_{X_0/S_0}^{0,
\red}\to\Pic_{X_1/S_0}^{0, \red})$$ by a torus~$T$ (here
$\Ker^{0, \red}$ denote the identity component of the kernel
endowed with its reduced structure).
\end{lemma}
\begin{proof} It follows from \cite[Lemma 4.1.2]{BS} and the proof of
\cite[Prop. 4.1.3]{BS} (see also \cite[3.2-3.5]{RA}). Actually, using
(\ref{sscomp}) it follows from (\ref{semisimp}) and the
representability of $\Pic_{X_i/S_0}$, since the fppf-sheaf
$(\pi_i)_*\G_{m, X_i}$ is representable by a torus and the resulting
extensions are representable, \eg by \cite[Prop. 17.4]{OG}.
\end{proof}

\begin{remark}\label{withouttorsion}
{\rm We describe the torus $T$ appearing in Lemma~\ref{reppic}.
Consider the homology of the following complex of tori (\cf
\eqref{sscomp} and
\eqref{semisimp})
$$(\pi_0)_*\G_{m,X_0}\to (\pi_1)_*\G_{m,X_1}\to (\pi_2)_*\G_{m,X_2}.$$
The torus $T$ is given by the quotient of $T_2:=\Ker^{0, \red}
((\pi_1)_*\G_{m,X_1}\to
(\pi_2)_*\G_{m,X_2})$ by $T_1:=\Im ((\pi_0)_*\G_{m,X_0}\to
(\pi_1)_*\G_{m,X_1})$. We also describe the cocharacter group of $T$ as
follows.
Assume that $S_0 =\Spec (k)$ and $k$ is algebraically closed. Let
$C_2$, $C_1$ and $C_0$ be the free abelian groups generated
by the connected (= irreducible) components of $X_2$, $X_1$ and $X_0$
respectively.
Let $$C_2\stackrel{d_2}{\longrightarrow} C_1 \stackrel{d_1}
{\longrightarrow} C_0 $$be the group homomorphisms defined by the
alternating sums of the faces maps of the simplicial scheme. Let
$$C^0 \stackrel{d^1}{\longrightarrow} C^1 \stackrel{d^2}
{\longrightarrow} C^2$$be the dual complex, actually given by the
cocharacters of the above complex of tori. The map $C_1\to C_0$
associates to a component $X_{1,j}$ of $X_1$ the element
$X_{0,h}-X_{0,i}$ where $d_0^1(X_{1,j})\subset X_{0,h}$ and
$d_1^1(X_{1,j})\subset X_{0,i}$. In particular, the image of
$C_1\to C_0$ is a direct summand of $C_0$. This implies that the kernel
of $(\pi_0)_*\G_{m,X_0}\to (\pi_1)_*\G_{m,X_1}$ is a torus. In
particular, the group of cocharacters of $T_1$ coincides with $\Im d^1$. Since
$\Ker d^2$ is the cocharacter group of $T_2$, we conclude that the
cocharacter group of $T$ is the
{\it free group} $\Ker d^2/\Im d^1$. }
\end{remark}

Let $(X_{\d}, Y_{\d})$ be a simplicial pair over $S_0$ such that
$X_{\d}$ is as above (\cf Lemma~\ref{reppic}) and $Y_{\d}\subseteq
X_{\d}$
is a reduced simplicial $S_0$-divisor with strict normal crossings.
Here we
have that every $Y_i$ is a divisor with strict normal
crossings in $X_i$, \ie $Y_i = \cup Y_{ij}$ such that $Y_{ij}$ is
smooth over $S_0$ and $Y_{ij}$ has
pure codimension $1$ in $X_i$.  For example, by de Jong's
theory of alterations \cite{DE}, every algebraic variety $V$ over
a perfect field $S_0 =\Spec(k)$ admits a proper hypercovering $V_{\d}\to
V$ such that $V_{\d}= X_{\d} - Y_{\d}$ for some $(X_{\d}, Y_{\d})$
as above, and moreover $V_0\to V$ is generically \'etale.
Consider the fppf-sheaf $\bPic_{[X_{\d} | Y_{\d}]/S_0}$
associated to the presheaf (\cf \cite[5.1.13.1]{Deligne})
$$T\leadsto  \HH^1_{Y_{\d}\times_{S_0} T}(X_{\d}\times_{S_0}
T,\G_{m}).$$
We obtain a canonical map (forgetting the support, \cf \cite[\S
4.2]{BS})
$$u \colon\bPic_{[X_{\d}| Y_{\d}]/S_0}\to \bPic_{X_{\d}/S_0}$$

Since $S_0= \Spec (k)$ and $k$ is a perfect field, we then
describe $u$ as follows.  Let $\bar{S_{0}}=
\Spec (\bar{k})$ be some algebraic closure of $k$.  Let $\Gal
(\bar{k}/k)$ be the Galois group.  Let $\bar{X_{\d}}$ and
$\bar{Y_{\d}}$ denote the base change to $\bar{S_{0}}$. Denote
\begin{equation}\label{Weil}
\bDiv_{\bar{Y_{\d}}}(\bar{X_{\d}})
\end{equation} the group of those Weil
divisors $D$ on $X_{0}\times_{S_0} \bar{S_{0}}$, supported on
$Y_{0}\times_{S_0} \bar{S_{0}}$, such that $d_0^*(D) = d^*_1(D)$
as divisors on $X_{1}\times_{S_0} \bar{S_{0}}$ (note that on $X_0$ and
$X_1$ Cartier and Weil divisors coincide).

The group $\bDiv_{\bar{Y_{\d}}}(\bar{X_{\d}})$ is naturally a
$\Gal (\bar{k}/k)$-module with respect to the Galois action on
Weil divisors.  Moreover, we have a natural identification of the
group of $\bar{S_{0}}$-points of $\bPic_{[\bar{X_{\d} } |\bar {
Y_{\d} } ]/\bar { S_0 } }$ with $\bDiv_{\bar{Y_{\d}}}(\bar{X_{\d}})$.
In fact, the analogous
spectral sequence (\ref{sscomp}) for the cohomology with supports
degenerates over $\bar{ S_0 }$ and yields the claimed group isomorphism
(see \cite[(33) p. 56]{BS}). We then obtain canonical $\Gal
(\bar{k}/k)$-equivariant maps
\begin{equation}\label{dpic}
\begin{CD} \bDiv_{\bar{Y_{\d}}}(\bar{X_{\d}})
@>{u^\prime}>>
\bPic(\bar{X_{\d}})\\
@V{\veq}VV@V{\veq}VV\\
\bPic_{[\bar{X_{\d} }|
\bar { Y_{\d} } ]/\bar { S_0 }}(\bar {S_0} ) @>>>
\bPic_{\bar{X_{\d} }/\bar { S_0 }}(\bar {S_0} ) \\
\end{CD}
\end{equation}
where $u^{\prime}$ is the canonical map sending a Weil divisor $D$
on $\bar{X_{0}}$ to the corresponding line bundle
$\cO_{\bar{X_{0}}}(D)$.
Let $\bPic_{[X_{\d} |Y_{\d}]/S_0}^{0, \red}$ be the inverse image of
$\bPic_{X_{\d}/S_0}^{0,\red}\subset \bPic_{X_{\d}/S_0}$ under the map
$u$.
Remark that $ \bPic_{\bar { X_{\d} }
/\bar { S_0 }}^{0, \red}$ is stable under the $\Gal (\bar{k}/k)$-action
on $ \bPic_{X_{\d}/S_0} \times_{S_0} \bar{S_{0 }}$.\\

The following definition is equivalent to \cite[Def. 4.2.1]{BS} and
\cite[Def. 4.1]{RA} but the proof of the independence of the choices
made, for positive characteristics, appears in our
Appendix~\ref{appendix} only.
\begin{definition} Let $(X_{\d}, Y_{\d})$ be a pair as above. Define
$$\Pic^{+}(X_{\d},
Y_{\d}) \df
[\bPic_{[X_{\d}| Y_{\d}]/S_0}^{0, \red}\by { u }
\bPic_{X_{\d}/S_0}^{0, \red}].$$ If $V_{\d}= X_{\d} - Y_{\d}$ is an
hypercovering of an algebraic variety $V$ over a perfect field (such that $V_0\to V$ is generically \'etale) we let  $$\Pic ^{+}(V)\df \Pic^{+}(X_{\d}, Y_{\d})$$
denote the {\rm cohomological Picard 1-motive}\, of $V$ (see the
Appendix~\ref{appendix} for the fact that $\Pic ^{+}$ is well defined
and contravariantly functorial).
\end{definition}

Conversely, we remark that the 1-motive $\Pic^{+}(X_{\d},Y_{\d})$ over
$S_{0}$ is determined by the $\Gal (\bar{k}/k)$-equivariant map
$$\bDiv_{\bar{Y_{\d}}} (\bar{X_{\d}})
\stackrel{u^{\prime}}{\llongrightarrow} \bPic(\bar{X_{\d}}) \cong
\bPic_{\bar{X_{\d}}/\bar{S_0}}(\bar{S_0}) $$over
$\bar{S_{0}}$, in particular $\Pic ^{+}(V)_{\bar{k}} =  \Pic
^{+}(V_{\bar{k}})$. We may and will denote $\bDiv_{\bar{Y_{\d}}}
^0(\bar{X_{\d}})$ the
subgroup of those divisors in $\bDiv_{\bar{Y_{\d}}} (\bar{X_{\d}})$
mapping to zero in
${\rm NS} (\bar{X_{\d}})\df \pi_0 (\bPic_{\bar{X_{\d}}/\bar{S_0}})$.\\

Now let $S_0\into S_n$ be a thickening defined by an ideal with
nilpotent divided powers. As above, we are mainly interested in
the case that the thickening $S_n = \Spec (\WW_{n+1}(k))$ is the
affine scheme defined by the Witt vectors of length $n+1$ (or
equivalently by $\WW (k)/p^{n+1}$). We then always get a lifting
$\Pic ^{+}(V)_n$ of the 1-motive $\Pic ^{+}(V)$ to $S_n$.

Remark that if $(X_{\d}, Y_{\d})$ itself lift (this is not the
case in general!) to a similar pair $(X_{\d}, Y_{\d})_n$ over
$S_n$ and $\bPic$ is representable over $S_n$, we may set
$$\Pic ^{+}(V)_n = \bPic ^{+}((X_{\d}, Y_{\d})_n).$$
Nevertheless, a formal lifting of the 1-motive $\Pic ^{+}(V)$ always exists.

\subsection{$\bPic^\natural$}\label{unipic}
Recall (see \cite[\S 4.5]{BS}) the construction of the functor
$\bPic_{X_{\d}/S_0}^{\natural}$. Here $S_0$ is the spectrum of a
perfect field and $X_{\d}$ is smooth and proper over $S_0$.

The group $\bPic^{\natural}(X_{\d})$ is given  by
isomorphism classes of pairs $(\ccL_{\d},\nabla_{\d})$, where
$\ccL_{\d}$ is a simplicial line bundle and $\nabla_{\d}$ is a
simplicial integrable connection
$$\nabla_{\d}\colon \ccL_{\d}\to \ccL_{\d}\otimes_{\cO_{X_{\d}}}
\Omega_{X_{\d}/S_0}^{1}.$$ There is a natural map
$$\bPic^{\natural}(X_{\d})\llongrightarrow
\HH^1(X_{\d},\cO^*_{X_{\d}}\by{\rm
dlog}\Omega^{1}_{X_{\d}/S_0}).$$ In characteristic zero this map
is an isomorphism and provides the group of $k$-points of the
universal extension of $\bPic^{0}$, whose Lie algebra is the first
De Rham cohomology (see \cite[4.5]{BS}), \ie the simplicial
$\natural$-Picard functor $\bPic_{X_{\d}/S_0}^{\natural}$,
obtained by sheafifying the functor
$$T\leadsto  \bPic^{\natural} (X_{\d}\times_{S_0} T)$$
with respect to the fppf-topology on $S_0$, yields an
exact sequence of Zariski (or \'etale) sheaves
\begin{equation}\label{natsimp}
0\to\pi_*\Omega^1_{X_{\d}/S_0} \to  \bPic^{\natural}_{X_{\d}/S_0}\to
\bPic_{X_{\d}/S_0}\to R^1\pi_*\Omega^1_{X_{\d}/S_0}
\end{equation}
and
\begin{equation}\label{derhamzero}
\Lie \bPic^{\natural , 0}_{X_{\d}/S_0} \cong \H^1_{\rm DR}
(X_{\d}/S_0)
\end{equation}
over $S_0$ of characteristic zero.

In positive characteristics one dimensional De Rham cohomology cannot be
recovered as above, \eg we only have an injection
$$\H^1_{\rm DR} (X_{\d}/S_0) \subseteq
\HH^1(X_{\d},\cO_{X_{\d}}\to\Omega^{1}_{X_{\d}/S_0})$$ determined
by the Oda subspace $\H^0(X_0, \Omega^{1}_{X_{0}/S_0})_{d=0}$ (\cf
\cite[p. 121]{Oda}).

However, for curves, $\Pic^{+}$ coincides with Deligne's definition of
motivic cohomology and its universal extension is given by the
$\natural$-Picard functor (see \cite[10.3.13]{Deligne} and \cf
\cite[5.2]{Oda}).
A similar case is given by $X_{\d}$ such that each $X_i$ has
connected components which are abelian schemes (\cf \cite[5.1]{Oda}).

Actually, whenever the Hodge-De Rham spectral sequence degenerates we
get
De Rham cohomology as above.
\begin{lemma}\label{repnat} The $\natural$-Picard sheaf
$\bPic^{\natural}_{X_{\d}/S_0}$
is representable, by a group scheme locally of finite type over
$S_0$. We have an extension
$$0\to\pi_*\Omega^1_{X_{\d}/S_0,d=0} \to  \bPic^{\natural,
0}_{X_{\d}/S_0}\to
\bPic_{X_{\d}/S_0}^{0, \red}\to 0$$
where we set $\Omega^*_{X_{\d}/S_0,d=0}$ according to Oda \cite[Def.
5.5]{Oda}.

If we further assume that the Hodge-De Rham
spectral sequence of the components $X_i$ degenerates then
(\ref{natsimp}) and
(\ref{derhamzero}) hold over $S_0$.
\end{lemma}
\begin{proof} Representability follows by a similar argument to the
proof of Lemma~\ref{reppic}, since $\Pic^{\natural}_{X_{i}/S_0}$
is representable, \eg by using \cite[Prop. 17.4]{OG}.

By applying the arguments in
\cite[4.1.2]{MazurMessing} and \cite[p. 62]{BS} to our setting we
can see that the indeterminacy in putting an integrable connection
on the trivial simplicial line bundle is the space of closed
simplicial 1-forms. Moreover, the obstruction map  in putting an
integrable connection on a simplicial line bundle defines  a
homomorphism from  $\bPic_{X_{\d}/S_0}$ to the vector group
$R^1\pi_*\Omega^1_{X_{\d}/S_0, d=0}$. Thus, when restricted to the
semiabelian scheme $\bPic^{0,\red}_{X_{\d}/S_0}$ it is trivial,
proving the exactness on the right of the sequence in the Lemma.
Note that here $\bPic^{\natural, 0}_{X_{\d}/S_0}$ is just the pullback
of $\bPic^{0,\red}_{X_{\d}/S_0}$ along $\bPic^{\natural}_{X_{\d}/S_0}
\to \bPic_{X_{\d}/S_0}$.

For the components $X_{i}$ such that all global 1-forms are closed
the arguments provided by \cite[4.5.1]{BS} apply; an easy spectral
sequence and Lie algebra computations yields the result (\cf
\cite[\S 4]{MazurMessing}).
\end{proof}

Note that if we consider the component $X_i$ of any smooth proper
simplicial scheme
$X_{\d}$ and let $A_i$ denote
$\Alb_{X_i/S_0} :=
(\Pic_{X_{i}/S_0}^0)^\vee$, the Albanese (abelian) scheme, we have that
the
following extension
$$0 \to (\pi_0)_* \Omega_{A_i/S_0}^1
\llongrightarrow\Pic^{\natural, 0}_{A_i/S_0} \llongrightarrow
\Pic_{X_{i}/S_0}^{0, \red}\to 0$$
is the universal extension of the abelian scheme $\Pic_{X_{i}/S_0}^{0,
\red}$.

\section{Log-crystalline cohomology}
We refer the reader to~\cite[\S5]{Kato} for the basics on
log-crystalline theory.

\subsection{Logarithmic structures}\label{Sdlog}
Let~$k$ be a perfect field of characteristic~$p$. For every~$n\in\NN$,
let~$\bigl(S_n,L_n,\gamma_n\bigr)$ be the triple consisting of
\begin{itemize}
\item[1)] $S_n:=\Spec\bigl(\WW_{n+1}(k)\bigr)$ the
spectrum of the Witt vectors of length~$n+1$ over~$k$;

\item[2)] $L_n$ the trivial  logarithmic structure on~$S_n$ defined
by the map (of multiplicative monoids) $\NN\mapsto \WW_{n+1}(k)$
given by $1 \mapsto 0$;

\item[\rm 3)] $\gamma_n$  the standard divided powers
structure on the ideal~$\cI_n:=(p)$ of\/~$\WW_n(k)$.
\end{itemize}
Assume these data to be compatible for varying~$n$.
As usual we view the category of schemes over~$S_n$ as the full
subcategory of log-schemes over~$S_n$ with logarithmic structure induced
by~$L_n$.

\begin{remark}\label{recallKato} {\rm (\cf \cite[\S1.5]{Kato})\enspace
Let~$X$ be a smooth projective scheme over~$S_0$. Let~$Y\subset X$
be a  divisor with normal crossings. Let~$M$ be the logarithmic
structure defined by~$Y$ on~$X$ as follows. We can cover~$X$ with
affine schemes~$\{U_i=\Spec(A_i)\}_i$, \'etale over $X$, such that
for each~$i$ there exists irreducible elements
$\{\pi_{i,\alpha}\}_{\alpha}\subset A_i$ such that~$Y\cap U_i$ is
defined by the ideal~$\prod_{\alpha}\pi_{i,\alpha}$. Then,
$$M_{U_i}=\left\{g\in A_i\bigl[\pi_{i,\alpha}^{\pm
1}\bigr]_{\alpha}\right\}$$and coincides with the logarithmic
structure associated to the pre-logarithmic structure
$$\prod_{\alpha} \NN \llongrightarrow A_i, \qquad
\bigl(a_{\alpha}\bigr)_{\alpha} \longmapsto \prod_{\alpha}
\pi_{i,\alpha}^{a_{\alpha}}.$$}
\end{remark}

\subsubsection{Log structures and crystals}\label{logcrysWeil}
We continue with the discussion in~\ref{recallKato}.
Let~$\bigl(U_i,T_i,M_{T_i},\iota,\delta\bigr)$ be an open set for
the logarithmic crystalline site of~$\bigl(X,M\bigr)$
over~$\bigl(S_n,L_n,\gamma_n\bigr)$ defined in~\ref{Sdlog}. Then,
$$\iota\colon \bigl(U_i,M_{U_i}\bigr) \hookrightarrow
\bigl(T_i,M_{T_i}\bigr)$$is an exact closed immersion of
logarithmic schemes defined by an ideal~$I_i$ with divided
powers~$\delta$; \cite[\S5.2]{Kato}. In particular, $I_i$ is a
nilpotent ideal and the log-structure on $U_i$ is the one induced
from $T_i$. Therefore, for every~$\alpha$ there
exists~$m_{i,\alpha}\in\Gamma\bigl(T_i,\cO_{T_i}\bigr)$ such
that~$\iota^*\bigl(m_{i,\alpha}\bigr)\in A_i^* \cdot
\pi_{i,\alpha}$. Any two elements~$m_{i,\alpha}$
and~$m_{i,\alpha}'$, having the same property, satisfy
$$m_{i,\alpha}'m_{i,\alpha}^{-1}\in\Gamma\bigl(T_i,
\cO_{T_i}^*\bigr).$$Hence, the ideal~$\bigl(m_{i,\alpha}\bigr)$
defines a Cartier divisor on~$T_i$ relative to~$S_n$ lifting the
Cartier divisor~$\bigl(\pi_{i,\alpha}\bigr)$ on~$U_i$. Hence, for
every element~$\pi$ in the group~$M^{\rm gp}$ associated to~$M$ we
get a well defined invertible sheaf
$$\cO_{T_i}\bigl[\pi^{-1}\bigr]$$ on~$T_i$ lifting the invertible
sheaf on~$U_i$ defined by~$\pi$.

Let $(X_{\d}, Y_{\d})$ be a simplicial pair over~$S_0$ such
that~$X_{\d}$
is a projective and smooth simplicial scheme over~$S_0$
and~$Y_{\d}\subseteq
X_{\d}$ is a simplicial divisor with normal crossings relative
to~$S_0$.
\smallskip
Identify the pair~$(X_i,Y_i)$ with the scheme~$X_i$ with the fine
logarithmic structure defined by~$Y_i$. Then, $(X_i,Y_i)$ lies
over~$(S_0,L_0)$ in the sense of log schemes and it is log-smooth
over~$(S_0,L_0)$~{\rm \cite[Ex 3.7]{Kato}}.
\begin{definition}\label{Xdlog}
A {\em simplicial logarithmic pair} over~$(S_0,L_0)$ is a
simplicial pair $(X_{\d}, Y_{\d})$ over~$S_0$ as above where
$(X_i,Y_i)$ is regarded as a scheme~$X_i$ with the fine
logarithmic structure defined by~$Y_i$.
\end{definition}

\subsection{$\check{\bf C}$ech coverings}\label{Chechnotation}
Let~$(X,M)$ be a logarithmic scheme log-smooth over~$(S_0,L_0)$.
Assume that~$X$ is projective and smooth (in the classical sense)
over~$S_0$. Let~$\{U_i\}_i$ be a covering family of~$X$ by affine
schemes \'etale over~$X$. We choose and fix a total order on the
set of indices~$\{i\}$. Each~$U_i$ lifts to an
open~$U_i\hookrightarrow V_i$ in the logarithmic  site of~$(X,M)$
over~$(S_n,L_n)$ so that~$V_i$ is smooth over~$S_n$. In
particular, $V_i$ is flat over~$S_n$ and~$\cI_n V_i$ has a unique
PD structure extending the one on~$\cI_n$.  Then, $\bigl\{U_i
\hookrightarrow V_i\bigr\}_i$ is a covering family in the
site~$\bigl((X,M)/(S_n,L_n)\bigr)^{\logcrys}$.   For any~$i<j$
let~$U_{i,j}:=U_i\cap U_j$. Since~$X$ is separated, $U_{i,j}$ is
an affine scheme. Let~$V^j_i$ (resp.~$V^i_j$) be the open
subscheme of~$V_i$ (resp.~$V_j$) defined by~$U_{i,j}$. They are
isomorphic. Define~$V_{i,j} :=V^j_i$. Then~$V_{i,j}$ is the
logarithmic divided power envelope of the thickening~$U_{i,j}
\subset V_{i,j}$~\cite[Def 5.4 ]{Kato}. Continuing in this fashion
one defines
$$U_{\underline{i}} \hookrightarrow  V_{\underline{i}} $$for
every $n+1$-uple of indices $\underline{i}=(i_0,\ldots,i_n)$. One
gets an hypercovering
$$\cdots \amalg_{i_0< i_1} U_{(i_0, i_1)} \rightrightarrows \amalg
U_{i}$$and a Leray spectral sequence
$$E_2^{s,t}=\H^t_\logcrys\bigl(\amalg_{i_0< \cdots < i_s}
U_{(i_0, \ldots, i_s)} /(S_n,L_n)\bigr) \Rightarrow \H^{s+t}_\logcrys
\bigl((X,M)/(S_n,L_n)\bigr).$$By~\cite[Thm 6.4]{Kato} we have
$$\H^t_\logcrys\bigl(U_{(i_0, \ldots, i_s)}/(S_n,L_n)\bigr)
\cong\H^t\bigl(V_{(i_0, \ldots, i_s)},\omega^{\star,\log}_{V_{(i_0,
\ldots, i_s)}/(S_n,L_n)}\bigr);$$where $\omega^{\star,\log}$ is the de Rham
complex of differentials defined in the logarithmic sense~\cite[\S
1.7]{Kato}.
\begin{lemma}\label{Chechcomp} The homology of the total complex
of the bicomplex $\oplus_{i_0< \cdots< i_s}\omega^{t,\log}_{V_{(i_0,
\ldots, i_s)}/(S_n,L_n)}$ is $\H^{s+t}_\logcrys
\bigl((X,M)/(S_n,L_n)\bigr)$.
\end{lemma}

\subsection{Crystalline versus log-crystalline}\label{fil2}
The notation is as in~\ref{Chechnotation}. Let~$(X,L_0)$ be the
scheme~$X$ with logarithmic structure defined by pull back
of~$L_0$. The morphism~$f\colon (X,M)\rightarrow (X,L_0)$ induces
a natural morphism of topoi~{\rm \cite[\S5.9]{Kato}}
$$f_\logcrys \colon \bigl((X,M)/(S_n,L_n)\bigr)^{\logcrys}
\llongrightarrow \bigl((X,L_0)/(S_n,L_n)\bigr)^{\logcrys}=
\bigl(X/S_n\bigr)^\crys.$$Let~$I$ be a sheaf on~$(X/S_n)^\crys$.
Let~$(U,M_U)\hookrightarrow (T,M_T)$ be an open
on~$(X,M)^\logcrys$ defined by a PD ideal,
then~$f_\logcrys^*(I)\bigl((U,M_U)\subset
(T,M_T)\bigr)=I\bigl(U\subset T\bigr)$. This implies
that~$f_\logcrys^*$ is exact and commutes with taking global
sections. In particular, let $\cO_{X}^\crys \rightarrow I^{\d}$
and~$\cO_{(X,M)}^\logcrys \rightarrow J^{\d}$ be injective
resolutions. By the exactness of~$f_\logcrys^*$ we get a
commutative diagram with exact rows
\begin{displaymath}
\begin{CD} 0 @>>> f_\logcrys^*\bigl(\cO_{X}^\crys \bigr)
@>>> f_\logcrys^*\bigl(I^{\d}\bigr) \\
&& @VVV  @VVV \\0  @>>> \cO_{(X,M)}^\logcrys @>>> J^{\d}.
\end{CD}
\end{displaymath}Taking global sections we obtain for
every~$i\in\NN$ homomorphisms $$\H^i_{\crys}\bigl(X/S_n\bigr)
\longrightarrow \H^i_{\logcrys}\bigl((X,M)/(S_n,L_n)\bigr).$$We
wish to determine kernel and cokernel of the homomorphism
for~$i=0$ and~$1$. To do that we use the spectral sequence of the
`$\check{\rm C}$ech cohomology' for the log-crystalline as in
Lemma~\ref{Chechcomp}. Since we are considering and comparing two
logarithmic structures on~$X$, we write a superscript~$\log$
whenever we consider~$(X,M)$ and we use no superscript whenever we
consider~$(X,L_0)$. Fix indices~$i_0<\cdots<i_s$.
Let~$\alpha\colon M_{V_{(i_0,\ldots, i_s)}} \rightarrow
\cO_{V_{(i_0,\ldots, i_s)}}$ be the logarithmic structure
on~$V_{(i_0,\ldots, i_s)}$.  For every $r$ and $n$ let
$$\Q_{(i_0,\ldots,
i_s)}^r:=\omega^{r,\log}_{V_{(i_0,\ldots,
i_s)}/(S_n,L_n)}/\omega^r_{(i_0,\ldots, i_s)/S_n}$$be the quotient of
the $r$-th wedge product of the logarithmic differentials
of~$\bigl(V_{(i_0,\ldots, i_s)},M_{V_{(i_0,\ldots, i_s)}}\bigr)$
relative to~$(S_n,L_n)$ by the $r$-th wedge product of the usual
K\"ahler differentials of~$V_{(i_0,\ldots, i_s)}$ relative
to~$S_n$. By possibly shrinking the open subschemes~$\{U_i\}$ and
by the assumptions of logarithmic smoothness, we may assume that
there exists~$\{v_1,\ldots,v_w\}$ global sections
of~$\cO_{(i_0,\ldots, i_s)}$ and~$\{m_1,\ldots,m_z\}$ global
sections of~$M_{V_{(i_0,\ldots, i_s)}}$ such that
$$\omega^{1,\log}_{V_{(i_0,\ldots,
i_s)}/(S_n,L_n)}=\cO_{V_{(i_0,\ldots, i_s)}}dv_1\oplus\cdots\oplus
\cO_{V_{(i_0,\ldots, i_s)}}dv_w\oplus \cO_{V_{(i_0,\ldots,
i_s)}}d\log\bigl(m_1\bigr)\oplus\cdots\oplus d\log\bigl(m_z
\bigr).$$In particular, $$\Q^1_{(i_0,\ldots, i_s)}\cong
\cO_{V_{(i_0,\ldots, i_s)}}/m_1\cO_{V_{(i_0,\ldots, i_s)}}\oplus
\cdots \oplus \cO_{V_{(i_0,\ldots, i_s)}}/ m_z\cO_{V_{(i_0,\ldots,
i_s)}}.$$Consider the following diagram
\begin{displaymath}
\begin{CD}
0 \to \oplus_i \cO_{V_i} @>>> \oplus_{i,j}\cO_{V_{i,j}} \oplus_i
\omega^1_{V_i/S_n} @>>> \oplus_{i,j,k} \cO_{V_{i,j,k}}\oplus_{i,j}
\omega^2_{V_{i,j}/S_n}\oplus_{i,j} \omega^1_{V_{i,j}/S_n}\\
\veq && @VVV @VVV  \\
0 \to \oplus_i \cO_{V_i} @>>> \oplus_{i,j}\cO_{V_{i,j}} \oplus_i
\omega^{1,\log}_{V_i/(S_n,L_n)} @>>> \oplus_{i,j,k}
\cO_{V_{i,j,k}}\oplus_{i,j}
\omega^{2,\log}_{V_{i,j}/(S_n,L_n)}\oplus_{i,j}
\omega^{1,\log}_{V_{i,j}/(S_n,L_n)}\\
& & @VVV @VVV  \\
&& 0 \to \oplus_i \Q_i^1 @>>> \oplus_{i,j} \Q_{i,j}^1 \oplus_i
\Q_i^2\\
\end{CD}
\end{displaymath}By Lemma~\ref{Chechcomp} the homology of the first row
computes $\H^*_\crys(X/S_n)$, while the homology of the second row
computes $\H^*_\logcrys\bigl((X,M)/(S_n,L_n)\bigr)$.

\begin{corollary}\label{HH0logcrys}
The group $\H^0_{\logcrys}\bigl((X,M)/(S_n,L_n))$ is equal to
$\H^0_{\crys}\bigl(X/S_n)$ and is a locally
free~$\cO_{S_n}$-module of rank equal to the number of the
geometric irreducible components of~$X \rightarrow S_0$.
\end{corollary}
\begin{proof}  By~\ref{fil2} the formation
of~$\H^0_{\crys}\bigl(X/S_n\bigr)$ and
of~$\H^0_{\logcrys}\bigl((X,M)/(S_n,L_n)\bigr)$ commute with field
extensions of~$S_0$. Hence, we may assume that~$X\rightarrow S_0$
is geometrically irreducible. Furthermore, both groups injects in
$$K_n:=\Ker\Bigl(\prod_i \Gamma\bigl(V_i,\cO_{V_i}\bigr) \rightarrow
\prod_{i<j} \Gamma\bigl(V_{i,j},\cO_{V_{i,j}}\bigr)\Bigr)$$and we
have inclusions~$\cO_{S_n}\hookrightarrow
\H^0_{\crys}\bigl(X/S_n\bigr)\hookrightarrow
\H^0_{\logcrys}\bigl((X,M)/S_n\bigr)$. It suffices to prove that
the composite $\cO_{S_n}\hookrightarrow K_n$ is an isomorphism.
Proceeding by induction on~$n$ and since for every~$n$ we have the
exact sequence $0\rightarrow K_0 \stackrel{p^n}{\longrightarrow}
K_{n+1} \longrightarrow K_n$, one concludes.
\end{proof}

\begin{proposition}\label{compaeH1crysandH1LoGcrys}
Let $M$ denote the log-structure defined by a strict normal crossing
divisor $Y$ in $X$ smooth and projective over $S_0$.
We have an exact sequence
$$0 \to \H^1_{\crys}(X) \to
\H^1_{\logcrys}(X, M) \to \Div_Y (X)\otimes_\ZZ \WW(k)
\to \H_\crys^2(X).$$The map $\Div_Y (X)\otimes_\ZZ \WW(k) \to
\H_\crys^2(X)$ is defined via the crystalline first Chern class
and its kernel coincides with $\Div_Y^0 (X)\otimes_\ZZ \WW(k)$.
\end{proposition}
\begin{proof} We let $\Q_{\d}^r$ denote the complex
$0 \to \oplus_i \Q_i^r \to \oplus_{i,j} \Q_{i,j}^r\to\cdots $
By the arguments in~\ref{fil2}, for every $n\in\NN$, we have an exact
sequence
\begin{equation}\label{algeqtozeroviaChern}
0 \to \H^1_{\crys}\bigl(X/S_n\bigr) \to
\H^1_{\logcrys}\bigl((X,M)/(S_n,L_n)\bigr) \to
\Ker\Bigl(\H^0(\Q_{\d}^1)\to \H^0(\Q_{\d}^2)\Bigr)  \to
\H_\crys^2(X/S_n).
\end{equation} The claimed exact sequence follows inspecting such exact
sequence
and taking limits over $n\in\NN$. Note that all groups appearing above
are
finitely generated $\WW_{n+1}(k)$-modules and, thus, satisfy the
Mittag-Leffler condition. We further may assume that $k=\bar{k}$ is
algebraically closed. For $j=1,\ldots,z$ the image of $d\log(m_j)$
via the derivation $\omega^{1,\log}_{V_i/(S_n,L_n)} \to
\omega^{2,\log}_{V_i/(S_n,L_n)}$ is $d\bigl( d\log (m_j)\bigr)=0$.
Hence, the map $\H^0(\Q_{\d}^1) \to \oplus_i \Q_i^2$ is zero and $
\Ker\left(\H^0(\Q_{\d}^1)\to \H^0(\Q_{\d}^2)\right)=
\H^0(\Q_{\d}^1)$. Since the log-structure $M$ is defined by a
strict normal crossing divisor $Y$ of $X$, if $\{Y_j\}$ is the set
of irreducible components of $Y$, it follows from
Corollary~\ref{HH0logcrys}
that $\H^0(\Q_{\d}^1)=\prod_j
\H^0_\crys(Y_j/S_n)=\Div_Y (X)\otimes_\ZZ \WW_{n+1}(k)$.
The first Chern class map $$c_1^\crys\colon
\H^1(X,\cO_X^*) \to \H_\crys^2(X/S_n)\qquad \Bigl(\hbox{{\rm resp. }}
c_1^\logcrys\colon \H^1(X,\cO_X^*) \to
\H_\logcrys^2(X/(S_n,L_n))\Bigr)$$is
defined taking the long exact sequence of cohomology groups
associated to the short exact sequence $0 \rightarrow 1+ p
\cO^\crys_{X/S_n} \rightarrow {\cO^\crys_{X/S_n}}^* \rightarrow
\cO^*_{X/S_0} \rightarrow 0$ and the map $\log\colon 1+ p
\cO^\crys_{X/S_n} \rightarrow \cO^\crys_{X/S_n}$ (analogously for
$c_1^\logcrys$). Comparing~$c_1^\crys$ and~$c_1^\logcrys$ via the
long exact sequence relating the crystalline and the
log--crystalline cohomology groups, one sees that the composite of
the projection $\H^1_{\logcrys}(X/S_n,{\cO^\logcrys_{X/S_n}}^*)
\to \H^1(X,\cO_X^*)$ and $c_1^\crys$ factors via the connecting
homomorphism $\H^0(Q_{\d}^1)  \to \H_\crys^2(X/S_n)$
of~\eqref{algeqtozeroviaChern}. As remarked in~\ref{logcrysWeil}
$\Div_Y(X)\to \H^1(X,\cO_X^*)$ factors via
$\H^1_{\logcrys}(X/S_n,{\cO^\logcrys_{X/S_n}}^*)$. Thus,
$c_1^\crys$ restricted to $\Div_Y(X)$ factors as  $\Div_Y(X) \to
\H^0(Q_{\d}^1) \to
\H_\crys^2(X/S_n)$.

By~\cite[II.6.8.0]{Illusie} and the comparison between crystalline
and de Rham-Witt cohomology for proper and smooth schemes over
perfect fields, the Chern class $c_1^\crys$ induces an {\it
injection} ${\rm NS}(X)\otimes_\ZZ \ZZ_p \into \H_\crys^2(X)$.
By~\cite[II.5.10.1]{Illusie} and~\cite[II.6.8.4]{Illusie} the map
${\rm NS}(X)\otimes_\ZZ \WW(k) \to \H_\crys^2(X)$ is injective as
well.  Since $\Div^0_Y (X)$ is the subgroup of $\Div_Y (X)$ of
divisors {\it algebraically equivalent to zero}, the conclusion
follows.\end{proof}

\subsection{Crystalline via simplicial log-crystalline}
\label{H1crys1}
Let $(X_{\d}, Y_{\d})$ be a simplicial pair and
$n\in\NN$. For every sheaf~$\cF_{\d}$ of
abelian groups over the simplicial logarithmic crystalline site
$\bigl(X_{\d}, Y_{\d}\bigr)/\bigl(S_n,L_n,\gamma_n\bigr)$ define
$$\HH^*_\logcrys\Bigl(\bigl(X_{\d},
Y_{\d}\bigr)/\bigl(S_n,L_n,\gamma_n\bigr), \cF_{\d}\Bigr)$$ as the
right derived functor of the following left exact functor\/
$$\cF_{\d} \longmapsto \Ker\left(\Gamma\Bigl(\bigl(X_0,
Y_0\bigr)/\bigl(S_n,L_n,\gamma_n\bigr),\cF_0\Bigr)
\stackrel{(d_0^1)^*-(d_1^1)^*}{\llongrightarrow}
\Gamma\Bigl(\bigl(X_1,
Y_1\bigr)/\bigl(S_n,L_n,\gamma_n\bigr),\cF_1\Bigr)\right).$$ Note
that each $\cF_i$ is a sheaf on the scheme~$X_i$ with the
logarithmic structure. Define\/
$$\HH^i_\logcrys\Bigl(\bigl(X_{\d},
Y_{\d}\bigr)/\bigl(S_n,L_n,\gamma_n\bigr)\Bigr)$$as
$\HH^i_\logcrys$ of the structure sheaf $\cO_{X_{\d}}$. We
write~$\HH^i_\logcrys\bigl((X_{\d}, Y_{\d})/S_n\bigr)$ whenever
the logarithmic structure and the divided powers
structures~$\bigl(L_n,\gamma_n\bigr)$ are the ones fixed in
\ref{Sdlog}. Define
$$\HH^i_\logcrys\bigl(X_{\d}, Y_{\d}\bigr):=\lim_{\infty\leftarrow n}
\HH^i_\logcrys\bigl((X_{\d}, Y_{\d})/S_n\bigr).$$
Denote $\H^i_\crys(V_{\d}/\WW (k))\df\HH^i_\logcrys\bigl(X_{\d},
Y_{\d}\bigr)$
where $V_{\d} = X_{\d} - Y_{\d}$ is the corresponding smooth simplicial
scheme.
\subsection{A spectral sequence, $\H^1_\logcrys$ and
$\H^0_\logcrys$}\label{fil1}

   From the spectral sequence
\begin{equation}\label{cryss}
{\rm E}_1^{p,q}=\H^q_\logcrys\bigl((X_p,Y_p)/S_n\bigr) \implies
\HH^{p+q}_\logcrys\bigl((X_{\d}, Y_{\d})/S_n\bigr).
\end{equation}
we deduce exact sequences describing~$\H^1_\crys$
and~$\H^0_\crys$:
\begin{multline}\label{spectre} 0\llongrightarrow
\frac{\Ker\Bigl(\H^0_{\logcrys} \bigl((X_1,Y_1)/S_n\bigr)
\rightarrow
\H^0_{\logcrys}\bigl((X_2,Y_2)/S_n\bigr)\Bigr)}{\Im\Bigl(\H^0_{\logcrys}
\bigl((X_0,Y_0)/S_n\bigr) \rightarrow
\H^0_{\logcrys}\bigl((X_1,Y_1)/S_n\bigr)\Bigr)} \llongrightarrow
\HH^1_{\logcrys}\bigl((X_{\d},Y_{\d})/S_n\bigr) \\
\llongrightarrow
\Ker\Bigl(\H^1_{\logcrys}\bigl((X_0,Y_0)/S_n\bigr) \rightarrow
\H^1_{\logcrys}\bigl((X_1,Y_1)/S_n\bigr)\Bigr)
\llongrightarrow \\
\frac{\Ker\Bigl(\H^0_{\logcrys} \bigl((X_2,Y_2)/S_n\bigr)
\rightarrow
\H^0_{\logcrys}\bigl((X_3,Y_3)/S_n\bigr)\Bigr)}{\Im\Bigl(\H^0_{\logcrys}
\bigl((X_1,Y_1)/S_n\bigr) \rightarrow
\H^0_{\logcrys}\bigl((X_2,Y_2)/S_n\bigr)\Bigr)}.
\end{multline}
The arrows are defined by alternating sums of pull-backs along the
faces of the simplicial structure (\cf \cite[\S 4]{BS}).
By Corollary~\ref{HH0logcrys} we may replace~$\H^0_{\logcrys}$
with~$\H^0_{\crys}$ everywhere.  Analogously, we have a spectral
sequence giving us $\H^0_\logcrys$ as
$$\HH^0_\logcrys\bigl(X_{\d}/S_n\bigr) = \Ker
\Bigl(\H^0_{\logcrys} \bigl(X_0/S_n\bigr) \rightarrow
\H^0_{\logcrys}\bigl(X_1/S_n\bigr)\Bigr).$$
In particular, $\HH^0_\logcrys\bigl((X_{\d},
Y_{\d})/S_n\bigr) = \HH^0_\crys\bigl(X_{\d}/S_n\bigr)$. Taking limits
over $n\in\NN$ one concludes that
$$\frac{\Ker\Bigl(\H^0_{\logcrys}
\bigl(X_1,Y_1\bigr) \rightarrow
\H^0_{\logcrys}\bigl(X_2,Y_2\bigr)\Bigr)}{\Im\Bigl(\H^0_{\logcrys}
\bigl(X_0,Y_0\bigr) \rightarrow
\H^0_{\logcrys}\bigl(X_1,Y_1\bigr)\Bigr)}$$ is a {\it free}
$\WW(k)$-module. In fact, in the notation of
Remark~\ref{withouttorsion},
the latter is identified with $(\Ker d^2/\Im d^1)\otimes \WW (k)$ and
$\Ker d^2/\Im d^1$ is a free abelian group.

\subsection{Properties of $\H^1_{\logcrys}\bigl(X_{\d},Y_{\d}
\bigr)$}\label{verygory} The $\WW_{n+1}(k)$-modules appearing
in~\ref{fil2}-\ref{fil1} are finitely generated and, thus, they
satisfy the Mittag-Leffler condition. In particular, we get a
canonical (decreasing) weight filtration~$W$ on our crystalline cohomology
$\H^1_\crys(V_{\d}/\WW (k)):=\HH^1_\logcrys \bigl(X_{\d},
Y_{\d}\bigr)$ as follows:
\begin{itemize}
\item $W_{\geq 2} :=\HH^1_\logcrys
\bigl(X_{\d},Y_{\d}\bigr)$;
\item $W_{1}:=\HH^1_\crys (X_{\d})$ which is a $\WW(k)$-submodule of
$\HH^1_\logcrys \bigl(X_{\d},Y_{\d}\bigr)$
by \ref{fil1}, \ref{HH0logcrys} and~\eqref{algeqtozeroviaChern};
\item $W_{0} :=\frac{\Ker\Bigl(\H^0_{\crys}
\bigl(X_1\bigr) \rightarrow
\H^0_{\crys}\bigl(X_2\bigr)\Bigr)}{\Im\Bigl(\H^0_{\crys}
\bigl(X_0\bigr) \rightarrow \H^0_{\crys}\bigl(X_1\bigr)\Bigr)}$
which is a free $\WW(k)$-submodule of $\HH^1_\crys (X_{\d})$
by~\eqref{spectre} and the discussion in~\ref{fil1};
\item $W_{<0}:=0$.
\end{itemize}

\noindent Moreover we have the following description of the graded pieces $\Gr_i^W$
\begin{itemize}
\item  $\Gr_2^W$ over $k = \bar k$, by \eqref{spectre}, \ref{HH0logcrys} and~\ref{compaeH1crysandH1LoGcrys}, is contained in the $\WW(k)$-module $\bDiv_{Y_{\d}}^0 (X_{\d}) \otimes_\ZZ
\WW(k)$: here $\bDiv_{Y_{\d}}^0 (X_{\d})$ are the divisors defined in \ref{simpic}.
Therefore, $\Gr_2^W$ is a finitely generated free $\WW(k)$-module;

\item $\Gr_1^W$ is
contained in $\Ker\Bigl(\H^1_{\crys}\bigl(X_0\bigr)\rightarrow
\H^1_{\crys}\bigl(X_1\bigr) \Bigr)$. In particular, $\Gr_1^W$  is
a finitely generated, free $\WW(k)$-module by \cite[Prop.
II.3.11]{Illusie}.
\end{itemize}

\noindent For $p>2$ it follows from the proof of Theorem~B${}^{\prime}$,
see \ref{isomorphicrys} and \ref{compactinjective}, that in fact
$$\Gr_1^W=\Ker\Bigl(\H^1_{\crys}\bigl(X_0\bigr)\rightarrow
\H^1_{\crys}\bigl(X_1\bigr) \Bigr)$$
and, by \ref{endproof}, also $\Gr_2^W=\bDiv_{Y_{\d}}^0 (X_{\d}) \otimes_\ZZ
\WW(k)$. Anyway, we can already
conclude that $\H^1_\crys(V_{\d}/\WW (k))$ is a finitely
generated free $\WW(k)$-module.\\

Next, we show that $\H^1_\crys(V_{\d}/\WW (k))$ is naturally
endowed with the structure of a filtered $F$-$\WW (k)$-module. The
filtration has been discussed above. We define Frobenius. Let
$(X_{\d}^{(p)},Y_{\d}^{(p)})$ be the simplicial scheme obtained
pulling back $(X_{\d},Y_{\d})$ via the Frobenius map $\Spec(k)\to
\Spec(k)$. Then, Frobenius on $(X_{\d},Y_{\d})$ defines a morphism
of  simplicial $k$-schemes $ (X_{\d},Y_{\d}) \to
(X_{\d}^{(p)},Y_{\d}^{(p)})$ and, hence, a $\WW(k)$-linear map
$$ \HH^1_\logcrys \bigl(X_{\d}^{(p)},Y_{\d}^{(p)}\bigr)
\cong  \HH^1_\logcrys \bigl(X_{\d}, Y_{\d}\bigr)\otimes_\sigma
\WW(k) \llongrightarrow \HH^1_\logcrys \bigl(X_{\d},
Y_{\d}\bigr)$$ ($\sigma=$Frobenius on $\WW(k)$), which respects
the weights filtration by functoriality. In conclusion, we obtain
a $\sigma$-linear map
$$F\colon \HH^1_\logcrys \bigl(X_{\d},
Y_{\d}\bigr)\llongrightarrow \HH^1_\logcrys \bigl(X_{\d},
Y_{\d}\bigr).$$

\section{Proof of Theorem B${}^\prime$}\label{proof}
Our main task here is to compute, given a variety $V$ over a
perfect field $k$ of characteristic~$p \geq 3$, the Lie algebra of
the universal vector extension of a lifting of $\Pic^{+}(V)$ to
the Witt vector. See~\ref{casep=2} for some remarks on the
characteristic~$2$ case.
\subsection{Intermediate functors}
First construct some simplicial variations on themes from
\cite{MazurMessing} and log-crystalline variations on themes from
\cite{BS}.

\subsubsection{The functor~$\bLogCrys$}\label{functorLogGrys}
Assume that $p\geq 3$. Let $(X_{\d},Y_{\d})$ be a simplicial
logarithmic pair, see \ref{Sdlog}. Thus\/~$X_{\d}$ is a simplicial
log-scheme over~$S_0$. Let
$$\bLogCrys_{(X_{\d},Y_{\d})/S_n}\subset
\bLogCrysno_{(X_{\d},Y_{\d})/S_n}$$be the sheaves over the fppf
site of~$S_n$ associated to the following presheaves. Let~$T$ be a
scheme flat and of finite presentation over~$S_n$. Consider the
logarithmic and divided powers structure on~$T$ induced from those
on~$S_n$. Let~$\bLogCrysno\bigl((X_{\d},Y_{\d})/T\bigr)$ be the
group of isomorphism classes of simplicial crystals of invertible
$\cO^\logcrys_{X_{\d}\times_{S_n} T/T}$-modules for the {\it
nilpotent} crystalline site of~$X_{\d}\times_{S_n} T$ relative
to~$ T$. There is a canonical pull back to the Zariski
site~$(X_{\d}\times_{S_n} T)_{\rm Zar}$. Let
$\bLogCrys\bigl((X_{\d},Y_{\d})/T\bigr)$ denote the subgroup of
sections of~$\bLogCrysno\bigl((X_{\d},Y_{\d})/T\bigr)$ which land
in $\bPic^{0, \red}_{X_{\d}/S_0}\bigl(T\times_{S_n} S_0\bigr)$
(\cf \ref{simpic}).

\noindent We might describe an element
of~$\bLogCrys\bigl((X_{\d},Y_{\d})/T\bigr)$ as the giving of a
pair $(L, \alpha )$
\begin{itemize}
\item[\rm i)] $L$ yields an element
of~$\LogCrys\bigl((X_0,Y_0)/T\bigr)$;
\item[\rm ii)] and the isomorphism $\alpha\colon
(d_0^1)^*(L)\isomarrow (d_1^1)^*(L)$ as elements
of~$\LogCrys\bigl((X_1,Y_1)/T\bigr)$ satisfying the cocycle
condition
$$\bigl((d_1^2)^*(\alpha)\bigr)^{-1}\circ
\bigl((d_2^2)^*(\alpha)\bigr)\circ\bigl((d_0^2)^*(\alpha)
\bigr)={\rm Id}.$$
\end{itemize}
Analogous description holds
for~$\bLogCrysno\bigl((X_{\d},Y_{\d})/T\bigr)$. For the sake of
notation, in the following, we sometimes omit reference to the
pair $(X_{\d},Y_{\d})$ if it is clear from the context.
\begin{lemma}\label{categorygroups} We have an exact sequence
of functors, respecting the group laws

\begin{multline*}0 \llongrightarrow
\frac{\Ker^\dag\Bigl(\pi_{1,*}^\logcrys \bigl(\GG_{m,X_1}\bigr)
\rightarrow \pi_{2,*}^\logcrys
\bigl(\GG_{m,X_2}\bigr)\Bigr)}{\Im\Bigl(\pi_{0,*}^\logcrys
\bigl(\GG_{m,X_0}\bigr) \rightarrow \pi_{1,*}^\logcrys
\bigl(\GG_{m,X_1}\bigr)\Bigr)  } \llongrightarrow
\bLogCrys_{(X_{\d},Y_{\d})/S_n} \\ \llongrightarrow
\Ker\Bigl(\LogCrys_{(X_0,Y_0)/S_n}\rightarrow
\LogCrys_{(X_1,Y_1)/S_n}\Bigr)\end{multline*} where $\Ker^\dag (-
)$ denotes the subsheaf of $\Ker (-)$ of those elements which land
in $\Ker^{0,\red} (-)$ over~$S_0$.
\end{lemma}\begin{proof} Here, for every
scheme~$T$ flat over~$S_n$ and any element~$(L,\alpha)$
of~$\bLogCrys\bigl((X_{\d},Y_{\d})/T\bigr)$, the canonical map is
$\bigl(L,\alpha\bigr)\mapsto L$. The kernel consists of the
automorphisms of~$(d_0^1)^*(\cO^\logcrys_{X_0\times T})\isomarrow
(d_1^1)^*\cO^\logcrys_{X_0\times T})$ satisfying the usual cocycle
condition  modulo the automorphisms of~$\cO^\logcrys_{X_0\times
T}$. This description proves the lemma with the correction $\Ker^\dag (-
)$ according to the description of the toric part of $\bPic^{0, \red}_{X_{\d}/S_0}$ in Remark~\ref{withouttorsion}.
\end{proof}

\subsubsection{Dual numbers}
Let~$Z$ be a log-scheme and $(I_Z,\gamma)$ a nilpotent ideal with
divided powers structure. Let~$Z[\varepsilon]:=\Spec\bigl(
\cO_Z[\varepsilon]/(\varepsilon^2)\bigr)$ be the scheme of dual
numbers over~$Z$. Then, $I_Z\cO_Z[\varepsilon]/(\varepsilon^2)$ is
endowed with a unique divided powers structure extending the one
on~$I_Z$; \cf \cite[Cor. 3.22]{BerthelotOgus}. Let~$L_n$ be the
unique  fine logarithmic structure extending the one on~$Z$. Let
$$i_Z\colon Z \rightarrow Z[\varepsilon]\qquad\hbox{{\rm
and}}\qquad j_Z\colon Z[\varepsilon] \rightarrow Z$$be the closed
immersion defined by~$\varepsilon=0$ and, respectively, the
natural projection map. These maps are compatible with logarithmic
and divided powers structures.

\subsubsection{The functor~$\InfDef$}\label{ccoommppaarree}
Let $$\InfDef\Bigl(\cO^\logcrys_{X_{\d}/S_n}\Bigr)$$be the group
of infinitesimal deformations of the structural
sheaf~$\cO^\logcrys_{X_{\d}/S_n}$ or, equivalently, the group of
isomorphism classes of pairs~$(\cM,\tau_\cM)$ where
\begin{itemize}
\item[\rm i)] $\cM$ is a
crystal of invertible
$\cO^\logcrys_{X_{\d}\times_{S_n}S_n[\varepsilon]/S_n[\varepsilon]}$-
modules,
\item[\rm ii)] $\tau_\cM\colon
i_{X_{\d},\logcrys}^*(\cM)\isomarrow \cO^\logcrys_{X_{\d}/S_n}$ is
an isomorphism of $\cO^\logcrys_{X_{\d}/S_n}$-modules.
\end{itemize}

\noindent The group structure is defined by the tensor product as
$\cO^\logcrys_{X_{\d}\times_{S_n}S_n[\varepsilon]/S_n[\varepsilon]}$-
modules.
The identity element is equal
to~$\cO^\logcrys_{X_{\d}\times_{S_n}S_n[\varepsilon]/S_n[\varepsilon]}$.
Multiplication by elements of\/~$1+\varepsilon \WW_n(k)$ induces
the structure of~$\WW_n(k)$-module. We might describe an element
of~$\InfDef\Bigl(\cO^\logcrys_{X_{\d}/S_n}\Bigr)$ as the giving of
a pair $(L, \alpha )$
\begin{itemize}
\item[\rm i)] $L$ yields an element
of~$\InfDef\Bigl(\cO^\logcrys_{X_0/S_n}\Bigr)$;
\item[\rm ii)] and $\alpha\colon
(d_0^1)^*(L)\isomarrow (d_1^1)^*(L)$ as elements
of~$\InfDef\Bigl(\cO^\logcrys_{X_1/S_n}\Bigr)$ (\ie as sheaves of
$\cO^\logcrys_{X_1/S_n[\varepsilon]}$-modules such
that~$i_{X_1,\logcrys}^*(\alpha)={\rm Id}$) satisfying the cocycle
condition $$\bigl((d_1^2)^*(\alpha)\bigr)^{-1}\circ
\bigl((d_2^2)^*(\alpha)\bigr)\circ\bigl((d_0^2)^*(\alpha)
\bigr)={\rm Id}.$$ \end{itemize}

\noindent Note that the functor $\InfDef$ is functorial
in~$X_{\d}$ and, in particular, the Frobenius $X_{\d}\to X_{\d}^{(p)}$
defines a $\sigma$-linear homomorphism
on~$\InfDef\bigl(\cO^\logcrys_{X_{\d}/S_n}\bigr)$.\\

\noindent We link $\InfDef$ to $\bLogCrysno$ as follows.
For every $i$ let $f_i\colon
\bigl(X_i[\varepsilon]/S_n[\varepsilon]\bigr)_\logcrys \to
\bigl(X_i/S_n\bigr)_\logcrys$ be the standard morphism of topoi:
for $(U \subset T,\delta,L)$ an object of the crystalline site of
$X_i$ relative to $S_n$ we let $f_i^{-1}(U \subset T,\delta,L)$ be
$U[\varepsilon]\subset T[\varepsilon]$. As in \cite[\S
II.1.5]{MazurMessing} we get an exact sequence $0\rightarrow
\cO_{X_i/S_n}^\logcrys \rightarrow
(\cO_{X_i[\varepsilon]/S_n[\varepsilon]}^\logcrys)^* \rightarrow
(\cO_{X_i/S_n}^\logcrys)^* \rightarrow 0$ compatibly with the simplicial
structure.
Let~$\Lie\left(\bLogCrysno_{(X_{\d},Y_{\d})/S_n}\right)$ be the
group of $S_n[\varepsilon]$--valued points
of\/~$\bLogCrysno_{(X_{\d},Y_{\d})/S_n}$ reducing to the identity
modulo~$\varepsilon$. We conclude that
$$\Lie\left(\bLogCrysno_{(X_{\d},Y_{\d})/S_n}\right)\isomarrow
\InfDef\Bigl(\cO^\logcrys_{X_{\d}/S_n}\Bigr).$$Analogously, let
$\Lie\left(\bLogCrys_{(X_{\d},Y_{\d})/S_n}\right)$ be the group of
$S_n[\varepsilon]$--valued points
of\/~$\bLogCrys_{(X_{\d},Y_{\d})/S_n}$ reducing to the identity
modulo~$\varepsilon$.  Then, we have a canonical, Frobenius
equivariant homomorphism
$$\Lie\biggl(\bLogCrys_{(X_{\d},Y_{\d})/S_n}\biggr)
\llongrightarrow
\Lie\left(\bLogCrysno_{(X_{\d},Y_{\d})/S_n}\right)\isomarrow
\InfDef\Bigl(\cO^\logcrys_{X_{\d}/S_n}\Bigr).$$

\begin{lemma}\label{cohomologygroups} We have a canonical
isomorphism
$$\InfDef\Bigl(\cO^\logcrys_{X_{\d}/S_n}\Bigr)\cong
\HH^1_{\logcrys}\Bigl(\bigl(X_{\d},Y_{\d}\bigr)/S_n\Bigr)$$
compatibly with Frobenius and the exact sequences:
\begin{multline*} 0\llongrightarrow
\frac{\Ker\Bigl(\Aut\bigl(\cO^\logcrys_{X_{1}/S_n[\varepsilon]}
\bigr) \rightarrow \Aut\bigl(\cO^\logcrys_{X_{2}/S_n[\varepsilon]}
\bigr)\Bigr)}{\Im\Bigl(\Aut\bigl(\cO^\logcrys_{X_{0}/S_n[\varepsilon]}
\bigr) \rightarrow \Aut\bigl(\cO^\logcrys_{X_{1}/S_n[\varepsilon]}
\bigr)\Bigr)  } \llongrightarrow \InfDef\Bigl(
\cO^\logcrys_{X_{\d}/S_n}\Bigr) \\
\llongrightarrow
\Ker\Bigl(\InfDef\bigl(\cO^\logcrys_{X_0/S_n}\bigr) \rightarrow
\InfDef\bigl(\cO^\logcrys_{X_1/S_n}\bigr)\Bigr)\end{multline*}and

\begin{multline*} 0\llongrightarrow
\frac{\Ker\Bigl(\H^0_{\logcrys}\bigl(X_1/S_n\bigr) \rightarrow
\H^0_{\logcrys}\bigl(X_2/S_n\bigr)\Bigr)}{\Im\Bigl(\H^0_{\logcrys}
\bigl(X_0/S_n\bigr) \rightarrow
\H^0_{\logcrys}\bigl(X_1/S_n\bigr)\Bigr)} \llongrightarrow
\HH^1_{\logcrys}\bigl((X_{\d},Y_{\d})/S_n\bigr) \\
\llongrightarrow
\Ker\Bigl(\H^1_{\logcrys}\bigl((X_0,Y_0)/S_n\bigr) \rightarrow
\H^1_{\logcrys}\bigl((X_1,Y_1)/S_n\bigr)\Bigr)\end{multline*}
where~$\Aut\bigl( \cO^\logcrys_{X_i/S_n[\varepsilon]} \bigr)$
consists of the automorphisms reducing to the identity
via~$i_{X_1,\logcrys}^*$.
\end{lemma}
\begin{proof} Using the identification $\Lie\left(\bLogCrysno
\bigl((X_{\d},Y_{\d})/T\bigr)\right)\isomarrow
\InfDef\Bigl(\cO^\logcrys_{X_{\d}/S_n}\Bigr)$, the argument is a
simplicial variant of  \cite[\S II.1.5]{MazurMessing}. Note that
Frobenius on
$(\cO_{X_i[\varepsilon]/S_n[\varepsilon]}^\logcrys)^*$ induces the
Frobenius on $\cO_{X_i/S_n}^\logcrys$. \end{proof}

\begin{lemma}\label{LocGrysofproduct} Let~$X$ and~$Y$ be projective
and smooth schemes over~$k$. Then,
$$\LogCrys_{X\times_kY/S_n}=\LogCrys_{X/S_n}\times
\LogCrys_{Y/S_n}.$$
\end{lemma}
\begin{proof} We proceed by induction on~$n$. If~$n=0$,
then~$S_0=\Spec(k)$ and for every logarithmic scheme~$Z$ smooth
and log-smooth over~$k$ the group of $\overline{k}$-points of the
functor~$\LogCrys_{Z/k}$ is the group of isomorphism classes of
invertible sheaves on~$Z\times_k \overline{k}$ algebraically
equivalent to~$0$ and endowed with a logarithmic integrable
connection. Hence, denoting by~$\omega^1_{Z/k}$ the sheaf of
logarithmic differentials, we have an exact sequence
$$0\llongrightarrow \H^0(Z,\omega_{Z/k}^1)_{d=0}\llongrightarrow
\LogCrys_{Z/k} \llongrightarrow \Pic^0_{Z/k}.$$Note that
$\Pic^{0,\red}_{X\times_k Y/k}\cong\Pic^{0,\red}_{X/k}\times_k\Pic^{0,\red}_{Y/k}$
(see \cite{BL}).
Let~$\pi_1$ and $\pi_2$ be the two projections from~$X\times_kY$ to $X$ and $Y$ respectively. Then,
$\omega^1_{X\times_k Y/k}=\pi_1^*(\omega^1_{X/k})\oplus
\pi_2^*(\omega^1_{Y/k})$. Since~$X$ and~$Y$ are projective we
have,
$\pi_{1,*}\bigl(\pi_1^*(\omega^1_{X/k})\bigr)=\omega^1_{X/k}$ and
$\pi_{2,*}\bigl(\pi_2^*(\omega^1_{Y/k})\bigr)=\omega^1_{Y/k}$.
Hence
$$\H^0(X\times_kY,\omega_{X\times_k Y/k}^1)_{d=0}\cong
\H^0(X,\omega_{X/k}^1)_{d=0}\oplus\H^0(Y,\omega_{Y/k}^1)_{d=0}.$$ We
conclude
that the lemma holds for~$n=0$. Suppose that the lemma is proven
for~$n$. Pulling back via the projections $X\times_k Y\rightarrow
X$ and $X\times_k Y\rightarrow Y$ we obtain the following
commutative diagram
\begin{displaymath}
\begin{CD} 0  & & 0\\
@VVV  @VVV\\
\H^1_\logcrys\bigl(X/k\bigr)\oplus
\H^1_\logcrys\bigl(Y/k\bigr)@>{\sim}>>\H^1_\logcrys\bigl(X\times_k
Y/k\bigr)
\\
@V{\cdot p^n}VV @V{\cdot p^n}VV\\
\LogCrys_{X/S_{n+1}}\times \LogCrys_{Y/S_{n+1}} @>>>
\LogCrys_{X\times_k
Y/S_{n+1}}\\
@VVV  @VVV\\
\LogCrys_{X/S_n}\times \LogCrys_{Y/S_n}@>{\sim}>>
\LogCrys_{X\times
Y/S_n}.\\
\end{CD}
\end{displaymath}
\noindent Indeed, the kernel of the bottom vertical right hand
side map is isomorphic to the crystals of invertible modules on
the crystalline site of~$X\times_k Y$ relative to~$S_{n+1}$
reducing to the trivial crystal~$\cO^\logcrys_{X\times_kY/S_n}$.
Using `Cech cohomology' by \ref{Chechnotation}, it can be
identified with the cohomology group~$\H^1_\logcrys\bigl(X\times_k
Y/k\bigr)\otimes_{k} (p^n\WW_{n+1}(k))$. One gets an analogous
description of the kernel  of the bottom vertical left hand side
map. The bottom horizontal isomorphism exists by inductive
hypothesis, while the top horizontal isomorphism follows from the
equality~$\H^1_\logcrys\bigl(X\times_k Y/k\bigr)=\Lie
\bigl(\LogCrys_{X\times_k Y/k}\bigr)$ and the $n=0$ case.
\end{proof}

\subsection{The case of abelian
varieties}\label{thecaseofabelianvarieties} This is Grothendieck's
remark \cite{MG}, \cf \cite{MazurMessing}. Let~$B_0$ be an abelian
variety over~$S_0$. Since the deformation functor of abelian
varieties is unobstructed, there exists an abelian scheme~$B_n$
over~$S_n$ lifting~$B_0$. Then, $B_n[\varepsilon]$ is a lifting
of~$B_0$, smooth over~$S_n[\varepsilon]$. Hence, the category of
crystals of invertible~$\cO^\crys_{B_0\times_{S_n}
S_n[\varepsilon]/S_n[\varepsilon]}$-modules over the nilpotent
crystalline site of $B_0\times_{S_0} S_0[\varepsilon]$ relative to
$S_n[\varepsilon]$ is equivalent to the category of line bundles
over~$B_n[\varepsilon]$ with integrable connection.
Let~$\E\bigl(B_n^\vee\bigr)$ be the universal extension of the
abelian scheme dual to~$B_n$. It actually classifies isomorphism
classes of line bundles over~$B_n$ algebraically equivalent to~$0$
and endowed with integrable connection. Hence, we have an
equivalence of functors over the fppf site of $S_n$, compatible
with group structures,
$$\E\bigl(B_n^\vee \bigr)\cong \Pic^{\natural, 0}_{B_n/S_n} \cong
\LogCrys_{B_0/S_n}.$$Taking $\Lie$ we get a natural isomorphism
of~$\cO_{S_n}$-modules
$$\TT^\crys(B_0)\tensor
\cO_{S_n}\cong \Lie\Bigl(\E\bigl(B_n^\vee \bigr) \Bigr) \cong
\Lie\left(\LogCrys_{B_0/S_n}\right) \cong
\InfDef\Bigl(\cO^\logcrys_{B_0/S_n}\Bigr).$$

\subsection{The compact case}
Here we assume that~$Y_{\d}=\emptyset$, thus all logarithmic
structures are trivial. We then write $\Crys$ for the functor
$\LogCrys$ defined in the previous sections. We also assume that
each irreducible component of~$X_0$, $X_1$ and~$X_2$ is
geometrically irreducible. For the general case, since $k$ is
assumed to be perfect we pass to a separable closure $ k^{\rm
sep}$ of $k$. Then, one easily verifies that the
morphisms~$\phi_i$ in~\ref{albPic}, and
consequently~$\phi_{\crys}^*$ in~\ref{compactinjective}, are
$\Gal\bigl(k^{\rm sep}/k\bigr)$-equivariant.

\subsubsection{Comparison via the Albanese}\label{albPic}
Let $i=0$,~$1$ or~$2$. Let~$X_i=:\amalg X_i^j$ be the irreducible
components of~$X_i$. Define
$$X_i^{[2]}:=\amalg_j X_i^j\times_k X_i^j,\qquad\hbox{{\rm and}}\qquad
A(X_i):=\amalg_j \Alb(X_i^j/k).$$For every~$j$ there exists a
unique morphism $\phi_i^j\colon X_i^j\times_k X_i^j\rightarrow
\Alb(X_i^j/k)$ such that for every~$x_i^j\in X_i^j(\bar{k})$ the
map~$\phi_i^j(x_i^j,\_)\colon X_i^j\rightarrow \Alb(X_i^j/k)$ is
the Albanese map sending~$x_i^j\mapsto 0$. Let $$\amalg_j
\phi_i^j=:\phi_i\colon X_i^{[2]}\llongrightarrow A(X_i)$$be the
induced map. By the functoriality of the Albanese construction it
is easily checked that the maps $d_0^1,d_1^1\colon X_1\rightarrow
X_0$ and $d_0^2, d_1^2, d_2^2\colon X_2\rightarrow X_1$ induce
maps $d_0^1,d_1^1\colon X_1^{[2]}\rightarrow X_0^{[2]}$ and
$d_0^2, d_1^2, d_2^2\colon X_2^{[2]}\rightarrow X_1^{[2]}$ and
maps $d_0^1,d_1^1\colon A(X_1)\rightarrow A(X_0)$ and $d_0^2,
d_1^2, d_2^2\colon A(X_2)\rightarrow A(X_1)$. Then,
$\bigl(X_i^{[2]},d_i^j\bigr)$ and~$\bigl(A(X_i),d_i^j\bigr)$ are
simplicial schemes over~$k$ and the maps $\left\{\phi_i\colon
X_i^{[2]}\rightarrow A(X_i)\right\}_{i=0,1,2}$ define a map of
$2$-truncated simplicial schemes. By functoriality of~$\Crys$ we
have the map
$$\phi^*\colon\Ker\Bigl(\Crys_{A(X_0)/S_n}\rightarrow
\Crys_{A(X_1)/S_n}\Bigr)\llongrightarrow
\Ker\Bigl(\Crys_{X_0^{[2]}/S_n}\rightarrow
\Crys_{X_1^{[2]}/S_n}\Bigr).
$$Let~$\A_0$ be the reduced kernel of $(d_0^1)^*-(d_1^1)^* \colon
\Pic^{0, \red}_{X_0/k}\rightarrow \Pic^{0,\red}_{X_1/k}$.
Let~$B_{0,n}$ and~$B_{1,n}$ be abelian schemes over~$S_n$
lifting~$\bigl(\Pic^{0,\red}_{X_0/S_0}\bigr)^\vee$
and~$\bigl(\Pic^{0,\red}_{X_1/S_0}\bigr)^\vee$. Let~$\A_n$ be an
abelian scheme over~$S_n$ lifting~$\A_0$. The  sequence
$$\bigl(\Pic^{0,\red}_{X_1/S_0}\bigr)^\vee
\longrightarrow \bigl(\Pic^{0,\red}_{X_0/S_0}\bigr)^\vee
\llongrightarrow \A_0^\vee\llongrightarrow 0$$is exact. Hence, we
get from~\ref{thecaseofabelianvarieties} a commutative diagram of
sheaves over the fppf site of $S_n$:
\begin{equation}
\begin{CD}
\Crys_{\A_0^\vee/S_n} @>>>
\Crys_{\bigl(\Pic^{0,\red}_{X_0/S_0}\bigr)^\vee/S_n} @>>>
\Crys_{\bigl(\Pic^{0,\red}_{X_1/S_0}\bigr)^\vee/S_n}\\
@V{\wr}VV @V{\wr}VV @V{\wr}VV\\ \E\bigl(A_n^\vee\bigr)@>{a_n}>>
\E\bigl(B_{0,n}\bigr)@>{b_n}>> \E\bigl(B_{1,n}\bigr).\\
\end{CD}
\end{equation}
\noindent Since~$\A_0 \subset \Pic^{0,\red}_{X_0/S_0}$ is a closed
immersion,~$a_0$ is a closed immersion. Since~$a_n\times_{S_n}
S_0$ is~$a_0$, then also~$a_n$ is a closed immersion. In
particular, $a_n$ is injective. Using Lemma~\ref{LocGrysofproduct}
to get the second equality we have
$$\Crys_{A(X_i)/S_n}=\prod_j
\Crys_{\Alb(X_i^j/k)/S_n}=\Crys_{\prod_j\Alb(X_i^j/k)/S_n}.$$ Note
that $\Crys_{\bigl(\Pic^0_{X_i/S_0}\bigr)^\vee/S_n}
=\Crys_{\prod_j\Alb(X_i^j/k)/S_n}$. In conclusion, we obtain
$$\Crys_{\A_0^\vee/S_n} \hookrightarrow
\Ker\Bigl(\Crys_{A(X_0)/S_n}\rightarrow
\Crys_{A(X_1)/S_n}\Bigr).$$For every~$i$ and~$j$ the composite of
the diagonal embedding~$X_i^j\rightarrow X_i^j\times_k X_i^j$
and~$\phi_i^j\colon X_i^j\times_k X_i^j \rightarrow \Alb(X_i/k)$
factors via the identity of~$\Alb(X_i/k)$. Hence, the homomorphism
$$\Crys_{\Alb(X_i^j/k)/S_n}\llongrightarrow \Crys_{X_i^j\times_k
X_i^j/S_n}\llongrightarrow \Crys_{X_i^j/S_n}$$is the zero map. By
Lemma~\ref{LocGrysofproduct} we have~$\Crys_{X_i^j\times_k
X_i^j/S_n}=\Crys_{X_i^j/S_n}\times \Crys_{X_i^j/S_n}$ and the map
to~$\Crys_{X_i^j/S_n}$ is the sum. Hence, the map of fppf sheaves
on $S_n$
$$(\phi_i^j)^*\colon
\Crys_{\Alb(X_i^j/k)/S_n}\llongrightarrow \Crys_{X_i^j\times_k
X_i^j/S_n}$$factors via~$({\rm Id},-{\rm
Id})\colon\Crys_{X_i^j/S_n}\rightarrow \Crys_{X_i^j/S_n}\times
\Crys_{X_i^j/S_n}$. Hence, the map $\phi^*$ induces a homomorphism
$$\phi^*\colon \Crys_{\A_0^\vee/S_n}\llongrightarrow
\Ker\Bigl(\Crys_{X_0/S_n}\rightarrow
\Crys_{X_1/S_n}\Bigr).$$Taking~$\Lie$ of these functors we get a
map $$\phi_\crys^*\colon
\H^1_\crys\bigl(\A_0^\vee/S_n\bigr)\rightarrow
\Ker\Bigl(\H^1_\crys\bigl(X_0/S_n\bigr)\rightarrow
\H^1_\crys\bigl(X_1/S_n\bigr)\Bigr)$$

\begin{proposition}\label{isomorphicrys} The induced map
$$\TT_\crys\bigl(\A_0\bigr)
\longby{\simeq} \Ker\Bigl(\H^1_\crys\bigl(X_0\bigr)\rightarrow
\H^1_\crys\bigl(X_1\bigr)\Bigr)
$$is an isomorphism.
\end{proposition}
\begin{proof}
Using Illusie's theory \cite[\S II.1]{Illusie} we can compute
crystalline cohomology of the smooth proper schemes $X_i$ as the
hypercohomology of De Rham Witt complexes. In particular, by
\cite[II.3.11.2]{Illusie} the groups $\H^1_\crys(X_i)$ and
$\H^1_\crys\bigl(\Alb(X_i/k)\bigr)$ coincide. Since the sequence
$\bigl(\Pic^{0,\red}_{X_1/S_0}\bigr)^\vee \rightarrow
\bigl(\Pic^{0,\red}_{X_0/S_0}\bigr)^\vee \rightarrow
\A_0^\vee\rightarrow 0$ is exact and the contravariant Dieudonn\'e
functor defines a fully faithful functor from the category of
abelian varieties to the category of  $F$-crystals, the associated
sequence $$0 \llongrightarrow\H^1_\crys\bigl(\A_0^\vee\bigr)
\llongrightarrow
\H^1_\crys\bigl(\bigl(\Pic^{0,\red}_{X_0/S_0}\bigr)^\vee\bigr)
\llongrightarrow
\H^1_\crys\bigl(\bigl(\Pic^{0,\red}_{X_1/S_0}\bigr)^\vee\bigr)$$
is exact. Since
$\H^1_\crys\bigl(\A_0^\vee\bigr)=\TT_\crys\bigl(\A_0\bigr)$, the
conclusion follows.
\end{proof}

\subsubsection{Souping up}\label{symplicialcomparison}
Let~$\G_0$ be the semiabelian scheme~$\bPic^{0,\red}_{X_{\d}/k}$.
Let\/~$\A_0$ be the abelian part of~$\G_0$. Fix~$n\in\NN$.
Let~$\G_n\rightarrow S_n$ be a lift of~$\G_0\rightarrow S_0$ as a
semiabelian scheme. Let~$\A_n$ be the abelian part of~$\G_n$. By
the discussion above, see \ref{albPic}, and the
identification~$\E(\A_n)\cong \Crys_{\A_0^\vee/S_n}$, see
\ref{thecaseofabelianvarieties}, we have a homomorphism of sheaves
of abelian groups on the fppf site of~$S_n$
$$\phi^*\colon \E(\A_n)\llongrightarrow
\Ker\Bigl(\Crys_{X_0/S_n}\rightarrow
\Crys_{X_1/S_n}\Bigr).$$Define $\E_n'$ as the fibred product
\begin{displaymath}\begin{CD} \E_n' @>>>
\bCrys_{X_{\d}/S_n}\\ @VVV @VVV\\
\E\bigl(\A_n\bigr) @>{\phi^*}>>
\Ker\Bigl(\Crys_{X_0/S_n}\rightarrow \Crys_{X_1/S_n}\Bigr)\\
\end{CD}\end{displaymath} \noindent By~\ref{categorygroups} the
kernel of the right vertical map is represented by a torus~$\T_n$
over~$S_n$.

\begin{lemma}\label{extensionSArEequal} We have
$\E_n'\simeq \E\bigl(\G_n\bigr)$ as extensions of\/ $\E(\A_n)$ by
$\T_n$.
\end{lemma}
\begin{proof}
By the discussion above we have a sequence $$0\llongrightarrow
\T_n\llongrightarrow \E_n'\llongrightarrow
\E\bigl(\A_n\bigr)\llongrightarrow 0.$$We claim that it is exact
on the right and, hence, exact. We have a homomorphism of abelian
sheaves
$$\Ker\Bigl(\Crys_{X_0/S_n}\rightarrow
\Crys_{X_1/S_n}\Bigr)\llongrightarrow \Bigl(\pi_{3,*}^\crys
\bigl(\GG_{m,X_3}\bigr)/\Im\bigl(\pi_{2,*}^\crys
\bigl(\GG_{m,X_2}\bigr) \bigr)\Bigr)$$associating to any crystal
$L$ of invertible modules on $(X_0/S_n)$, such that there exists
an isomorphism $\alpha\colon(d_0^1)^*(L)\cong (d_1^1)^*(L)$, the
class of $\bigl((d_1^2)^*(\alpha)\bigr)^{-1}\circ
\bigl((d_2^2)^*(\alpha)\bigr)\circ\bigl((d_0^2)^*(\alpha) \bigr)$.
Note that such a class is trivial if and only if there exists
$\alpha$ satisfying the cocycle condition, \ie if and only if
$\alpha$ lifts to an element of $\bCrys_{X_{\d}/S_n}$. Thus, to
prove the exactness it suffices to prove that the composite
$\Upsilon_n\colon \E(\A_n)\to \T_n':=\pi_{3,*}^\crys
\bigl(\GG_{m,X_3}\bigr)/\Im\bigl(\pi_{2,*}^\crys
\bigl(\GG_{m,X_2}\bigr) \bigr)$ is trivial. We proceed by
induction on $n$. For $n=0$ the triviality follows from $\Hom
(\GG_{a,k},\GG_ {m,k})=0$ and $\Hom (\A_0,\GG_{m,k})=0$. If
$\Upsilon_n$ is trivial, then $\Upsilon_{n+1}$ factors as
$\E(\A_0)\to \Ker\bigl(\T_{n+1}'\to \T_n'\bigr)$. The group on the
right hand side is a sum of $\GG_{a,k}$'s. Since $\E(\A_0)$ is the
universal $\GG_{a,k}$-extension of $\A_0$, we have
$\Hom(\E(\A_0),\GG_{a,k})=\{0\}$. This concludes the inductive
step and proves the claimed exactness.\\

Interpreting~$\Crys_{\  /S_0}$ as isomorphism classes of invertible
sheaves endowed with integrable connection, $\E_0'$ is the fibred
product of~$\bPic^{\natural,0}_{X_{\d}/k}$
and~$\Pic^{\natural,0}_{\A_0/k}$ over the kernel
of~$\Pic^{\natural,0}_{X_0/k}\rightarrow
\Pic^{\natural,0}_{X_1/k}$. It follows from the definition
of~$\G_0$ that~$\E_0'\cong \E\bigl(\G_0\bigr)$.

Since~$\E_n'$ is an extension of representable sheaves over the
fppf site of~$S_n$, \cite[Prop. 17.4]{OG} guarantees that~$\E_n'$
is itself representable. Moreover, the base change of~$\E_n'$
to~$S_0$ is isomorphic to~$\E_0'$. It follows from the crystalline
nature of~$\E\bigl(\G_n\bigr)$ (see \ref{gory}) that there exists
a unique isomorphism of extensions~$\E\bigl(\G_n\bigr)\cong \E'_n$
lifting the isomorphism~$\E_n\times_{S_n} S_0\cong \E_0'\cong
\E\bigl(\G_0\bigr)$ over~$S_0$.
\end{proof}

\subsubsection{The comparison}\label{cOMP}
In conclusion, given a lift~$\G_n\rightarrow S_n$
of~$\G_0\rightarrow S_0$, we get a canonical homomorphism
$$\E\Bigl(\G_n\Bigr)
\llongrightarrow \bCrys_{X_{\d}/S_n}.$$Taking $\Lie$ and
using~\ref{ccoommppaarree} and~\ref{cohomologygroups}, we get
functorial homomorphisms of~$\cO_{S_n}$-modules
\begin{equation}
\begin{CD}
\Lie\Bigl(\E\bigl(\G_n\bigr)\Bigr)
@>>> \InfDef\Bigl( \cO^\crys_{X_{\d}/S_n}\Bigr) \\
@V{\wr}VV @V{\wr}VV\\ \Tcrys(\G)\tensor_{\WW(k)} \WW_n(k)& &
\HH^1_{\crys}\bigl(X_{\d}
/S_n\bigr).\\
\end{CD}
\end{equation}

\begin{proposition}\label{compactinjective}
For $p\geq 3$ the inverse limit of the above maps defines an
isomorphism
$$\phi_{\crys}^*\colon \Tcrys(\G_0)
\longby{\simeq} \HH^1_{\crys}\bigl(X_{\d}\bigr)(1)$$respecting the weight
filtrations and the actions of Frobenius in~\ref{propTcrys} and
in~\ref{verygory}.
\end{proposition}
\begin{proof} By~\ref{propTcrys} and~\ref{verygory} the two $\WW(k)$-modules admit a weight filtration, preserved by $\phi_\crys^*$. Hence, it suffices to
prove that it is an isomorphism on the various graded pieces. For
$\Gr_1$ use~\ref{isomorphicrys} and for $\Gr_0$ use \ref{fil1} (\cf
\ref{withouttorsion}).
\end{proof}

\subsection{On the case $p=2$}\label{casep=2} For $p=2$ there is a
problem in defining the functor~$\LogCrys$ in~\ref{functorLogGrys}
since the divided power structure we consider on $p\WW_{n+1}(k)$
is {\it not} nilpotent. Following \cite{MazurMessing} we may 
replace $\LogCrys$ with the sheaf
$$\barbLogCrys_{(X_{\d},Y_{\d})/S_n}$$associated to the
following presheaf over the fppf site of~$S_n$. Let~$T$ be a
scheme over~$S_n$ flat and of finite presentation. Consider the
logarithmic and divided power structure on~$T$ induced from those
on~$S_n$. Consider the group of isomorphism classes of simplicial
crystals of invertible $\cO^\logcrys_{X_{\d}\times_{S_n}
T/T}$-modules for the {\it full} crystalline site
of~$X_{\d}\times_{S_n} T$ relative to~$ T$ endowed with a
trivialization over $T\times_S S_0$; compare with \cite[\S
11.1]{MazurMessing}.\\

By \cite[Ex. 4.14 \& Cor. 6.8 ]{BerthelotOgus} the category of
crystals of invertible $\cO^\logcrys_{X_{\d}\times_{S_0}
T/T}$-modules for the full crystalline site of $X_{\d}\times_{S_0}
T$ relative to~$T$ is equivalent to the category of invertible
simplicial sheaves on~$X_{\d}$ with an integrable, {\it
quasi-nilpotent} connection. In particular, if $\A_0$ is an
abelian variety, $\barLogCrys_{\A_0/S_n}$ is the formal group
$\overline{\E}(\A_n^\vee)$ of the universal extension of the
dual of a lifting $\A_n$ of $\A_0$ to $S_n$. Note
that~$\barLogCrys_{\A_0/S_n}$ and~$\E(\A_n^\vee)$ have the
same Lie algebras.\\

For $p\geq 3$, if we consider the {\it nilpotent} crystalline
site,  the category of crystals of invertible
$\cO^\logcrys_{X_{\d}\times_{S_n} T/T}$-modules is equivalent to
the category of  invertible simplicial sheaves on~$X_{\d}$ with an
integrable connection (no nilpotence assumption). Thus, we get a
natural inclusion $\barLogCrys \subset \LogCrys$ inducing an
isomorphism at the level of tangent spaces $\Lie$. By \cite[Ex.
4.14]{BerthelotOgus}  $\barbLogCrys_{X_{\d}/S_n}$ is the fiber
product of $\bLogCrys_{X_{\d}/S_n}$ and
$\barbLogCrys_{X_{\d}/S_0}$
over $\bLogCrys_{X_{\d}/S_0}$. \\

Unfortunately, this is not enough to get Theorem B${}^{\prime}$ for $p=2$
since~\ref{comparisonnoncompact} is lacking. 

\subsection{The general case}
See \ref{simpic} for the general definition of the Picard
$1$-motive $\M_0 = \Pic^+ (V)$ associated to the algebraic
$k$-variety $V$. Let\/ $\bDiv^0_{\bar Y_{\d}} \bigl(\bar
X_{\d}\bigr)$ be the group defined in (\ref{Weil}) over $\bar k$
along with a Galois action; that is an \'etale sheaf on~$S_0=\Spec
(k)$ and identified with the discrete part of $\Pic^+ (V)$.
\smallskip
\begin{definition}
\noindent Let~$n\in\NN$.  Define $$\v_n\colon \X \llongrightarrow
\bLogCrys\bigl((X_{\d},Y_{\d})/S_n\bigr)$$to be the homomorphism
associating to a local section~$D$ of\/~$\X$ the
crystal~$\cO^\logcrys_{X_{\d}/S_n}\bigl[D^{-1}\bigr]$ of
invertible $\cO^\logcrys_{X_{\d}/S_n}$-modules via the procedure
in~\ref{logcrysWeil}.
\end{definition}

\begin{proposition}\label{comparisonnoncompact}
For $n\in\NN$ there is a unique homomorphism
$$\psi_n\colon\E\bigl(\M_n\bigr)_{\G_n}   \llongrightarrow
\bLogCrys_{(X_{\d},Y_{\d})/S_n}$$ such that
\begin{itemize} \item[\rm a)] for varying~$n$ the $\psi_n$ are
compatible, and
\item[\rm b)] $\psi_n$ makes the following diagram commute
\begin{equation}
\begin{CD}
\E\bigl(\G_n\bigr) @>>> \bCrys_{X_{\d}/S_n}\\
@VVV @VVV \\ \E\bigl(\M_n\bigr)_{\G_n}   @>\psi_n>>
\bLogCrys_{(X_{\d},Y_{\d})/S_n}\\@A{\u_n}AA @A{\v_n}AA\\ \X
@>=>> \X,\\
\end{CD}
\end{equation}where the morphism of top line is defined
in~\ref{cOMP}.
\end{itemize}
\end{proposition}
\begin{proof} The $\psi_n$ are clearly obtained after an \'etale
covering of~$S_0$. We may
then assume that~$\X$ is a constant group scheme. Note
that~$\M_n=\bigl[ \X\rightarrow \G_n\bigr]$. Proceed by induction
on~$n$. If~$n=0$, the claim is granted by \ref{cOMP} and the following
commutative diagram
\begin{equation}
\begin{CD} \E\bigl(\M_0\bigr)_{\G_0}   @>\psi_0>>
\bLogCrys_{(X_{\d},Y_{\d})/S_0}\\@V{p_0}VV @V{q_0}VV\\ \G_0
@>=>> \bPic^{0, \red}_{X_{\d}/S_0}. \\
\end{CD}
\end{equation}
\noindent Note that~$\G_0$ is the semiabelian part of the
universal extension of~$\M_0$, the kernel of the map $p_0$ is
$\Ext (\M_0, \GG_a)^{\vee}$, and the diagram is provided by
universality as soon as $q_0$ is a vector group extension of
$\M_0$. To see this, note that $$\bPic^{\natural,
0}_{X_{\d}/S_0}\into \bPic^{\natural -\log,
0}_{(X_{\d},Y_{\d})/S_0}= \bLogCrys_{(X_{\d},Y_{\d})/S_0}$$ and the
kernel of the surjection $q_0$ is the group scheme representing
connections on the structure sheaf of the simplicial
scheme~$X_{\d}$ with logarithmic poles along the simplicial
divisor~$Y_{\d}$. Hence, it is isomorphic to the vector group
scheme
$$\HH^0\Bigl(X_{\d},\Omega^1_{X_{\d}/S_0}\bigl( \log (Y_{\d}\bigr)
\Bigr)_{d=0}$$
by the arguments in the proof of Lemma~\ref{repnat}.
The universal property of~$\E\bigl(\M_0\bigr)$ yields the claimed
$\psi_0$ uniquely.

Suppose that~$\psi_n$ has been defined for~$n<N$. Since we have a
non-canonical splitting
$$\E\bigl(\M_N\bigr)_{\G_N}\isomarrow\E\bigl(\G_N\bigr)
\times_{S_N} \E\bigl[\X\rightarrow 0\bigr]$$we certainly can
lift~$\psi_{N-1}$ to a homomorphism $$\psi_N'\colon
\E\bigl(\M_N\bigr)_{\G_N} \llongrightarrow
\bLogCrys_{(X_{\d},Y_{\d})/S_N}$$such that the upper square in~(b)
commutes. The possible liftings~$\psi_N'$ are a principal
homogeneous space under $$\Hom\left(\X,
\Lie\Bigl(\bLogCrys_{(X_{\d},Y_{\d})/S_N}\Bigr)\right).$$Hence,
there is a unique~$\psi_N'$ such that the lower square in~(b)
commutes.
\end{proof}

\subsubsection{Concluding}\label{endproof}
By taking~$\Lie$ of the map~$\psi_n$ in
Proposition~\ref{comparisonnoncompact}, we get a $\WW_n(k)$-linear
homomorphism
$$\Tcrys \bigl(\M_0\bigr)\otimes_{\WW(k)}\WW_n(k) \llongrightarrow
\HH^1_\logcrys\Bigl(\bigl(X_{\d}, Y_{\d}\bigr)/S_n\Bigr)(1)$$as
claimed in Theorem~B${}^{\prime}$, compatible for varying~$n$,
preserving the weight filtrations~$W$ and the action of Frobenius.
Due to Proposition~\ref{compactinjective} in order to conclude the
proof of Theorem~B${}^{\prime}$, it suffices to check that the
above map is an isomorphism after taking the inverse limit over
$n\in\NN$ on the $\Gr^W_0$ parts. Thus we are left to deal with the
case of $k =\bar k$ as all maps are clearly Galois equivariants.
By~\ref{verygory} we have $\iota\colon \Gr_2^W\left(\HH^1_\logcrys\bigl(
X_{\d}, Y_{\d}\bigr)\right)\subseteq \bDiv_{Y_{\d}}^0
(X_{\d})\otimes_\ZZ \WW(k)$. Furthermore $\bDiv_{Y_{\d}}^0
(X_{\d})\otimes_\ZZ \WW(k)=\Gr^W_0\left(\Tcrys
\bigl(\M_0\bigr) \right)$ by definition of $\M_0 = \Pic^+(V)$, see \ref{simpic}. By Proposition~\ref{comparisonnoncompact} we obtain a map
$$\eta\colon \Gr^W_0\left(\Tcrys \bigl(\M_0\bigr) \right) \to
\Gr^W_0\left(\HH^1_\logcrys\bigl(X_{\d}, Y_{\d}\bigr)(1)\right).$$We now show that $\iota\circ\eta$ is the identity modulo $p$. This suffices to conclude the proof of Theorem~B${}^{\prime}$.
Using Section \ref{fil2} for the case $S_n=S_0=\Spec(k)$ we get, by construction,  the following commutative diagram
$$
\begin{CD}   \Gr^W_0\left(\HH^1_\logcrys\bigl(X_{\d}, Y_{\d}\bigr)(1)\right)
@>\iota>> \bDiv_{Y_{\d}}^0 (X_{\d})\otimes_\ZZ \WW(k)\\@VVV @VVV\\ 
\HH^0\Bigl(X_{\d},\Omega^1_{X_{\d}/S_0}\bigl( \log (Y_{\d}\bigr)
\Bigr)_{d=0}/\HH^0\Bigl(X_{\d},\Omega^1_{X_{\d}/S_0}
\Bigr)_{d=0} @>{\rm Res}>>  \bDiv_{Y_{\d}} (X_{\d})\otimes_\ZZ k,\\
\end{CD}
$$where ${\rm Res}$ is the map defined by taking residues and $\iota$ is the inclusion above. By the proof of Proposition \ref{comparisonnoncompact} for every $D\in \bDiv_{Y_{\d}}^0 (X_{\d})=\X$ we have $\psi_0(u_0(D))=(\cO_{X_{\d}}(-D),d)$, where $d$ is the canonical connection on $\cO_{X_{\d}}(-D)$. Thus, the  map $\bar\psi_0\colon\E(\M_0)_{\G_0}/\E(\G_0)=\X \tensor \GG_a \rightarrow \bPic^{\natural -\log,
0}_{(X_{\d},Y_{\d})/S_0}/\bPic^{\natural,
0}_{X_{\d}/S_0}$ induced by $\psi_0$ is the unique $\A^1$-linear map sending $D\otimes 1$ to the class of $ (\cO_{X_{\d}}(-D),d)$ and the induced map on Lie algebras   
yields the lower horizontal arrow in the following commutative diagram
$$
\begin{CD} \Gr^W_0\left(\Tcrys \bigl(\M_0\bigr) \right)@>\eta>>  \Gr^W_0\left(\HH^1_\logcrys\bigl(X_{\d}, Y_{\d}\bigr)(1)\right) \\@VVV @VVV\\ 
\Gr^W_0\left(\Tcrys \bigl(\M_0\bigr) \right)\otimes \ZZ/p\ZZ @>\eta_k>> \HH^0\Bigl(X_{\d},\Omega^1_{X_{\d}/S_0}\bigl( \log (Y_{\d}\bigr)
\Bigr)_{d=0}/\HH^0\Bigl(X_{\d},\Omega^1_{X_{\d}/S_0}
\Bigr)_{d=0}.\\
\end{CD}
$$Furthermore, using the description of $\bar\psi_0$ given above, one verifies that ${\rm Res}\circ \eta_k$  is simply the inclusion $ \bDiv_{Y_{\d}}^0(X_{\d})\otimes_\ZZ k\into  \bDiv_{Y_{\d}} (X_{\d})\otimes_\ZZ k $. Patching the two diagrams we conclude that $\iota\circ \eta$ is the identity modulo $p$ as claimed.

\begin{corollary}\label{independenceH1crysV} Let $V$ be an algebraic variety over a perfect field $k$ of characteristic $p>2$.
The filtered $F$-$\WW (k)$-module $\H^1_\crys(V/\WW (k))\df\H^1_\crys(V_{\d}/\WW(k))$
is independent of the choice of the hypercovering $V_{\d}\to V$
(subject to the condition that $V_0\to V$ is generically \'etale).
\end{corollary}
\begin{proof} Let $V_{\d}\to V$ and $V_{\d}'\to V$ be two such hypercoverings.
We may assume that there is a map $\varphi_{\d}
\colon V_{\d}'\to V_{\d}$ of $V$-simplicial schemes compatibly
with the normal crossing boundaries, \ie which is the restriction
of a map $X_{\d}'\to X_{\d}$ (\cf Appendix~\ref{appendix}). By Theorem~B${}^{\prime}$, the Appendix~\ref{appendix} and functoriality  we get an induced isomorphism
$\TT_{\crys}(\Pic^+(V))\cong \H^1_\crys\bigl(V_{\d}/\WW(k)\bigr) \longby{\cong}
\H^1_\crys\bigl(V_{\d}'/\WW(k)\bigr)\cong \TT_{\crys}(\Pic^+(V))$  of filtered $F$-$\WW (k)$-modules.
\end{proof}

\appendix
\section{Appendix}\label{appendix}
We provide details showing that the cohomological Picard 1-motive
$\Pic^{+}(V)$ is independent of the choices made, \ie it is
independent of choices of hypercoverings and compactifications, over
perfect fields.
In characteristic zero this is provided by \cite[Prop. 2.5]{BM},
see also \cite[Remark 4.4.4]{BS}. However, the argument for positive
characteristics is slightly more involved.\\

Suppose we are given two such smooth proper hypercoverings
$f\colon V_{\d}\to V$ and $f'\colon V_{\d}'\to V$ which admit
smooth compactifications with normal crossing boundaries, $V_{\d}
= X_{\d} - Y_{\d}$ and $V_{\d}' =  X_{\d}' - Y_{\d}'$ and assume 
that $V_0 \rightarrow V$ and $V_0'\rightarrow V$  are generically
\'etale. We refer to De Jong's theory of alterations over perfect fields 
for their existence,  see \cite[Thm. 4.1]{DE}. Note that, as usual (\cf \cite[\S
6.2]{Deligne}), we can always find a  third one mapping to both
(\cf 5.1.7 and 5.2.4 in Expos\'e V bis of \cite{SGA4}). We then can
assume that there is a map $\varphi_{\d} \colon V_{\d}'\to V_{\d}$
of $V$-simplicial schemes compatibly with the normal crossing
boundaries, \ie it is the restriction of a proper map
$\varphi_{\d}\colon  X_{\d}'\to X_{\d}$. 

Then, by pulling back (see \cite[6.2]{BS}), we get a map of 1-motives
$\varphi^*\colon
\Pic^{+}(X_{\d}, Y_{\d})\to \Pic^{+}(X_{\d}', Y_{\d}')$ and  thus,
equivalently, we have a $\Gal (\bar{k}/k)$-equivariant map of
complexes $$\bar \varphi^*\colon [\bDiv_{\bar
Y_{\d}}^{0}(\bar X_{\d})\to \bPic^{0,\red}( \bar
X_{\d})]\to  [\bDiv_{\bar Y_{\d}'}^{0}(\bar X_{\d}')\to
\bPic^{0,\red}( \bar X_{\d}')].$$ We are then reduced to deal with
$k = \bar k$ since our constructions are natural enough to be
automatically compatible with  $\Gal (\bar{k}/k)$-actions. The
proof is divided in two steps. First of all, in
Section~\ref{independenceatl},
we show that, for $\ell\neq p={\rm char}(k)$, the map $\varphi^*$
induces an
isomorphism of the $\ell$-adic realization of the above
$1$-motives. This implies that $\varphi^*$ is an isogeny of $p$-th
power order. The argument follows closely \cite[Thm. 4.4.3]{BS}. Then,
in Section~\ref{A2}, we prove that $\varphi^*$ is an isomorphism: here,
the techniques
used are specific to characteristic $p$.

\subsection{$\ell$-adic realizations}\label{independenceatl}
Let $p= {\rm char} (k)$ and $k = \bar k$. Recall
\cite[10.1.10]{Deligne} for the definition of  $\ell$-adic
realization $\TT_{\ell}(\M)$ of a 1-motive $\M$ over $k$  for
$\ell\neq p$; this is given by taking the inverse limit over $\nu$
of the inverse system $\M[\ell^{\nu}]$ as defined in our section
\ref{Mpinfty}. From (\ref{exactseq}) it yields an exact sequence
(\cf  \cite[\S 1.3 ]{BS})
\begin{equation}\label{etexseq}
0\llongrightarrow  \prod_{\ell\neq p}\TT_{\ell}(\G) \llongrightarrow
\prod_{\ell\neq p}\TT_{\ell}(\M)
\llongrightarrow  \prod_{\ell\neq p}\TT_{\ell}(\X) \llongrightarrow 0
\end{equation}
for  $\M=[\X\by{u} \G]$. Note that for a complex of abelian groups $\cC
= [\cF\to\cG]$ such that ${}_{\ell^{\nu}}\cF=0$  and $\cG/\ell^{\nu}=0$
(for all $\nu>0$) $\TT_{\ell}(\cC)$ can be defined along with an
exact sequence
$$0\to \liminv{\nu}{{}_{\ell^{\nu}}\cG}\to \TT_{\ell}(\cC) \to
\liminv{\nu}{\cF/\ell^{\nu}}\to \liminv{\nu}^1{{}_{\ell^{\nu}}\cG}$$
explaining (\ref{etexseq}).
The following is a modification of \cite[Prop. 1.3.1]{BS}.
\begin{lemma}\label{faith} Let $p= {\rm char} (k)$ and $k = \bar k$.
The functor $\prod_{\ell\neq p}\TT_{\ell}$ from the category of
1-motives to abelian groups is faithful. Moreover, if $\M\to \M'$
is a map of 1-motives such that $$\prod_{\ell\neq
p}\TT_{\ell}(\M)\cong  \prod_{\ell\neq p}\TT_{\ell}(\M')$$ then
$\M\to \M'$ is an isogeny of $p$-power order.\end{lemma}
\begin{proof}
Consider $\M=[\X\by{u} \G]$, $\M'=[\X'\by{u'}\G']$ and
$f\colon\M\to \M'$. By making use of the exact sequence
(\ref{etexseq}) we can see that it is enough to check faithfulness
separately for maps of semi-abelian schemes or lattices. Since
torsion points coprime to $p$ are Zariski dense in a semi-abelian
scheme over $k=\bar k$, $ \prod_{\ell\neq p}\TT_{\ell}(f)=0$
implies $f=0$ for morphisms $f$ between semi-abelian schemes.
Moreover, $ \prod_{\ell\neq p}\TT_{\ell}(\X[1])=\prod_{\ell\neq p}
\X\otimes
\ZZ_{\ell}$ is clearly faithful.\\

If $\M\to \M'$ induces an isomorphism
$\prod_{\ell\neq p}\TT_{\ell}(\M)\cong \prod_{\ell\neq
p}\TT_{\ell}(\M')$ then by (\ref{etexseq}) we have that
$\prod_{\ell\neq p}\TT_{\ell}(\G)$ injects into $\prod_{\ell\neq
p}\TT_{\ell}(\G')$ and $\prod_{\ell\neq p}\TT_{\ell}(\X[1])$ surjects
onto
$\prod_{\ell\neq p}\TT_{\ell}(\X'[1])$, therefore we have an exact
sequence
$$0\to \X''\to \X\to \X'\to \cA\to 0$$
where $\cA$ is finite and killed by a power of $p$.
Moreover by the snake lemma applied to the resulting diagram given by
(\ref{etexseq}) we get that
$$\prod_{\ell\neq p}\TT_{\ell}(\X''[1])\cong \prod_{\ell\neq p}
\frac{\TT_{\ell}(\G')}{\TT_{\ell}(\G)}.$$
      Since $\prod_{\ell\neq p}\TT_{\ell}(\G)\into
\prod_{\ell\neq p}\TT_{\ell}(\G')$ we have that $\cB = \Ker(\G\to \G')$
is a finite group.
Let $\cB\{p\}$ be the sub-group of $p^{n}$-torsion elements for some
$n\gg 0$ then
$$\cB/\cB\{p\}\cong \prod_{\ell\neq
p}\frac{\TT_{\ell}(\G/\cB)}{\TT_{\ell}(\G)}\into \prod_{\ell\neq
p}\frac{\TT_{\ell}(\G')}{\TT_{\ell}(\G)}.$$
Thus $p^{n}\cB=0$, for some $n\gg 0$, since $\cB/\cB\{p\}$ injects into
$\prod_{\ell\neq p}\TT_{\ell}(\X''[1])$ which is torsion free. If we
let $\G''$ denote the cokernel of the map $\G\to \G'$, we then get the
following
exact sequence of complexes
$$0\to [\X''\to \cB]\to [\X\to \G]\to [\X'\to \G']\to [\cA\to \G'']\to
0.$$
Applying $\prod_{\ell\neq p}\TT_{\ell}$ we have
$$\prod_{\ell\neq p}\TT_{\ell}([\X''\to \cB])\to\prod_{\ell\neq
p}\TT_{\ell}(\M)\longby{\cong}\prod_{\ell\neq
p}\TT_{\ell}(\M')\to\prod_{\ell\neq p}\TT_{\ell}([\cA\to\G''])$$
where $\TT_{\ell}([\X''\to \cB])=\TT_{\ell}(\X''[1])$ and
$\TT_{\ell}([\cA\to\G''])= \TT_{\ell}([\G''])$ since
${}_{\ell^{\nu}}\cB=0$  and $\cA/\ell^{\nu}=0$ (for all $\nu>0$
and $\ell\neq p$) respectively. Therefore the composition of the
induced maps here-above is the zero  map as well as an
isomorphism. Thus $\prod_{\ell\neq  p}\TT_{\ell}(\X''[1])=
\prod_{\ell\neq p}\TT_{\ell}(\G'')=0$ whence $\X''= \G''=0$, \ie
$\G\to \G'$ is an isogeny with kernel the finite group $\cB$, $\X
\to \X'$ is injective with cokernel the finite group $\cA$ and we
can find  a positive integer $\bar\nu$ such that
$p^{\bar\nu}\cA=p^{\bar\nu}\cB=0$.\end{proof}

Now just apply Lemma~\ref{faith} to our map $\varphi^*$ and obtain
an induced isomorphism
$$\prod_{\ell\neq p}\TT_{\ell}(\varphi^*)\colon \prod_{\ell\neq
p}\TT_{\ell}(\Pic^{+}(X_{\d}, Y_{\d}))\longby{\simeq} \prod_{\ell\neq
p}\TT_{\ell}(\Pic^{+}(X_{\d}', Y_{\d}'))$$
by cohomological descent. In fact, the same arguments in the proof of
Theorem~4.4.3 in \cite{BS} applies here (see also \cite[7.2]{BS} for
compatibility with Galois actions) provided that $\ell \neq p$, and
therefore
$$\TT_{\ell}(\Pic^{+}(X_{\d}', Y_{\d}')) = \H^1_{\et} (V, \ZZ_{\ell}(1))
= \TT_{\ell}(\Pic^{+}(X_{\d}, Y_{\d})).$$ We then have that
$\varphi^*$ is an isogeny of $p$-power order.

\begin{corollary} If $\M\to \M'$
is a map of 1-motives such that $$\prod_{\ell\neq
p}\TT_{\ell}(\M)\cong  \prod_{\ell\neq p}\TT_{\ell}(\M')
\qquad\hbox{{\rm and}}\qquad \TT_\crys(\M) \cong
\TT_\crys(\M'),$$then $\M\to\M'$ is an isomorphism.
\end{corollary}
\begin{proof} By Lemma~\ref{faith} we know that $\M \to\M'$ is a
$p$-power isogeny. In particular, $\G \to \G'$ is a $p$-power
isogeny and the map $\X\to\X'$ is injective with $p$-power
cokernel. By definition the crystalline realization
$\TT_\crys(\M)$ of $\M$ is the covariant Dieudonn\'e module of the
$p$-divisible group $\M[p^\infty]$. Since the Dieudonn\'e functor
is fully faithful, we deduce that $\M[p^\infty]\to \M'[p^\infty]$
is an isomorphism. Note that we have an exact sequence $0\to
\G[p^\infty] \to\M[p^\infty] \to \X[p^\infty] \to 0$.  Since the
map $\X[p^\infty]=\X\otimes \QQ_p/\ZZ_p \to
\X'[p^\infty]=\X'\otimes \QQ_p/\ZZ_p$ is injective, we deduce from
the above exact sequence that it is also surjective and that
$\G[p^\infty]\to \G'[p^\infty]$ is an isomorphism. The conclusion
follows.
\end{proof}

\subsection{$p$-adic realization}\label{A2}
Let $\M:=[\X \rightarrow \G]:=\Pic^{+}(X_{\d}, Y_{\d})$ and $\M':=[\X'
\rightarrow \G']:=\Pic^{+}(X_{\d}, Y_{\d})$. Let $\G$ be an
extension of the abelian variety $\A$ by the torus $\T$. Let $\G'$
be an extension of the abelian variety $\A'$ by the torus $\T'$.
All cohomology groups in the sequel are considered for the
fppf-topology.

\begin{lemma}\label{Milnegeneralized} Let $\pi
\colon X_{\d}\to \Spec(k)$ be the structural morphism. Let $\cF$ be
a finite commutative group scheme over $k$. Let
$\cF^\vee:=\HOM(\cF,\GG_m)$ be the Cartier dual of $\cF$. The natural
map of $\fppf$-sheaves over $k$
$$R^1\pi_{*}(\cF_{X_{\d}}) \longrightarrow \HOM
(\cF^\vee,\bPic_{X_{\d}/k}),$$defined by push-forward of simplicial
$\cF$-torsors via elements of $\cF^\vee$, is an
isomorphism.\end{lemma}

\begin{proof} By \cite[Prop. III.4.16]{Milne} the maps
$R^1(\pi_{i})_*(\cF_{X_{i}})\rightarrow \HOM (\cF^\vee,\Pic_{X_i/k})$
are
isomorphisms for every $i$. Thus, the map
$$\Ker\left(R^1(\pi_{0})_*(\cF_{X_{0}}) \rightarrow
R^1(\pi_{1})_*(\cF_{X_{1}})\right) \longrightarrow
\HOM\left(\cF^\vee,\Ker\bigl(\Pic_{X_0/k}\rightarrow
\Pic_{X_1/k}\bigr)\right)$$is an isomorphism. Furthermore,
$(\pi_{i})_*(\cF_{X_{i}})=\HOM(\cF^\vee,(\pi_{i})_*(\GG_m))$.
By \cite[Lemma III.4.17]{Milne} we have $\EXT^1(\cF^\vee,\GG_m)=0$
for every $i$. Then, $\HOM(\cF^\vee, \_)$ preserves exact sequences
of tori. Using the spectral sequences for
$R^1\pi_{*}(\cF_{X_{\d}})$ and for $R^1\pi_{*}(\GG_{m, {X_{\d}}})$
the lemma follows.\end{proof}

\begin{corollary}\label{consequenceMilnegeneralized} We have
$\HH^1(X_{\d},\ZZ_p)\cong\Hom(\mu_{p^\infty},\G)$.\end{corollary}
\begin{proof}  By Lemma~\ref{Milnegeneralized}, the map
$\HH^1(X_{\d},\ZZ_p)\to \Hom(\mu_{p^\infty},\bPic_{X_{\d}/k})$ is an
isomorphism. If $\cF$ is a finite $k$-group scheme or a discrete
group scheme, we have $\Hom(\mu_{p^\infty},\cF)=0$. Let
$K:=\Ker\bigl(\Pic_{X_0/k}\to \Pic_{X_1/k}\bigr)$ and let
$\cF:=K/K^{0,\red}$. Thus, $\Hom(\mu_{p^\infty},K^{0,\red}) \to
\Hom(\mu_{p^\infty},K)$ is an isomorphism. Hence, the map
$\Hom(\mu_{p^\infty},\G)\hookrightarrow
\Hom(\mu_{p^\infty},\bPic_{X_{\d}/k})$ is an isomorphism, as claimed.
\end{proof}

\begin{lemma}\label{ZZpinj}
Let $X$ be a normal $k$--scheme and let $U \subset X$
be an open dense subscheme. The map $H^1(X,\ZZ/p^n\ZZ) \rightarrow
H^1(U,\ZZ/p^n\ZZ)$ is injective.\end{lemma}
\begin{proof} Let $Y \to X$
be a $\ZZ/p^n\ZZ$--torsor. Since the map $Y\to X$ is \'etale and $X$
is normal, $Y$ is normal. Suppose that $Y\vert_U\rightarrow U$ is
trivial, \ie $Y\vert_U=\amalg_{i=1}^{p^n} U$. In particular, the
normalization~$Y':=\amalg_{i=1}^{p^n} X$ of~$X$ in~$Y\vert_U$ is the
trivial $\ZZ/p^n\ZZ$-torsor. The normality of~$Y$ implies that
$Y\cong Y'$ as $\ZZ/p^n\ZZ$--torsors over~$X$. Thus, $Y$ is the
trivial torsor as claimed.
\end{proof}

\begin{lemma}\label{ZZpntorsion} The group
$$\lim_{\infty \leftarrow n} \left(\H^1(V_i,\ZZ/p^n
\ZZ)/\H^1(X_i,\ZZ/p^n \ZZ)\right)$$is torsion free as
$\ZZ_p$-module. The map
$$\alpha_i\colon \lim_{\infty\leftarrow n}\left(\H^1(V_i,\ZZ/p^n \ZZ)/
\H^1(X_i,\ZZ/p^n \ZZ)\right)\longrightarrow \lim_{\infty\leftarrow
n}\left( \H^1(V_i',\ZZ/p^n \ZZ)/\H^1(X_i',\ZZ/p^n \ZZ)\right)$$is
injective.\end{lemma}

\begin{proof}  If $Z\subset X_i$ is an
irreducible divisor, we let $R_Z$ be a complete dvr which is an
extension of the completed local ring $\hat{\cO}_{X_i,Z}$ of $X_i$ at the
generic point of $Z$ with residue field equal to an
algebraic closure $k_Z$ of the fraction field $k(Z)$ of $Z$ and with maximal
ideal generated by that of $\hat{\cO}_{X_i,Z}$. Let
$K_Z$ be the fraction field of $R_Z$ and let $m_Z$ be the maximal
ideal of $R_Z$. For every $r\in\NN$ let $U_{Z,(r)}$ be the group
of units of $R_Z$ congruent to $1$ modulo $m_Z^r$. Let
$W_{Z,(r)}:=U_{Z,(1)}/U_{Z,(r)}$ and $W_Z:=\lim_{\infty\leftarrow
r}W_{Z,(r)}$. For every integer $1 \leq i\leq r-1$ prime to $p$
let $r_i$ be the smallest integer such that $p^{r_i}\geq r/i$. By
\cite[\S V.9 Prop.~9]{SeGrAlg} the group $W_{Z,(r)}$ is isomorphic
to the product of truncated Witt vectors $\WW_{r_i}(k_Z)$ over the
integers $1\leq i\leq r-1$ prime to $p$. By purity of the branch
locus \cite[X.3.4]{SGA2}  a $\ZZ/p^n\ZZ$-torsor over $V_i$ which
is unramified at the generic points of $Y_i$ is unramified. Thus,
$\H^1(V_i,\ZZ/p^n \ZZ)/\H^1(X_i,\ZZ/p^n \ZZ)$ is contained in
$\oplus_j
\H^1(K_{Y_{ij}},\ZZ/p^n\ZZ)/\H^1(R_{Y_{ij}},\ZZ/p^n\ZZ)$. By class
field theory, \cf \cite[\S XV.2]{CoLoc}, the latter is isomorphic
to $ \oplus_j \Hom(W_{Y_{ij}},\ZZ/p^n\ZZ)$. In particular, the map
$$\Psi_{X_i}\colon \lim_{\infty\leftarrow n}\left(\H^1(V_i,\ZZ/p^n \ZZ)/
\H^1(X_i,\ZZ/p^n \ZZ)\right)\longrightarrow \oplus_j
\Hom(W_{Y_{ij}},\ZZ_p)$$is injective. Let $Y_{ih}'$ be an
irreducible component of $Y_i'$ dominating $Y_{ij}$. Let ${\bf
Norm}_{jh}$ be the homomorphism $W_{Y_{ih}'}\rightarrow
W_{Y_{ij}}$ defined by the norm map from $K_{Y_{ih}'}$ to
$K_{Y_{ij}}$. Then, the map from $\H^1(V_i,\ZZ/p^n
\ZZ)/\H^1(X_i,\ZZ/p^n \ZZ)$ to $\H^1(V_i',\ZZ/p^n
\ZZ)/\H^1(X_i',\ZZ/p^n \ZZ)$ is compatible with the map obtained
applying $\Hom(\_, \ZZ/p^n\ZZ)$ to $\prod_{j,h}{\bf Norm}_{jh}$.
By class field theory the quotient of $W_{Y_{ih}'}\rightarrow
W_{Y_{ij}}$ is isomorphic to the wild inertia of the abelianized
Galois group of the extension $K_{{Y_{ih}'}}/K_{Y_{ij}}$; see
\cite[Cor. XV.2.3]{CoLoc}. In particular, it is a finite
$p$-group. Applying $\Hom(\_,\ZZ_p)$ we conclude that the map
$\alpha_i$ is injective as claimed.\end{proof}

\begin{corollary}\label{mup} The map
$\Hom(\mu_{p^\infty},\G)\rightarrow \Hom(\mu_{p^\infty},\G')$ is
an isomorphism.\end{corollary}

\begin{proof} Since $V_{\d}$ and $V_{\d}'$ are hypercoverings of $V$, we
have cohomological descent for the \'etale cohomology of constant sheaves. 
In particular, $\HH^i(V_{\d},\ZZ/p^n\ZZ)\cong
\H^i(V,\ZZ/p^n\ZZ)\cong \HH^i(V'_{\d},\ZZ/p^n\ZZ)$ for every $n$.
Hence, the induced map $\HH^i(V_{\d},\ZZ/p^n\ZZ)\rightarrow
\HH^i(V'_{\d},\ZZ/p^n\ZZ)$ is an isomorphism. Taking inverse
limits over $n\in\NN$ we have a commutative diagram with
exact rows (using  Lemma~\ref{ZZpinj})

\begin{displaymath}
\begin{CD}
0  @>>>  \HH^1(X_{\d},\ZZ_p) @>>> \HH^1(V_{\d},\ZZ_p) @>>>
\H^1(V_0,\ZZ_p)/\H^1(X_0,\ZZ_p) \\ 
& &  @VVV@V{\wr}VV @V{\alpha_0}VV\\
0  @>>>  \HH^1(X_{\d}',\ZZ_p) @>>> \HH^1(V_{\d}',\ZZ_p) @>>>
\H^1(V_0',\ZZ_p)/\H^1(X_0',\ZZ_p) \\
\end{CD}
\end{displaymath}
By Lemma~\ref{ZZpntorsion} the map $\alpha_0$ is injective. Thus,
the left vertical map is an isomorphism. By
Corollary~\ref{consequenceMilnegeneralized} such map coincides with the
map $\Hom(\mu_{p^\infty},\G)\rightarrow \Hom(\mu_{p^\infty},\G')$.
This proves the claim.\end{proof}

\begin{lemma}\label{alphapZZpZZ} The maps
$\Hom(\alpha_p, \Pic_{X_0/k}) \rightarrow
\Hom(\alpha_p,\Pic_{X_0'/k})$ and  $\Hom(\ZZ/p \ZZ,\Pic_{X_0/k})
\rightarrow \Hom(\ZZ/ p \ZZ,\Pic_{X_0'/k})$ are injective.
Thus $\Hom(\alpha_p,\G) \rightarrow
\Hom(\alpha_p,\G')$ and $\Hom(\ZZ/p \ZZ,\G) \rightarrow \Hom(\ZZ/
p \ZZ,\G')$ are injective.\end{lemma}

\begin{proof}  For any finite
commutative $k$-group scheme $\cF$, let $\cF^\vee$ be the Cartier dual
of $\cF$. By \cite[Prop. III.4.16]{Milne} we have
$\Hom(\cF,\Pic_{X_0/k})\cong \H^1(X_0,\cF^\vee)$ (and analogously for
$X_0'$). For $\cF=\alpha_p$ or $\cF=\ZZ/p\ZZ$, the latter group is
identified with a subgroup of $\H^0(X_0,\Omega^1_{X_0/k})$ (resp.
$\H^0(X_0',\Omega^1_{X_0'/k})$); \cf \cite[Prop. III.4.14]{Milne}.
Note that $\H^0(X_0,\Omega^1_{X_0/k}) \subset \Omega^1_{k(X_0)/k}$
and that $\H^0(X_0',\Omega^1_{X_0'/k}) \subset
\Omega^1_{k(X_0')/k}$. Furthermore, the map $X_0'\rightarrow X_0$
is generically separable by construction; thus,
$\Omega^1_{k(X_0)/k} \subset \Omega^1_{k(X_0')/k}$. This proves
the first claim. Since any homomorphism from $\alpha_p$ or $\ZZ/
p\ZZ$ to a torus is trivial, it suffices to prove the second
assertion for $\A$ in place of $G$ and $\A'$ in place of $\G'$. By
construction $\A$ is the reduced kernel of the map $\Pic^0_{X_0/k}
\rightarrow \Pic^0_{X_1/k}$ and $\A'$ is the reduced kernel of the
map $\Pic^0_{X_0'/k} \rightarrow \Pic^0_{X_1'/k}$. Thus, it suffices
to prove the Lemma for $\Pic^{0,\red}_{X_0/k}$ instead of $\G$ and
$\Pic^{0,\red}_{X_0'/k}$ instead of $\G'$. The conclusion
follows.\end{proof}

\begin{corollary}\label{independencesemiableian} The homomorphism
$\G \rightarrow \G'$ is an isomorphism.
\end{corollary}
\begin{proof} By Section~\ref{independenceatl} for every $\ell\neq p$
the
homomorphism of Tate modules ${\bf T}_\ell(\G)\rightarrow {\bf
T}_\ell (\G')$ is an isomorphism. This implies that the map
$\G\rightarrow \G'$ is an isogeny with kernel $\cF$ of $p$-power
order. By Corollary~\ref{mup} and Lemma~\ref{alphapZZpZZ}, $\cF$ does
not
contain any subgroup isomorphic to $\mu_p$ or $\alpha_p$ or
$\ZZ/p\ZZ$. Hence, $\cF=0$.\end{proof}

\begin{proposition}\label{indeP1Mot} The homomorphism
$\varphi^*\colon \M \rightarrow \M'$ is an isomorphism.
\end{proposition}
\begin{proof} By Section~\ref{independenceatl} the map $\varphi^*$
is an isogeny of $p$-power order. In particular, the lattices $\X$
and $\X'$ have the same rank and $\X'/\X$ is killed by a power of
$p$. By Corollary~\ref{independencesemiableian} the  map induced by
$\varphi^*$ on the semiablian parts is an isomorphism. It then
suffices to prove that the map $\X/p\X \to \X'/p\X'$ is injective.
Let $\cC \df [\bDiv_{Y_{\d}}(X_{\d})\to  \bPic (X_{\d})]$ as a complex
of abelian groups.
As in \ref{independenceatl}, denote
$$\TT_{\ZZ/p\ZZ}(\cC):=\frac{\left\{(D,\ccL_{\d})\in
\bDiv_{Y_{\d}}(X_{\d})\times \bPic (X_{\d}) \vert \eta_{\d}\colon
\ccL_{\d}^p\cong
\cO_{X_{\d}}(-D)\right\}}{\left\{\bigl(pD,\cO_{X_{\d}}(-
D)\bigr)\right\}}.
$$It sits in an exact sequence $$0
\longrightarrow \TT_{\ZZ/p\ZZ}(\bPic(X_{\d}))\longrightarrow
\TT_{\ZZ/p\ZZ}(\cC) \longrightarrow
\bDiv_{Y_{\d}}(X_{\d})/p\bDiv_{Y_{\d}}(X_{\d}) \longrightarrow
0.$$ Let $\rho_p\colon \TT_{\ZZ/p\ZZ}(\cC) \to
\HH^1(V_{\d},\mu_p)$ be the map
$\rho_p\bigl((D,\ccL_{\d})\bigr):=(\ccL_{\d},\eta_{\d})\vert_{V_{\d}}$
defined via simplicial Kummer theory, \cf \cite[\S 4.4]{BS}. As in
loc.~cit. the induced map $\TT_{\ZZ/p\ZZ}(\bPic(X_{\d}))\to
\HH^1(X_{\d},\mu_p)$ is an isomorphism. Since $\mu_p(X_i)=0$, the map
$\HH^1(X_{\d},\mu_p)\to \H^1(X_0,\mu_p)$ is injective. Note that
$\rho_p\bigl((D,\ccL_{\d})\bigr)$ is a simplicial $\mu_p$-covering
of $X_{\d}$ and the associated covering of $X_0$ is ramified
exactly over the support of $D$. Thus, the inverse image of
$\HH^1(X_{\d},\mu_p)$ via $\rho_p$ is
$\TT_{\ZZ/p\ZZ}(\bPic(X_{\d}))$ and the induced map
$\bDiv_{Y_{\d}}(X_{\d})/p\bDiv_{Y_{\d}}(X_{\d}) \to
\H^1(V_0,\mu_p)/\H^1(X_0,\mu_p)$ is injective. Hence, in the
following commutative diagram
\begin{displaymath}
\begin{CD}
0@>>>\TT_{\ZZ/p\ZZ}(\bPic(X_{\d})) @>>>
\TT_{\ZZ/p\ZZ}(\cC)@>>>
\bDiv_{Y_{\d}}(X_{\d})/p\bDiv_{Y_{\d}}(X_{\d})@>>>0\\
&&@VVV @VVV@VVV\\
0@>>> \H^1(X_0,\mu_p) @>>> \H^1(V_0,\mu_p)@>>>
\H^1(V_0,\mu_p)/H^1(X_0,\mu_p)@>>>0\\
\end{CD}
\end{displaymath}
the left and right vertical arrows are injective. Therefore, we
conclude that
\begin{equation}\label{exactMwithV0}
\TT_{\ZZ/p\ZZ}(\cC) \into \H^1(V_0,\mu_p).
\end{equation}
is an injective map.\\

Note that we assumed that $V_0 \rightarrow V_0'$ is a dominant map of
normal schemes over $k$ which is generically finite and separable, see
\cite[Thm. 4.1]{DE}. Let $k(V_0)$ and $k(V_0')$ be the function fields
of $V_0$ and $V_0'$ respectively.
Consider the diagram
\begin{displaymath}
\begin{CD}
\H^1(V_0,\mu_p) @>>> \H^1(k(V_0),\mu_p)@>{\sim}>>
k(V_0)^*/\bigl(k(V_0)^*)^p \\
@VVV @VVV @VVV\\
\H^1(V_0',\mu_p) @>>> \H^1(k(V_0'),\mu_p)@>{\sim}>>
k(V_0')^*/\bigl(k(V_0')^*)^p\\
\end{CD}
\end{displaymath}
By normality the horizontal arrows are injective. By the
separability assumption the right vertical arrow is injective.
Thus, the induced map $\H^1(V_0,\mu_p)\to \H^1(V_0',\mu_p)$ is
injective.
Let $\cC' \df [\bDiv_{Y_{\d}'}(X_{\d}')\to  \bPic (X_{\d}')]$.

Using (\ref{exactMwithV0}) for $\cC$ and
$\cC'$, we obtain that
\begin{equation}\label{injectpicv}
\TT_{\ZZ/p\ZZ}(\cC) \into
\TT_{\ZZ/p\ZZ}(\cC')
\end{equation} is injective.\\

Let $\TT_{\ZZ/p\ZZ}(\M)$ be (as in Section~\ref{Mpinfty}) the  group of
$k$-valued $p$-torsion points of $\M$: it sits in the exact sequence
\eqref{exactseq}.
By construction it maps
to $\TT_{\ZZ/p\ZZ}(\cC)$. Let ${\rm NS}(X_{\d})\df
\pi_0(\bPic_{X_{\d}/k})$. It is a discrete group and it
coincides with $\bPic_{X_{\d}/k}(k)/ \bPic^{0,\red}_{X_{\d}/k}(k)$.
The map $\TT_{\ZZ/p\ZZ}(\G) \to \TT_{\ZZ/p\ZZ}(\bPic (X_{\d}))$ is
injective. Since $\X$ is the fiber product of
$\bDiv_{Y_{\d}}(X_{\d}) \to \bPic(X_{\d})$ and
$\bPic^{0,\red}(X_{\d})\to  \bPic(X_{\d}) $, the kernel
$\cF$ of $\X/p\X \to \bDiv_{Y_{\d}}(X_{\d})/ p \bDiv_{Y_{\d}}(X_{\d})$
is a
subgroup of ${\rm NS}(X_{\d})$ (as well as $\cF'$ injects in ${\rm NS}
(X_{\d}')$). By Corollary~\ref{independencesemiableian}
the map $\bPic^{0,{\rm red}}(X_{\d}) \to \bPic^{0,{\rm
red}}(X_{\d}')$ is an isomorphism. The
$p$-power torsion of $\bPic(X_{\d})$ injects into $\Pic(X_{0})$, as
$\mu_p(X_i)=0$ (as well as for $X_{\d}'$). Then,  by
Lemma~\ref{alphapZZpZZ}, the
map $\bPic(X_{\d})\to \bPic(X_{\d}')$ is injective on the
$p$-power torsion. Thus, the map ${\rm NS}(X_{\d})\to {\rm
NS}(X'_{\d})$ is injective on $p$-power torsion: in particular, $\cF$
injects into $\cF'$. Therefore, the kernel of
$\TT_{\ZZ/p\ZZ}(\M) \to \TT_{\ZZ/p\ZZ}(\cC)$ is contained in
the kernel of
$\TT_{\ZZ/p\ZZ}(\M') \to \TT_{\ZZ/p\ZZ}(\cC')$. By
\eqref{injectpicv} we then get that the map
$\TT_{\ZZ/p\ZZ}(\M) \to \TT_{\ZZ/p\ZZ}(\M')$ is injective. Since by
Corollary~\ref{independencesemiableian}
the map $\TT_{\ZZ/p\ZZ}(\G) \to \TT_{\ZZ/p\ZZ}(\G')$ is an
isomorphism, we conclude that the map $\X/p\X\to \X'/p\X'$ is
injective as claimed.
\end{proof}

\vfill

\noindent
F. Andreatta\\
Dipartimento di Matematica Pura ed Applicata\\
Universit\`a degli Studi di Padova\\
Via G. Belzoni, 7\\
35131 - Padova\\
E-Mail: {\tt fandreat@math.unipd.it}\\[1cm]
L. Barbieri Viale\\
Dipartimento di Metodi e Modelli Matematici\\
Universit\`a degli Studi di Roma {\it La Sapienza}\\
Via A. Scarpa, 16\\
00161 - Roma\\
E-Mail: {\tt Luca.Barbieri-Viale@uniroma1.it}

\vfill


\begin{thebibliography}{}

\bibitem{BM}{\sc L. Barbieri-Viale, A. Rosenschon {\rm and} M. Saito}:
Deligne's conjecture on 1-motives, {\it Annals of Math.}  {\bf 158} N. 2
(2003) 593-633.

\bibitem{BS}{\sc L. Barbieri-Viale {\rm and} V. Srinivas}:
Albanese and Picard 1-motives, M\'emoires de la Soci\'et\'e
Math\'ematique de France {\bf 87} Paris, 2001.

\bibitem{BerMess}{\sc P. Berthelot {\rm and} W. Messing}: Th\'eorie de
Dieudonn\'e cristalline: I, in {\it Ast\'erisque} {\bf 63}
Soci\'et\'e Math\'ematique de France, Paris (1979) 17-37.

\bibitem{BerthelotOgus}{\sc P. Berthelot {\rm and} A. Ogus}: Notes
on crystalline cohomology, {\it Mathematical Notes Princeton
University Press} {\bf 21}, Princeton, 1978.

\bibitem{BL}{\sc S. Bosch, W. Lutkebohmert {\rm and} M. Raynaud}:
N\'eron Models, Springer {\it Ergebnisse der Math.} {\bf 21}
Heidelberg, 1990.

\bibitem{DE}{\sc A.J. de Jong}: Smoothness, semistability and
alterations, {\it Publ. IHES} {\bf 83} (1996) 51-93.

\bibitem{Deligne}{\sc P. Deligne}: Th\'eorie de Hodge III {\it Publ.
Math.}\, IHES {\bf 44} (1974) 5-78.

\bibitem{FO} {\sc J.-M. Fontaine}: Groupes $p$-divisibles sur les corps
locaux, {\it Ast\'erisque} {\bf 47-48} Soci\'et\'e Math\'ematique de
France, Paris, 1977.

\bibitem{FJ} {\sc J.-M. Fontaine {\rm and} K. Joshi}: Notes on
1-Motives, work in progress (Preliminary Version: May, 1996).

\bibitem{Illusie}{\sc L. Illusie}: Complexe de De Rham-Witt et
cohomologie cristalline, {\it Ann. scient. \'Ec. Norm. Sup.} $4^e$
s\'erie t. {\bf 12} (1979) 501-661.

\bibitem{MG}{\sc A. Grothendieck}: Groupes de Barsotti-Tate et
cristaux de Dieudonn\'e, {\it S\'eminaire de math\'ematiques
sup\'ereures} {\bf 45} Montreal, 1974.

\bibitem{SGA2}{\sc A. Grothendieck {\rm et al.}}: SGA2 - Cohomologie
locale des
faisceaux coh\'erents et th\'eor\`emes de Lefschetz locaux et globaux
(1962),
North-Holland, Amsterdam, 1968.

\bibitem{SGA4}{\sc A. Grothendieck {\rm et al.}}: SGA4 -  Th\'eorie des
topos
et cohomologie \'etale des sch\'emas (1963-64), Springer LNM {\bf 269
270
305}, Heidelberg, 1972-73.

\bibitem{SGA7}{\sc A. Grothendieck {\rm et al.}}: SGA7 - Groupes de
monodromie
en g\'eom\'etrie alg\'ebrique (1967-68), Springer LNM {\bf 288 340},
1972-73.

\bibitem{Kato} {\sc K. Kato}: Logarithmic structures of
Fontaine-Illusie, in {\it Proceedings of JAMI conference}, (1988),
191-224.

\bibitem{MazurMessing} {\sc B. Mazur {\rm and} W. Messing}: Universal
extensions and one dimensional crystalline cohomology,  Springer LNM
{\bf 370} Berlin-Heidelberg-New York, 1974.

\bibitem{Messing} {\sc W. Messing}: The crystals associated to
Barsotti-Tate groups with applications to abelian schemes, Springer LNM
{\bf 264} Berlin-Heidelberg-New York, 1972.

\bibitem{Milne}{\sc J. S. Milne}:  \'Etale cohomology, {\it Princeton
University
Press}, Princeton, New Jersey, 1980.

\bibitem{Oda} {\sc T. Oda}: The first De Rham cohomology group and
Dieudonn\'e modules, {\it Ann. scient. \'Ec. Norm. Sup.} $4^e$
s\'erie t. {\bf 2} (1969) 63-135.

\bibitem{OG}{\sc F. Oort}: Commutative group schemes, Springer
LNM {\bf 15}, 1966.

\bibitem{Pet} {\sc D. Petrequin}: Classes de Chern et classes de cycles
en cohomologie rigide, (Ph D Thesis) {\it Bull. Soc. Math. France} {\bf
     131}  (2003),  no. 1, 59--121.

\bibitem{RA}{\sc N. Ramachandran}: Duality of Albanese and Picard
1-motives, {\it K-Theory} {\bf 22} (2001) 271-301.

\bibitem{Ray}{\sc M. Raynaud}: $1$-motifs et monodromie
g\'eom\'etrique, Expos\'e VII, {\it Ast\'erisque} {\bf 223} Soci\'et\'e
Math\'ematique de France, Paris (1994) 295-319.

\bibitem{Shi}{\sc A. Shiho}: Crystalline fundamental groups. II. Log
convergent cohomology and rigid  cohomology,  {\it J. Math. Sci. Univ.
Tokyo} {\bf  9} (2002),  no. 1, 1--163.

\bibitem{SeGrAlg}{\sc J.-P. Serre}: Groupes alg\'ebriques et corps de
classes, (Publications de l'Universit\'e de Nancago, No. VII) {\it
Hermann}, Paris, 1959.

\bibitem{CoLoc}{\sc J.-P. Serre}: Corps locaux, (Deuxi\`eme \'edition,
Publications de l'Universit\'e de Nancago, No. VIII) {\it Hermann},
Paris,  1968.

\bibitem{Tsu}{\sc N. Tsuzuki}: Cohomological descent of rigid
cohomology for proper coverings,  {\it Invent. Math.} {\bf  151}
(2003),  no. 1, 101--133.

\end{thebibliography}
\end{document}